\documentclass[11pt]{article}
\usepackage{fullpage}
\usepackage{natbib}
\usepackage{graphicx} 

\usepackage[utf8]{inputenc} 
\usepackage[T1]{fontenc}    
\usepackage[colorlinks,allcolors=blue,backref=page]{hyperref}     
\renewcommand*{\backref}[1]{}
\renewcommand*{\backrefalt}[4]{%
    \ifcase #1 (Not cited.)%
    \or        (Cited on page~#2.)%
    \else      (Cited on pages~#2.)%
    \fi}
\usepackage{url}            
\usepackage{booktabs}       
\usepackage{amsfonts}       
\usepackage{nicefrac}       
\usepackage{microtype}      
\usepackage{xcolor}         

\usepackage{amsthm}
\usepackage{amsmath}
\usepackage{cleveref}

\usepackage{amsfonts}
\usepackage{amssymb}
\usepackage{xpatch}
\usepackage{enumitem}

\usepackage{algorithm}
\usepackage{algpseudocode}

\usepackage{graphicx}
\usepackage{thmtools,thm-restate}
\usepackage{pifont}%

\newcommand{\ip}[2]{\left\langle #1,#2\right\rangle}
\newcommand{\norm}[2][]{\left\|#2\right\|_{#1}}
\newcommand{\eps}{\varepsilon}
\usepackage{multirow}
\usepackage{wrapfig}
\usepackage{tikz}
\tikzstyle{vertex}=[circle, draw, inner sep=0pt, minimum size=6pt]
\usetikzlibrary{arrows}

\usepackage{xcolor}
\usepackage{colortbl} 
\usepackage{tabularx} 
\usepackage{makecell}  
\usepackage{tabstackengine} 
\usepackage{adjustbox} 
\usepackage{footnote}
\usepackage{caption}
\makesavenoteenv{tabular}
\usepackage{titletoc} 
\usepackage{comment}
\usepackage{float}

\newtheorem{theorem}{Theorem}
\newtheorem{lemma}{Lemma}

\newtheorem{property}{Property}
\newtheorem{corollary}[theorem]{Corollary}
\newtheorem{proposition}[theorem]{Proposition}
\newtheorem{definition}{Definition}
\newtheorem{assumption}{Assumption}
\newtheorem{problem}{Problem}
\newtheorem{fact}{Fact}

\usepackage{thm-restate}

\DeclareMathOperator\AndKeyWord{and}

   \DeclareMathOperator*{\argmin}{argmin}

\DeclareMathOperator\prox{prox}
\DeclareMathOperator\bprox{bprox}

\newcommand{\R}{\mathbb{R}}

\newcommand{\E}{\mathbb{E}}

\newcommand{\todo}[1]{}

\newcommand{\greg}{$p$-GREG}

\newcommand{\odc}{$\mathrm{ODC}_p$~}
\usepackage{thm-restate}

\usepackage{makecell}
\usepackage{array}
\usepackage{amsmath}
\usepackage{booktabs}

\usepackage{algpseudocode}
\usepackage{chngcntr} 
\usepackage[top=1in, right=1in, left=1in, bottom=1in]{geometry}

\title{Bias-Reduced and Parallel Heavy-Tailed Stochastic\\ Convex Optimization in $\ell_p$ Geometries}

\author{
\begin{tabular}{c@{\hspace{2cm}}c@{\hspace{2cm}}c}
Zhengyun Nie\thanks{Authors listed in alphabetical order.}
& Anant Shyam\footnotemark[1] & Brian Bullins \\
Purdue University & Purdue University & Purdue University \\
\texttt{nie34@purdue.edu} &
\texttt{shyama@purdue.edu} & \texttt{bbullins@purdue.edu}
\end{tabular}
}

\date{}
\begin{document}
\maketitle
\thispagestyle{empty}
\begin{abstract}
    We consider the fundamental problem of parallel stochastic convex optimization (SCO), whose aim is to approximately minimize a convex objective when given the ability to query a moment-bounded stochastic gradient oracle multiple times in parallel per round. Whereas most previous results have focused on the variance-bounded Euclidean ($\ell_2$-norm) setting, meaning the objective is Lipschitz with respect to the $\ell_2$-norm over an $\ell_2$-bounded feasible set in $\R^d$, in this work we study parallel SCO in $\ell_p$ settings for $2 \le p < \infty$ under a heavy-tailed noise model. 
    
    To this end, we begin by introducing a new class of $\ell_p$-optimal-distance-convergence (ODC$_p$) algorithms, drawing on prior stochastic bias-reduced techniques~\cite{asi2021stochastic}, which combine with multilevel Monte Carlo methods to efficiently estimate the optimum of uniformly convex objectives.  As a consequence of our \odc algorithms, we first establish rates for  projection-efficient $\ell_p$-Lipschitz optimization, composite convex minimization, and efficient $\ell_p$-proximal point estimation. 
        
    Further combining our \odc algorithms with a new $\ell_p$-generalized Gaussian-based stochastic gradient estimator (\greg), along with advances in $\ell_p$ ball-oracle acceleration~\cite{adil2026convex,contreras2025non}, we provide two efficient $\ell_p$ parallel algorithms, EpochDA and $\ell_p$-AC-SA, which achieve optimization error $\eps > 0$ with at most $\tilde{O}(d^{\frac{p}{2p + 2}} \eps^{-\frac{p}{p + 1}})$ query depth, improving upon the sequential heavy-tailed stochastic query depth of $O(\eps^{-p})$ for $\eps \leq \tilde{O}(d^{-1/(2p)})$. In addition, we establish dimension-dependent lower bounds of \(\Omega\bigl(\min\{\varepsilon^{-p},d^{p/(3p-2)}\}\bigr)\) adaptive rounds for nonsmooth convex optimization over a subset of the unit \(\ell_p\) ball, with polynomially many parallel local-oracle queries per round. Our results improve upon the \(\Omega(\min\{\varepsilon^{-p},\varepsilon^{2/3}d^{1/3}\})\) bounds of \cite{diakonikolas2019lower}, when adapted to the same feasible set.
\end{abstract}
\clearpage
\newpage
\tableofcontents
\thispagestyle{empty}

\clearpage
\newpage 
\pagestyle{plain} 
\pagenumbering{arabic}
\section{Introduction}
Stochastic methods are an essential part of modern optimization, with applications spanning machine learning, statistics, and theoretical computer science. The foundational work of Nemirovski and Yudin~\cite{nemirovskij1983problem} introduced the \emph{oracle model} for convex optimization, whereby the goal is to understand, for a particular oracle class, how many queries are needed to minimize a function to sufficient accuracy. Often considered in the case of stochastic convex optimization (SCO) is the stochastic first-order oracle, which for a convex function $f : \R^d \to \R$ queried at a point $x$, provides an unbiased estimate of the subgradient with bounded $q$-th moment ($q \in (1,2]$), that is, $\widehat{\nabla} f(x)$ such that $\E[\widehat{\nabla} f(x)] \in \partial f(x)$ and $\E[\|\widehat{\nabla} f(x)\|_q^q] \leq \sigma^q$ for $\sigma > 0$.

An additional line of work, initiated by Nemirovski~\cite{nemirovski1994parallel}, has further considered the setting of \emph{parallel} convex optimization, in which we may query an oracle multiple times in parallel for each round of interaction, and where we denote the \emph{query depth} as the total number of such rounds.
While several works have studied first-order methods for parallel settings in the Euclidean ($\ell_2$-norm) case~\cite{bubeck2019complexity,carmon2023resqueing,duchi2012randomized}, meaning $f$ is Lipschitz continuous with respect to the $\ell_2$-norm over an $\ell_2$-bounded feasible set, parallel SCO has been less explored in non-Euclidean settings that deal with $\ell_p$ geometries for $2 \leq p < \infty$. For the $p=\infty$ case, Nemirovski~\cite{nemirovski1994parallel} established $\Omega((\frac{d}{\ln(2Kd)})^{1/3}\ln(1/\eps))$ lower bounds on the number of rounds required to achieve an $\varepsilon$-accurate solution, where $K$ denotes the number of parallel queries in each round, while Diakonikolas and Guzmán~\cite{diakonikolas2019lower} provided lower bounds of the form \(\Omega(\min\{\varepsilon^{-p},\varepsilon^{2/3}d^{1/3}\})\) for $p \geq 2$.\footnote{Diakonikolas and Guzmán~\cite{diakonikolas2019lower} also established $\Omega(\eps^{-2})$ and $\Omega(\eps^{-2/3})$ lower bounds for the $1 < p < 2$ and $p=1$ cases, respectively.} Meanwhile, work by Chakrabarty, Graur, Jiang, and Sidford~\cite{chakrabarty2023parallel} presented algorithms that achieve $O(d^{1/3}\eps^{-2/3})$ parallel stochastic query depth for $\ell_\infty$-Lipschitz convex problems, based on leveraging additional Gaussian smoothing for this setting, thereby matching the lower bounds for constant accuracy, though their results rely on stronger bounded second-moment (i.e., variance) assumptions that are standard when analyzing stochastic first-order methods.
Moreover, $\ell_p$-based optimization methods have been shown to be broadly useful across a variety of contexts including training deep neural networks \cite{luo2025stacey, massena2026sgd, pethick2026use}, $p$-norm diffusion and flow problems \cite{adil2020faster, fountoulakis2020p, kelner2014almost,liu2020faster}, and stochastic optimization under heavy-tailed noise assumptions \cite{fatkhullin2025can,song2026can,pmlr-v178-vural22a}.

In this work, we study parallel stochastic $\ell_p$-Lipschitz convex optimization under bounded $q$-th moment assumptions (for $q = \frac{p}{p-1}$), for which we provide new query depth complexities that improve upon previous sequential results for certain ranges of solution accuracy. Underlying our results are two key technical contributions:
\begin{enumerate}
    \item First, we introduce a class of $\ell_p$-optimal-distance-convergence (ODC$_p$) algorithms, which we then show may be used to efficiently estimate the optimum of uniformly convex objectives. Although initially inspired by ideas on stochastic bias-reduced gradient methods~\cite{asi2021stochastic}, we show how extending this type of approach to generalized $\ell_p$-based Moreau envelopes leads to several technical complications that need to be addressed, including the fact that the gradient function for the generalized Moreau envelope is a nonlinear mapping. 

\item Second, we develop and analyze a new $\ell_p$-generalized Gaussian reweighting scheme for gradient estimation, which broadly extends previous approaches based on standard Gaussian convolution. As a result of these tools, we are able to handle several basic challenges that arise when attempting to estimate the gradients of these $\ell_p$-based Moreau envelopes that are central to our techniques.
\end{enumerate} 

By combining these techniques with recent algorithms for $\ell_p$ ball-oracle acceleration~\cite{adil2026convex,contreras2025non}, we arrive at our improved stochastic query depth complexities for parallel $\ell_p$-Lipschitz SCO. In addition, we establish new dimension-dependent lower bounds that improve upon the previous results of \cite{diakonikolas2019lower}, under certain $\ell_p$ domain restrictions.

\subsection{$\ell_p$ Stochastic Bias-Reduced Gradient Methods}
We begin by discussing the setting of bias-reduced gradient methods. 
When considering the problem of minimizing a function $f: \mathcal{X} \rightarrow \mathbb{R}$ given access to a stochastic subgradient estimator $\hat{\nabla} f$ that is unbiased and has bounded variance, a natural question is to ask whether this estimator can be transformed into a nearly unbiased estimator of the true minimizer $x_* = \argmin_{x \in \mathcal{X}} f(x)$.  For strongly convex $f$, \cite{asi2021stochastic} showed how to do so by proposing an estimator $\hat{x}_*$ such that for any $\delta > 0$, it holds that \begin{flalign*}
    &||\mathbb{E}\hat{x}_* - x_*||_2 \le \delta \hspace{0.5em}\AndKeyWord{}\hspace{0.5em} \mathbb{E}||\hat{x}_* - \mathbb{E}\hat{x}_*||_2^2 = O\left(\frac{G^2}{\mu^2} \log\left(\frac{G}{\mu \delta}\right)\right),
\end{flalign*} where $\mathbb{E}\|\widehat{\nabla} f(x)\|_2^2\le G^2$. Given the importance of these estimators and the extensive number of applications that these are used for such as estimating proximal points \cite{parikh2014proximal, monteiro2013accelerated}, minimizing the maximal loss \cite{carmon2021thinking}, and gradient-efficient composite optimization \cite{lan2016gradient} and the prevalence of uniformly convex optimization contexts such as $\ell_p$ regression \cite{bullins2020highly, adil2026convex, adil2019fast, jambulapati2022improved} and online learning \cite{lei2020convergence}, it is natural to ask whether such estimators can be developed for the $\ell_p$ setting. 

In this work, we first answer this in the affirmative, and provide a class of $\ell_p$-optimal-distance-convergence (ODC$_p$) algorithms that, when combined with multilevel Monte Carlo methods \cite{giles2015multilevel, blanchet2015unbiased}, efficiently estimate the optimum of a uniformly convex function with low bias and bounded $q$-th moment in the \emph{dual} $\ell_q$ norm, where $q = \frac{p}{p-1}$.  Further, we use these estimators to yield improved rates for projection-efficient $\ell_p$-Lipschitz optimization, derive complementary rates to \cite{asi2021stochastic} for composite minimization, and develop efficient methods for estimating proximal points, where the proximal point of $y$ is $P_{f, \lambda}(y) = \argmin_{x \in \mathbb{R}^d}\{f(x) + \frac{\lambda}{p} ||x - y||_p^p \}$.   

In Section \ref{section:biased-reduced}, we describe our \odc algorithms, which are part of a larger set of ODC algorithms that take as input an unbiased subgradient estimator which satisfies a given moment bound, and provides a (hopefully) efficient convergence guarantee for convergence to $x_*$. In the Euclidean setting, \cite{asi2021private} defines an ODC algorithm to minimize a strongly convex function $F = f + \psi$, where $f$ is convex and $\psi$ is simple such that the algorithm's output $x = $ ODC$(\hat{\nabla}f, \psi, T)$ satisfies $\mathbb{E}||x - x_*||_2^2 \le c\frac{G^2}{\mu^2 T}$ for some constant $c$. Moreover, such algorithms are widely found in the literature \cite{hazan2011beyond, rakhlin2011making, asi2021private}. Given the advantages that these algorithms have provided, we develop analogous ODC$_p$ algorithms for $p \ge 2$, for minimizing a $\mu$-uniformly convex function $F = f + \psi$, where $f$ is convex. As in \cite{asi2021stochastic}, we leverage the multilevel Monte Carlo (MLMC) technique of \cite{blanchet2015unbiased} to construct a low-bias, bounded-dual-moment estimator of the $\ell_p$-Moreau gradient. However, unlike in the Euclidean-case, when considering $p > 2$, the map from the proximal point to the Moreau gradient is nonlinear. Therefore, we depart from previous analyses, by applying the MLMC construction directly to approximate the Moreau gradient, rather than first constructing a low-bias estimator of the proximal point and then using this to compute the Moreau gradient estimation. We prove that for any $\lambda$-uniformly convex function $F = f + \psi$ where $f$ is convex and iteration budget $T_{\max}$, our ODC$_p$ algorithms provide good estimators $g'$ of the true gradient $g_*$ in a reasonably fast manner (see Lemma \ref{lem:morgrad_mlmc} for more details). 
\paragraph{EpochDA (Section \ref{subsection:odc_p_and_epochda}).} Next, we provide a concrete ODC$_p$ algorithm, EpochDA (Algorithm \ref{alg:epochda_odc}), which extends Algorithm 3 of \cite{nesterov2009primal} to a composite setting using an $\ell_p$-based uniformly convex prox function, and will be useful in achieving efficient $\ell_p$ parallel rates (see Section \ref{section: parallel_upper_bounds} for more details). Within each epoch, it applies a dual-averaging solver to the regularized objective, incorporating stochastic gradient estimates into successive updates while controlling the iterates through the global prox function centered at $\overline{x}$. It uses geometrically increasing epoch lengths and decreasing step sizes, averaging the iterates at the end of each epoch to initialize the next. 
\paragraph{Projection-Efficient $\ell_p$-Lipschitz Optimization
(Section \ref{subsection:projection-efficient}).}
In stochastic nonsmooth optimization, constrained proximal
or projection computations may be substantially more expensive
than stochastic subgradient evaluations
\cite{asi2021stochastic}.
In \Cref{subsection:projection-efficient}, we study this
tradeoff in the $\ell_p$ setting for $p>2$, where $f$ is
convex and accessed
through an unbiased stochastic subgradient oracle
$\widehat{\nabla}f$.
In the Euclidean setting, \cite{asi2021stochastic}
combine Moreau envelope smoothing with accelerated gradient
descent to obtain $O(GR/\epsilon)$ expensive projections
while retaining near-optimal stochastic subgradient complexity.
For $p>2$, we combine our $\ell_p$-Moreau envelope gradient
estimator with accelerated updates using $p$-power regularization,
inspired by \cite{adil2026convex}.
Our analysis uses the smoothness of the envelope with respect
to $\ell_p$, while all constrained updates use $p$-power
regularization.
The resulting algorithm requires
$O_p((GR/\epsilon)^{2p/(p+2)})$ expensive constrained updates
and $\widetilde O_p((GR/\epsilon)^{p(p-1)})$ stochastic
subgradient evaluations.
Thus, compared with direct stochastic mirror descent,
our method reduces the number of expensive constrained
updates at the cost of a larger stochastic subgradient bound.
\paragraph{Convex Composite Minimization (Section~\ref{subsection:convex_composite}).}
We consider the convex composite optimization problem
\begin{flalign*}
    &\min_{x \in \mathcal X} \Psi(x) := \Lambda(x) + f(x),
\end{flalign*}
where $\Lambda$ is $(L,\kappa)$-weakly smooth with respect to
$\|\cdot\|_p$, $f$ satisfies \Cref{Assumption}, and there
exists a minimizer $x^\star$ such that
$\|x_0-x^\star\|_p \le R$.
A closely related setting was studied in
\cite{diakonikolas2024complementary}, who considered composite
minimization with a weakly smooth component and a uniformly
convex component. Their generalized AGD+ framework yields
an iteration complexity of
$O\left((LR^\kappa/\epsilon)^{
\frac{p}{\kappa p-p+\kappa}}\right)$
when specialized to $\ell_p$ geometry without exploiting
uniform convexity.
Building on the stochastic composite AGD framework of
\cite{asi2021stochastic}, we develop an $\ell_p$ extension
using Bregman proximal steps and our reduced-bias gradient
estimators. Our algorithm
(Algorithm~\ref{alg:composite_convex_opt_alg_label})
achieves the same dependence on $\epsilon$ in the number
of $\nabla\Lambda$ evaluations, up to logarithmic factors,
while allowing $f$ to be nonsmooth and accessed through
a stochastic subgradient oracle.

\paragraph{$\ell_p$ Accelerated Proximal Point Methods (Section \ref{subsection:lp_proximal_point}). } Here, we develop Monteiro-Svaiter-type $\ell_p$ accelerated proximal point methods which leverage our reduced-bias Moreau envelope gradient estimators. Inspired by \cite{asi2021stochastic}, we present a stochastic accelerated $\ell_p$ proximal point method (Algorithm \ref{alg:stochastic_accelerated_prox_point}), that combines adaptive proximal regularization with an accelerated first-order update. At each iteration, the algorithm selects a regularization parameter $(\lambda_{k + 1})$, from a designed NextLambda algorithm (Algorithm \ref{alg:next_lambda_alg}), and computes an acceleration weight ($a_{k + 1}$), which determines the cumulative weight ($A_{k + 1}$) and the interpolation point ($y_k$) between the current iterate and the auxiliary sequence. We obtain the next iterate ($x_{k + 1}$) by approximately solving a regularized problem centered at $y_k$, with the algorithm terminating when the cumulative accelerated weight reaches a prescribed threshold or the iteration budget is exhausted.

\subsection{Parallel Optimization in $\ell_p$ Geometries}
\paragraph{Parallel Upper Bounds.} In Section \ref{section: parallel_upper_bounds}, we consider the following problem, simplified for the purposes of the introduction. We assume that there is a convex function $f: \mathbb{R}^d  \rightarrow \mathbb{R}$ whose minimizer lies in an $\ell_p$ ball of radius $R$, which can be queried through a stochastic gradient oracle $g$, where $\mathbb{E} g \in \partial f$ and $\mathbb{E}||g||_q^q \le L^q$, for some $L > 0$. Moreover, throughout the text, unless otherwise specified, $q= \frac{p}{p - 1}$. Our goal is to minimize $f$ up to an expected additive error $\eps_{\mathrm{opt}}$. Furthermore, we assume a setting where multiple queries can be submitted simultaneously. Along with the query depth, we also compute the computational depth and total number of queries to the stochastic gradient oracle to characterize the efficiency of our algorithms. As in \cite{carmon2023resqueing}, the \emph{computational depth} measures the parallel depth of the algorithm, assuming that each oracle query has $O(1)$ depth, and the \emph{total number of queries} refers to the total number of queries made to the stochastic gradient oracle over all rounds. 
\paragraph{Generalized Gaussian.} A key technical tool that we introduce to achieve the parallel rates for our EpochDA and $\ell_p$-AC-SA solvers is an $\ell_p$-Gaussian Reweighting for Gradient Estimation (\greg) estimator. Beginning with a generalized $\ell_p$-norm Gaussian as studied by \cite{bassily2021non}, our estimator computes an unbiased estimate of a convolution of our objective $f$ with this $\ell_p$ Gaussian, which furthermore offers the potential for parallelism.
\paragraph{EpochDA Based Solver.}
Our first algorithm (Theorem~\ref{theorem:epochda_parallel_guarantee}) is based on an epoch-based stochastic solver using our \greg~estimators within the ball acceleration framework. By implementing the stochastic gradient queries in parallel, our algorithm achieves an oracle query depth of $\widetilde O(d^{\frac{p}{2(p+1)}}\epsilon_{\rm opt}^{-\frac{p}{p+1}})$ when $L_p=R_p=1$. This improves upon the $O(\epsilon_{\rm opt}^{-p})$ oracle query depth of stochastic $\ell_p$ mirror descent in the regime $\epsilon_{\rm opt}\leq d^{-\frac{1}{2p}}$, up to logarithmic factors.
\paragraph{$\ell_p$-AC-SA Based Solver.}
Our second algorithm (Theorem~\ref{theorem:p-parallel-acsa}) replaces EpochDA with an accelerated counterpart based on $\ell_p$-AC-SA, achieving a computational depth of $\widetilde O(d^{\frac{p}{2(p+1)}}\epsilon_{\rm opt}^{-\frac{p}{p+1}}+d^{\frac{p}{2(p+2)}}\epsilon_{\rm opt}^{-\frac{2p}{p+2}})$ when $L_p=R_p=1$. This improves upon both the oracle query depth and computational depth of stochastic $\ell_p$ mirror descent in the regime $\epsilon_{\rm opt}<d^{-\frac{1}{2p}}$, up to logarithmic factors. Compared with our EpochDA-based solver, our $\ell_p$ AC-SA variant further reduces computational depth while maintaining the same oracle query depth and total query complexity up to logarithmic factors.

\begin{table*}[t]
\centering
\small
\setlength{\tabcolsep}{8pt}
\renewcommand{\arraystretch}{2.5}
\begin{tabular}{@{}lcc@{}}
\toprule
Method & Query depth & Computational depth\\
\midrule

\makecell[l]{Stochastic Mirror Descent\\
(\cite{pmlr-v178-vural22a}, Theorem 6)}
&
$\begin{aligned}\left(\frac{L_p R_p}{\eps_{\mathrm{opt}}}\right)^{p}\end{aligned}$
&
$\begin{aligned}\left(\frac{L_p R_p}{\eps_{\mathrm{opt}}}\right)^{p}\end{aligned}$
\\[20pt]

\makecell[l]{BallAccel + EpochSGD\\
(\cite{carmon2023resqueing}, Theorem 1)}
&
$\begin{aligned}d^{1/3}\left(\frac{L_2 R_2}{\eps_{\mathrm{opt}}}\right)^{2/3}\end{aligned}$
&
$\begin{aligned}d^{1/3}\left(\frac{L_2 R_2}{\eps_{\mathrm{opt}}}\right)^{2/3} +\left(\frac{L_2 R_2}{\eps_{\mathrm{opt}}}\right)^{2}\end{aligned}$
\\[20pt]

\makecell[l]{BallAccel + AC-SA\\
(\cite{carmon2023resqueing}, Theorem 2)}
&
$\begin{aligned}d^{1/3}\left(\frac{L_2 R_2}{\eps_{\mathrm{opt}}}\right)^{2/3}\end{aligned}$
&
$\begin{aligned}
&d^{1/3}\left(\frac{L_2 R_2}{\eps_{\mathrm{opt}}}\right)^{2/3} +d^{1/4}\left(\frac{L_2 R_2}{\eps_{\mathrm{opt}}}\right)
\end{aligned}$
\\[20pt]

\cite{jambulapati2024closing} (Theorem 2)
&
$\begin{aligned}d^{1/3}\left(\frac{L_2 R_2}{\eps_{\mathrm{opt}}}\right)^{2/3}\end{aligned}$
&
$\begin{aligned}
&d^{1/3}\left(\frac{L_2 R_2}{\eps_{\mathrm{opt}}}\right)^{2/3}
\end{aligned}$
\\[20pt]

\makecell[l]{$\ell_p$ BallAccel + EpochDA\\
(Theorem \ref{theorem:epochda_parallel_guarantee})}
&
$\begin{aligned}d^{\frac{p}{2p+2}}\left(\frac{L_p R_p}{\eps_{\mathrm{opt}}}\right)^{\frac{p}{p+1}}\end{aligned}$
&
$\begin{aligned}
&d^{\frac{p}{2p+2}}\left(\frac{L_p R_p}{\eps_{\mathrm{opt}}}\right)^{\frac{p}{p+1}}+d^{\frac{p(p-2)}{2(p+1)}}
\left(\frac{L_p R_p}{\eps_{\mathrm{opt}}}\right)^{\frac{p(2p - 1)}{p+1}}
\end{aligned}$
\\[20pt]

\makecell[l]{$\ell_p$ BallAccel + $\ell_p$-AC-SA\\
(Theorem \ref{theorem:p-parallel-acsa})}
&
$\begin{aligned}d^{\frac{p}{2p+2}}\left(\frac{L_p R_p}{\eps_{\mathrm{opt}}}\right)^{\frac{p}{p+1}}\end{aligned}$
&
$\begin{aligned}
&d^{\frac{p}{2p+2}}\left(\frac{L_p R_p}{\eps_{\mathrm{opt}}}\right)^{\frac{p}{p+1}}+d^{\frac{p}{2p + 4}}\left(\frac{L_p R_p}{\eps_{\mathrm{opt}}}\right)^{\frac{2p}{p + 2}}
\end{aligned}$
\\

\bottomrule
\end{tabular}
\caption{ Query and computational depths to achieve $\eps_{\mathrm{opt}}$-accurate solution, for general $L_2, R_2, L_p, R_p > 0$, where $L_2$ and $L_p$ represent $\ell_2$ and $\ell_p$ Lipschitz constants respectively and $R_2$ and $R_p$ represent the $\ell_2$ and $\ell_p$ norm bounds on the feasible set (omitting additional polylogarithmic factors).}
\label{tab:comparison}
\end{table*}

\begin{table*}[t]
\centering
\small
\setlength{\tabcolsep}{8pt}
\begin{tabular}{@{}lc@{}}
\toprule
Method & Total query complexity \\
\midrule

\makecell[l]{Stochastic Mirror Descent\\
(\cite{pmlr-v178-vural22a}, Theorem 6)}
&
$\begin{aligned}\left(\frac{L_p R_p}{\eps_{\mathrm{opt}}}\right)^{p}\end{aligned}$
\\[20pt]

\makecell[l]{BallAccel + EpochSGD\\
(\cite{carmon2023resqueing}, Theorem 1)}
&
$\begin{aligned}d^{1/3}\left(\frac{L_2 R_2}{\eps_{\mathrm{opt}}}\right)^{2/3} +\left(\frac{L_2 R_2}{\eps_{\mathrm{opt}}}\right)^{2}
\end{aligned}$
\\[20pt]

\makecell[l]{BallAccel + AC-SA\\
(\cite{carmon2023resqueing}, Theorem 2)}
&
$\begin{aligned}
&d^{1/3}\left(\frac{L_2 R_2}{\eps_{\mathrm{opt}}}\right)^{2/3} +\left(\frac{L_2 R_2}{\eps_{\mathrm{opt}}}\right)^{2}
\end{aligned}$
\\[20pt]

\cite{jambulapati2024closing} (Theorem 2)
&
$\begin{aligned}
&d^{1/3}\left(\frac{L_2 R_2}{\eps_{\mathrm{opt}}}\right)^{2/3} +\left(\frac{L_2 R_2}{\eps_{\mathrm{opt}}}\right)^{2}
\end{aligned}$
\\[20pt]

\makecell[l]{$\ell_p$ BallAccel + EpochDA\\
(Theorem \ref{theorem:epochda_parallel_guarantee})}
&
$\begin{aligned}
&d^{\frac{p}{2p+2}}\left(\frac{L_p R_p}{\eps_{\mathrm{opt}}}\right)^{\frac{p}{p+1}}+d^{\frac{p(p-2)}{2(p+1)}}
\left(\frac{L_p R_p}{\eps_{\mathrm{opt}}}\right)^{\frac{p(2p - 1)}{p+1}}
\end{aligned}$
\\[20pt]

\makecell[l]{$\ell_p$ BallAccel + $\ell_p$-AC-SA\\
(Theorem \ref{theorem:p-parallel-acsa})}
&
$\begin{aligned}d^{\frac{p}{2p+2}}\left(\frac{L_p R_p}{\eps_{\mathrm{opt}}}\right)^{\frac{p}{p+1}}
+d^{\frac{p}{2p + 4}}\left(\frac{L_p R_p}{\eps_{\mathrm{opt}}}\right)^{\frac{2p}{p + 2}}
+d^{\frac{p(p-2)}{2(p+1)}}
\left(\frac{L_p R_p}{\eps_{\mathrm{opt}}}\right)^{\frac{p(2p - 1)}{p+1}}\end{aligned}$
\\
\bottomrule
\end{tabular}
\caption{ Total query complexities to achieve $\eps_{\mathrm{opt}}$-accurate solution, for general $L_2, R_2, L_p, R_p > 0$, where $L_2$ and $L_p$ represent $\ell_2$ and $\ell_p$ Lipschitz constants respectively and $R_2$ and $R_p$ represent the $\ell_2$ and $\ell_p$ norm bounds on the feasible set (omitting additional polylogarithmic factors).}
\label{tab:comparison}
\end{table*}
\paragraph{Parallel Lower Bounds.}
Several prior works have established fundamental limitations on parallelizing nonsmooth convex optimization. \cite{balkanski2018parallelization} developed an adaptivity-based lower-bound framework, showing that polynomially many parallel queries cannot improve upon sequential convergence rates for up to \(\widetilde{\Omega}(d^{1/3})\) rounds in Euclidean geometry. Subsequently, \cite{bubeck2019complexity}  strengthened the Euclidean lower bound to \(\widetilde{\Omega}(d^{1/2})\) by introducing a wall construction. For general \(\ell_p\) geometries with \(p\ge 2\), \cite{diakonikolas2019lower} established a lower bound of \(\widetilde{\Omega}(\min\{\varepsilon^{-p},\varepsilon^{2/3}d^{1/3}\})\). Building on the proof framework of \cite{balkanski2018parallelization} and incorporating the wall function adapted to \(\ell_p\) geometry, we establish an improved lower bound of
\[
\widetilde{\Omega}_p\left(\min\left\{\varepsilon^{-p},d^{\frac{p}{3p-2}}\right\}\right)
\]
for nonsmooth convex optimization over a subset of the unit \(\ell_p\) ball, with polynomially many parallel local-oracle queries per round. Our result improves the dimension-dependent rate obtained by adapting the construction of \cite{diakonikolas2019lower} to the same feasible set, while recovering the known \(d^{1/2}\) exponent when \(p=2\). However, our current upper bound retains a dimension exponent of \(1/2\), whereas our lower bound has exponent \(p/(3p-2)<1/2\). Closing this gap and determining the optimal parallel complexity in general \(\ell_p\) geometry remains an open question.

\subsection{Paper Organization and Notation}
In Section \ref{section:gg_intro}, we introduce an $\ell_p$-generalized Gaussian distribution and present the \greg~estimator, which estimates the gradient of a convolution of $f$ with this $\ell_p$-generalized Gaussian. In Section \ref{section:biased-reduced}, we introduce our \odc algorithms, and discuss their applications to projection-efficient (Section \ref{subsection:projection-efficient}) and convex composite (Section \ref{subsection:convex_composite}) optimization, in addition to efficient $\ell_p$-proximal point estimation (Section \ref{subsection:lp_proximal_point}).

In Section \ref{section: parallel_upper_bounds}, we present two algorithms, EpochDA and $\ell_p$-AC-SA and present upper bounds on their query depths, computational depths, and total number of gradient queries. In Section \ref{section:parallel_lower_bounds}, we present our lower bounds for the number of adaptive rounds needed to acheive at most $\eps$ error. 
\newline \newline 
Throughout the paper, we use $\tilde{O}(\cdot)$ to hide polylogarithmic factors in the parameters and will refer to $p$ and $q$ as Hölder conjugate exponents (i.e., $\frac{1}{p} + \frac{1}{q} = 1$). Moreover, unless otherwise specified, we consider $p \in [2, \infty)$. We use $\mathbb{B}_p(\overline{x}, r)$ to denote an $\ell_p$ ball around $\overline{x}$ of radius $r$, i.e., the set $\{x' \in \mathbb{R}^d: ||x' - \overline{x}||_p \le r \}$. Furthermore, for $a \in (0, 1]$, let $\mathrm{Geom}(a)$ denote the geometric distribution with parameter $a$. We also refer to $\mathcal{N}(\mu, \sigma^2 I)$ as a multivariate Gaussian distribution with mean $\mu \in \mathbb{R}^d$ and covariance $\sigma^2 I$, where $I$ denotes the identity matrix in $\mathbb{R}^{d \times d}$. 
\newline \newline 
We say that a function $f: \mathbb{R}^d \rightarrow \mathbb{R}$ is $L$-Lipschitz with respect to the $\ell_p$-norm if, for all $x, y \in \mathbb{R}^d$, it holds that $|f(x) - f(y)| \le L||x - y||_p$. We say that $f$ is $(L, \kappa)$-weakly smooth with respect to $||\cdot ||_p$ for $L > 0$ and $\kappa \in (1, 2]$, if its gradient is $(L, \kappa - 1)$-Hölder continuous, i.e., \begin{flalign*}
    &\forall x, y \in \mathbb{R}^d, ||\nabla f(x) - \nabla f(y)||_q \le L ||x - y||_p^{\kappa - 1}.
\end{flalign*} We say that $f$ is convex if for all $x, y \in \mathbb{R}^d$ and $\lambda \in [0, 1]$, $f(\lambda x + (1 - \lambda)y) \le \lambda f(x) + (1 - \lambda) f(y)$. Moreover, we say that $f$ is uniformly convex with convexity parameter $\mu$ and order $p$ if for all $x, y \in \mathbb{R}^d$, it holds that \begin{flalign*}
    &f(y) \ge f(x) + \langle \nabla f(x), y - x\rangle + \frac{\mu}{p} ||y - x||_p^p.
\end{flalign*} 

\section{Framework}
\label{section:gg_intro}
\subsection{Generalized $\ell_p$ Gaussian}
First, we state the density of our Generalized $\ell_p$ Gaussian ($\mathcal{GG}_p$) \cite{bassily2021non}. 
 For any $X \sim \mathcal{GG}_p(\mu, \sigma^2)$, its density is represented as  \[\Bigg(\frac{(2\Gamma(1 + 1/p))^d}{\Gamma(1 + d/p)}(2\sigma^2)^{d/2}\Gamma(1 + d/2)\Bigg)^{-1}\exp\bigg(-\frac{\|x-\mu\|_p^2}{2\sigma^2}\bigg).\] We will use $\gamma_{\rho, p}$ to denote the probability density function of $\mathcal{GG}_p(0, \rho^2)$. Now, we define the convolution operation for our distribution. \begin{definition}\label{def:convolution}($\mathcal{GG}_p$-Convolution) For a function $f: \mathbb{R}^d \rightarrow \mathbb{R}$, we denote its convolution with a $\mathcal{GG}_p$ random variable of mean $0$ and variance $\rho^2$ by \begin{flalign*}
    &\hat{f}_{\rho,p}(x) := \mathbb{E}_{y \sim \mathcal{GG}_p(0, \rho^2)}f(x + y) = \int_{y \in \mathbb{R}^d} f(x - y) \gamma_{\rho, p}(y) dy.
\end{flalign*}
    
\end{definition}
\subsection{\greg~Estimator}
\label{section:grad_estimator_definition}
We now present our \greg~($\ell_p$-Generalized Gaussian Reweighting for Gradient Estimation) estimator for the $\ell_p$ norm, which is related to the ReSQue estimator~\cite{carmon2023resqueing}, as well as the estimator in~\cite{bubeck2019complexity}. 
\begin{definition} (\greg~Estimator) 
\label{defn:estimator definition}
Let $\overline{x} \in \mathbb{R}^d$ and let $f: \mathbb{R}^d \rightarrow \mathbb{R}$ be convex. Suppose we have a gradient estimator $g: \mathbb{R}^d \rightarrow \mathbb{R}^d$ such that $\mathbb{E}[g] \in \partial f$. We define the \emph{\greg}~estimator of radius $\rho$ as the following: \begin{flalign*}
    &\tilde{\nabla}_{\overline{x}}^g \hat{f}_{\rho,p}(x) := \frac{\gamma_{\rho, p}(x - \overline{x} - \xi)}{\gamma_{\rho, p}(\xi)} g(\overline{x} + \xi), \xi \sim \mathcal{GG}_p(0, \rho^2)
\end{flalign*} where $\xi$ is sampled first and $g$ is independently queried at $\overline{x} + \xi$. 
\end{definition} 
We next prove that if $g$ is an unbiased estimate of $\partial f$ and has a certain moment bound, then our estimator is also unbiased for the convolved function and enjoys a similar moment bound. 

\begin{lemma}\label{lem:GREG-q moment}
    The following statements hold about our estimator presented in Definition \ref{defn:estimator definition}: \begin{itemize}
        \item $\mathbb{E}[\tilde{\nabla}_{\overline{x}}^g \hat{f}_{\rho, p}(x)] = \nabla \hat{f}_{\rho,p}(x)$
        \item If $\mathbb{E}||g(z)||_q^q \le G^q$ for all $z \in \R^d$ and $x \in \mathbb{B}_{p}(\overline{x},\rho)$, 
        then $\mathbb{E}||\tilde{\nabla}_{\overline{x}}^g \hat{f}_{\rho, p}(x)||_q^q \leq 3G^q$ 
    \end{itemize}
\end{lemma}
\begin{proof}
To prove the first statement, we note that 
\begin{flalign*}
\mathbb{E}[\tilde{\nabla}_{\overline{x}}^g \hat{f}_{\rho, p}(x)] = \int \frac{\gamma_{\rho, p}(x - \overline{x} - \xi)}{\gamma_{\rho, p}(\xi)} \mathbb{E}[g(\overline{x} + \xi)] \gamma_{\rho, p}(\xi) d\xi &= \int \gamma_{\rho, p}(x - \overline{x} - \xi) \mathbb{E}[g(\overline{x} + \xi)] d\xi \\ 
&= \int \E 
\left[g(\overline{x} + \xi)\right] \cdot \gamma_{\rho, p}(x - \overline{x} - \xi) d\xi = \nabla \hat{f}_{\rho,p}(x)
\end{flalign*} To show the second statement, let $v=x-\bar x$. By definition,
\[\left\|\widetilde \nabla_{\bar x}^{g}\widehat f_{\rho,p}(x)\right\|_q^q=\left(
\frac{\gamma_{\rho,p}(v-\xi)}
{\gamma_{\rho,p}(\xi)}\right)^q
\|g(\bar x+\xi)\|_q^q.\]
Taking expectation first with respect to the randomness of $g$ and using
$\mathbb E_g\|g(z)\|_q^q\le G^q$, we get \begin{align*}
\mathbb E\left\|\widetilde \nabla_{\bar x,p}^{g}\widehat f_{\rho,p}(x)\right\|_q^q &\leq G^q
\mathbb E_{\xi\sim \gamma_{\rho,p}} \left[\left(
\frac{\gamma_{\rho,p}(v-\xi)}{\gamma_{\rho,p}(\xi)}
\right)^q\right] = G^q \int\frac{\gamma_{\rho,p}(v-\xi)^q}{\gamma_{\rho,p}(\xi)^{q-1}}d\xi \\
&= G^q \bigg(\frac{(2\Gamma(1 + 1/p))^d}{\Gamma(1 + d/p)}(2\rho
^2)^{d/2}\Gamma(1 + d/2)\bigg)^{-1}\int\exp \Big( \frac{\|\xi\|_p^2 - p \|v - \xi\|_p^2}{2(p-1)\rho^2}\Big)d\xi
\end{align*} Use the fact that $f(x) = \frac{1}{2}\|x\|_p^2$ is $(p-1)$-smooth with respect to the $p$ norm, by Theorem 5.8 of \cite{beck2017first}, we have
\[f(\lambda x + (1-\lambda)y)
\geq \lambda f(x) + (1-\lambda)f(y)
- \frac{p-1}{2}\lambda(1-\lambda)\|x-y\|^2_p,\quad \forall x,y\in \mathbb{R}^d,\ \lambda\in[0,1].\]
Let $\lambda = (p-1)/p$, $x = \xi - \frac{1}{\lambda} v$ and $y = \xi$ one can get \begin{align*}
f(\xi - v) &= f(\lambda x + (1-\lambda)y)
\geq \lambda f(x) + (1-\lambda)f(y)
- \frac{p-1}{2}\lambda(1-\lambda)\|x-y\|_p^2 \\
&= \frac{p-1}{p}f(\xi - \frac{1}{\lambda}v) + \frac{1}{p}f(\xi) - \frac{(p-1)^2}{2p^2}\|\frac{p}{p-1}v\|_p^2
\end{align*}
By rearranging and multiplying both sides by $p/(p-1)$, this yields
\[\frac{\|\xi\|_p^2}{2(p-1)} - \frac{p\|\xi - v\|_p^2}{2(p-1)} \leq -\frac{\|\xi - \frac{1}{\lambda}v\|_p^2}{2} + \frac{p\|v\|_p^2}{2(p-1)}\]
Using this inequality we have
\begin{align*}
\mathbb{E}\bigl\|\widetilde{\nabla}_{\bar{x}}^g
\widehat f_{\rho,p}(x)\bigr\|_q^q &\leq G^q \bigg(\frac{(2\Gamma(1 + 1/p))^d}{\Gamma(1 + d/p)}(2\rho^2)^{d/2}\Gamma(1 + d/2)\Bigg)^{-1}\exp\Big( \frac{p\|v\|_p^2}{2(p-1)\rho^2}\Big) \int\exp\Big( -\frac{\|\xi - \frac{1}{\lambda}v\|_p^2}{2\rho^2}\Big)d\xi \\
&\leq G^q\exp\Big( \frac{p\|v\|_p^2}{2(p-1)\rho^2}\Big) \leq 3G^q
\end{align*}
where the last inequality uses that $\|v\|_p \leq \rho$.
\end{proof} 

We also prove the following bounds on higher-order moments for two important quantities. 
\begin{lemma}
\label{lemma:GGguaranteespart2}
If $x, x' \in \mathbb{B}_{\overline{x}}(\rho/m)$ for $m \ge 2$, then \begin{flalign*}
        &\mathbb{E}_{\xi \sim \mathcal{GG}_p(0, \rho^2)}\bigg[\bigg(\frac{\gamma_{\rho, p}(x - \overline{x} - \xi)}{\gamma_{\rho, p}(\xi)}\bigg)^m\bigg] \leq e^{p-1}\\
        &\mathbb{E}_{\xi \sim \mathcal{GG}_p(0, \rho^2)}\bigg[\bigg|\frac{\gamma_{\rho, p}(x - \overline{x} - \xi) - \gamma_{\rho, p}(x' - \overline{x} - \xi)}{\gamma_{\rho, p}(\xi)}\bigg|^m\bigg] \leq e^{p-1}\left(\frac{C_p \sqrt{m}\|x-x'\|_p}{\rho}\right)^m
    \end{flalign*} for some constant $C_p > 0$.    
\end{lemma}

\begin{proof}
Denote $v = x - \bar{x}$ and $v' = x'-\bar{x}$. By definition,
    \begin{align*}
        \mathbb{E}_{\xi \sim \mathcal{GG}_p(0, \rho^2)}\bigg[\bigg(\frac{\gamma_{\rho, p}(x - \overline{x} - \xi)}{\gamma_{\rho, p}(\xi)}\bigg)^m\bigg] &= \mathbb{E}_{\xi \sim \mathcal{GG}_p(0, \rho^2)}\bigg[\exp\bigg(\frac{m(\|\xi\|_p^2 - \|v - \xi\|_p^2)}{2\rho^2}\bigg)\bigg] \\
        &\leq \exp\bigg(\frac{m(p-1)\|v\|_p^2}{2\rho^2}\bigg) \cdot \mathbb{E}_{\xi \sim \mathcal{GG}_p(0, \rho^2)}\bigg[\exp\bigg(\frac{m}{\rho^2}\langle \nabla \frac{1}{2}\|\xi\|_p^2,v\rangle\bigg)\bigg]
    \end{align*} where the inequality is by the smoothness of $\frac{1}{2}\|\cdot\|_p^2$. Since $\|v\|_p \leq \rho/m$, we have, with Fact \ref{fact:expectedexpbound},
    \begin{align*}
        \mathbb{E}_{\xi \sim \mathcal{GG}_p(0, \rho^2)}\bigg[\bigg(\frac{\gamma_{\rho, p}(x - \overline{x} - \xi)}{\gamma_{\rho, p}(\xi)}\bigg)
        ^m\bigg] &\leq \exp\bigg(\frac{m(p-1)\|v\|_p^2}{2\rho^2}\bigg) \cdot \mathbb{E}_{\xi \sim \mathcal{GG}_p(0, \rho^2)}
        \bigg[\exp\bigg(\frac{m}{\rho^2}\langle \nabla \frac{1}{2}\|\xi\|_p^2
        ,v\rangle\bigg)\bigg] \\
        &\leq \exp\bigg(\frac{p-1}{2m}\bigg) \exp\bigg(\frac{ (p-1)m^2}{2\rho^2}\|v\|_p^2\bigg) \leq \exp\bigg(\frac{p-1}{2m} + \frac{p-1}{2}\bigg) \\
        &\leq e^{p-1}
    \end{align*}  For the second inequality, we define the line segment
    $v_t:=v'+t(v-v')$ for $t\in[0,1]$. By the fundamental theorem of calculus,
    \[\frac{\gamma_{\rho, p}(v-\xi)-\gamma_{\rho, p}(v'-\xi)}{\gamma_{\rho, p}(\xi)}
    = \int_0^1 \frac{\gamma_{\rho, p}(v_t-\xi)}{\gamma_{\rho, p}(\xi)}
    \frac{\langle \nabla\Phi(\xi-v_t),\,v - v'\rangle}{\rho^2}\,dt,\]
    where
    $\Phi(z):=\frac12\|z\|_p^2$. Taking absolute values and applying Jensen's inequality (convexity of $|\cdot|^m$) gives
    \[\left|\frac{\gamma_{\rho, p}(v-\xi)-\gamma_{\rho, p}(v'-\xi)}{\gamma_{\rho, p}(\xi)}
    \right|^m\leq \int_0^1
    \left(\frac{\gamma_{\rho, p}(v_t-\xi)}{\gamma_{\rho, p}(\xi)}\right)^m\left(
    \frac{|\langle \nabla\Phi(\xi-v_t),\,v - v'\rangle|}{\rho^2}\right)^m dt\] Taking expectations with respect to $\xi\sim \mathcal{GG}_p(0,\rho^2)$ and applying Fubini's theorem (as the integral is nonnegative) 
    we get \[\begin{aligned}
    \mathbb{E}_{\xi}
    \left[\left|\frac{\gamma_{\rho,p}(v-\xi)-\gamma_{\rho,p}(v^\prime-\xi)}
    {\gamma_{\rho,p}(\xi)}
    \right|^m\right]\le\int_0^1 \mathbb{E}_{\xi}\!\left[
    \left(\frac{\gamma_{\rho,p}(v_t-\xi)}{\gamma_{\rho,p}(\xi)}\right)^m    \left(\frac{|\langle\nabla\Phi(\xi-v_t),v-v^\prime\rangle|}
    {\rho^2}\right)^m\right]\,dt .
    \end{aligned}
    \] Now, let $z = \xi - v_t$, which implies that $\xi = z + v_t$. Since $\Phi$ is even, it holds that $\gamma_{\rho,p}(v_t-\xi)=\gamma_{\rho,p}(-z)=\gamma_{\rho,p}(z)$. Therefore, letting  $w = v - v' = x - x'$, \begin{flalign*}
        &\mathbb{E}_{\xi} \left[\left(\frac{\gamma_{\rho, p}(v_t - \xi)}{\gamma_{\rho, p}(\xi)}\right)^m \left(\frac{|\langle \nabla \Phi(\xi - v_t),  v - v' \rangle|}{\rho^2}\right)^m\right] = \int_{\mathbb{R}^d} \frac{\gamma_{\rho, p}(\xi) \gamma_{\rho, p}\left(v_t - \xi\right)^m}{\gamma_{\rho, p}\left(\xi\right)^m} \left(\frac{|\langle \nabla \Phi(\xi - v_t), w \rangle|}{\rho^2}\right)^m d \xi \\   &\qquad\qquad\qquad\qquad= \int_{\mathbb{R}^d} \gamma_{\rho, p}(z) \left(\frac{\gamma_{\rho, p}(z)}{\gamma_{\rho, p}(z + v_t)}\right)^{m - 1} \left(\frac{|\langle \nabla \Phi(z), w \rangle |}{\rho^2}\right)^m dz \\ &\qquad\qquad\qquad\qquad= \mathbb{E}_{z \sim \mathcal{GG}_p(0, \rho^2)}\left[\exp\left(\frac{(m - 1) (\Phi(z + v_t) - \Phi(z))}{\rho^2}\right)\left(\frac{|\langle \nabla \Phi(z), w \rangle |}{\rho^2}\right)^m\right] 
    \end{flalign*}  By the $(p - 1)$-smoothness of $\Phi$, it holds that \begin{flalign*}
        &\Phi(z + v_t) - \Phi(z) \le \langle \nabla \Phi(z), v_t \rangle + \frac{p - 1}{2} ||v_t||_p^2. 
    \end{flalign*} Therefore \begin{flalign*}
        &\mathbb{E}_{z \sim \mathcal{GG}_p(0, \rho^2)}[\exp(\frac{(m - 1) (\Phi(z + v_t) - \Phi(z))}{\rho^2})(\frac{|\langle \nabla \Phi(z), w \rangle |}{\rho^2})^m]  \\ 
        &\le \exp(\frac{(m - 1)(p - 1)}{2\rho^2} ||v_t||_p^2) \cdot \mathbb{E}_z [\exp(\frac{m - 1}{\rho^2} \langle \nabla \Phi(z), v_t \rangle )(\frac{|\langle \nabla \Phi(z), w \rangle |}{\rho^2})^m] 
    \end{flalign*} Letting \begin{flalign*}
        &A = \frac{m - 1}{\rho^2} \langle \nabla \Phi(z), v_t \rangle, X = \frac{\langle \nabla \Phi(z), w \rangle}{\rho^2}
    \end{flalign*} we have that by Cauchy-Schwarz, $\mathbb{E}[e^A |X|^m] \le \sqrt{\mathbb{E}[e^{2A}]} \sqrt{\mathbb{E}[|X|^{2m}]}$. We bound the two terms separately. For the first term, using Fact \ref{fact:expectedexpbound}, it holds that \begin{flalign*}
        &\mathbb{E}[e^{2A}] = \mathbb{E}_z [\exp(\frac{\langle \nabla \Phi(z), 2(m - 1) v_t \rangle}{\rho^2})] \le \exp(\frac{2(p - 1)(m - 1)^2}{\rho^2} ||v_t||_p^2) \le \exp(\frac{2(p - 1) (m - 1)^2}{m^2})
    \end{flalign*} where we used $||v_t||_p \le (1 - t) ||v'||_p + t||v||_p \le \rho/m$. Therefore, \begin{flalign*}
        &\sqrt{\mathbb{E}[e^{2A}]} \le \exp(\frac{(p - 1)(m - 1)^2}{m^2})
    \end{flalign*} We now bound the second term. First, we show that $X$ is sub-Gaussian. To prove this, we first prove that $Y(\xi) = \langle \nabla \Phi(\xi), v - v' \rangle$ is sub-Gaussian. First, since $\gamma_{\rho, p}$ is even and $\nabla \Phi$ is odd, it can be shown that $Y(-\xi) = -Y(\xi)$, which means that $\mathbb{E}[Y(\xi)] = 0$. Again, letting $w = v - v'$, and using Fact \ref{fact:expectedexpbound} with $u = \lambda \rho^2 w$ for $\lambda \in \mathbb{R}$, it holds that \begin{flalign*}
        &\mathbb{E}_{\xi} [\exp(\lambda Y)] = \mathbb{E}_{\xi}\left[\exp\left(\frac{1}{\rho^2} \langle \nabla \Phi(\xi), \lambda \rho^2 w \rangle \right)\right] \le \exp(\frac{p - 1}{2\rho^2} \lambda^2 \rho^4 ||w||_p^2) 
    \end{flalign*} where we use the fact that $||w||_p \le 2\rho/m$. By point 4 of proposition 2.5.2 of \cite{vershynin2022high}, we have that $Y$ is sub-gaussian with $K_4^2 := \frac{(p - 1)  \rho^2||w||_p^2}{2}$. By Markov's inequality, we have that \begin{flalign*}
        &\Pr(Y \ge t) = \Pr(e^{\lambda Y} \ge e^{\lambda t}) \le e^{-\lambda t} \mathbb{E}[e^{\lambda Y}] \le \exp(K_4^2 \lambda^2 - \lambda t)
    \end{flalign*} Minimizing the exponent over $\lambda \ge 0$ gives that $\Pr(Y \ge t) \le \exp(-t^2/4K_4^2)$. Applying the same argument with $-Y$ instead of $Y$, it holds that \begin{flalign*}
        &\Pr(|Y| \ge t) \le 2\exp(-\frac{t^2}{4K_4^2}) := 2\exp(-\frac{t^2}{K_1^2})
    \end{flalign*} so property 1 of proposition 2.5.2 of \cite{vershynin2022high} holds with $K_1 := 2K_4$. Now, it also holds that for all $k \ge 2$, \begin{flalign*}
        &\mathbb{E}|Y|^k \le \int_{0}^{\infty}kt^{k - 1} \Pr(|Y| \ge t) dt \le \int_{0}^\infty 2kt^{k - 1} \exp(-\frac{t^2}{K_1^2})dt = kK_1^k \Gamma(\frac{k}{2}) \le k K_1^k \cdot (\frac{k}{2})^{\frac{k}{2}}
    \end{flalign*} which implies that $(\mathbb{E}|Y|^k)^{1/k} \le K_1 \sqrt{2k} := K_2 \sqrt{k}$ by property 2 of proposition 2.5.2 of \cite{vershynin2022high}, so $K_2 = K_1 \sqrt{2} = 2K_4 \sqrt{2}$. By this property, it holds that \begin{flalign*}
        &\sqrt{\mathbb{E}|X|^{2m}} = \frac{\sqrt{\mathbb{E}|Y|^{2m}}}{\rho^{2m}} \le \frac{(K_2 \sqrt{2m})^m}{\rho^{2m}} = \frac{(4K_4 \sqrt{m})^m}{\rho^{2m}} \le (\frac{4\sqrt{p - 1} \rho ||v - v'||_p}{\sqrt{2}} \sqrt{m})^m \frac{1}{\rho^{2m}} = (\frac{C_p\sqrt{m} ||v - v'||_p}{\rho})^m
    \end{flalign*} where $C_p = 2\sqrt{2(p - 1)}$. Therefore, \begin{flalign*}
        &\mathbb{E}[e^A |X|^m] \le \exp(\frac{(p - 1)(m - 1)^2}{m^2})(\frac{C_p \sqrt{m} ||w||_p}{\rho})^m. 
    \end{flalign*} Thus, \begin{flalign*}
        &\mathbb{E}_{z \sim \mathcal{GG}_p(0, \rho^2)}[\exp(\frac{(m - 1) (\Phi(z + v_t) - \Phi(z))}{\rho^2})(\frac{|\langle \nabla \Phi(z), w \rangle |}{\rho^2})^m]  \\ 
        &\le \exp(\frac{(m - 1)(p - 1)}{2\rho^2} ||v_t||_p^2) \cdot \exp(\frac{(p - 1)(m - 1)^2}{m^2})(\frac{C_p \sqrt{m} ||w||_p}{\rho})^m \\ 
        &\le \exp(\frac{(m - 1)(p - 1)}{2m^2}) \exp(\frac{(p - 1)(m - 1)^2}{m^2})(\frac{C_p \sqrt{m} ||w||_p}{\rho})^m \\ 
        &= \exp((p - 1) \frac{(m - 1)(m - \frac{1}{2})}{m^2}) (\frac{C_p \sqrt{m} ||w||_p}{\rho})^m \\ 
        &\le e^{p - 1} (\frac{C_p \sqrt{m} ||v - v'||_p}{\rho})^m.
    \end{flalign*} Integrating over $t \in [0, 1]$ completes the proof. \qedhere
\end{proof}

\section{$\ell_p$ Stochastic Bias-Reduced Methods}
\label{section:biased-reduced}
In this section, we introduce a class of $\ell_p$-optimal-distance-convergence (ODC$_p$) algorithms, detailing the required convergence guarantee and giving a concrete algorithm for computing a low-bias estimator of an order-$p$ uniformly convex function $F$. 
We then use these algorithms in conjunction with multilevel Monte Carlo (MLMC) methods to use them for various applications.  We assume that $F$ can be represented in the form $F = f + \psi$, where $\psi$ is a simple function and $f$ satisfies Assumption \ref{Assumption}. 
\begin{assumption} 
\label{Assumption}
The function $f:\mathcal{X} \rightarrow \R$ is convex, where $\mathcal{X}$ is a closed convex set. We assume that there is an unbiased subgradient estimator $\hat{\nabla}f$ that has bounded $q_{th}$-moment $\E[\|\hat{\nabla}f(x)\|_q^q] \leq G^q$ for all $x \in \mathcal{X}$.
\end{assumption}
\subsection{ODC$_p$ algorithms}
\label{subsection:odc_p_and_epochda}
 
Our estimator uses an arbitrary optimization routine for $F$ as a black-box component, provided that it achieves a sufficiently fast convergence rate toward $\argmin_{x \in \mathcal{X}} F(x)$. 
We describe the required convergence property below. 
\begin{definition}[\odc algorithm]
\label{defn:odc_alg}
For $p\geq 2$, an $\ell_p$-optimal-distance-convergence (\odc) algorithm takes as input an unbiased stochastic subgradient oracle $\hat{\nabla}f$ for a function $f$ satisfying \Cref{Assumption}, a simple function $\psi$, and an oracle budget $T\geq 1$. Suppose that $F=f+\psi$ is $\mu$-uniformly convex of order $p$ with minimizer $x^\star$. The algorithm returns a point $x$ using at most $T$ oracle calls, satisfying
\[
\mathbb{E}\|x-x^\star\|_p^p
\leq \frac{cG^q}{\mu^q T^{1/(p-1)}},
\]
where $c>0$ is a constant.
\end{definition}

As in the work of~\cite{asi2021stochastic}, we leverage the multilevel
Monte Carlo (MLMC) technique of Blanchet and Glynn~\cite{blanchet2015unbiased}
to construct a low-bias, bounded-dual-moment estimator of the
$\ell_p$-Moreau gradient from an $\mathrm{ODC}_p$ algorithm, where we define the $\ell_p$-Moreau envelope to be
\[
f_{p,\lambda}(y) := \min_x\{f(x) + \frac{\lambda}{p}\|x-y\|_p^p\}\]
For $p>2$, however, the map from the proximal point to the Moreau gradient is nonlinear.  We therefore apply the MLMC construction directly to the Moreau-gradient approximations, rather than first constructing a low-bias estimator of the proximal point.
Fix $y\in\mathcal X$ and $\lambda>0$, and define $\psi(x):=
\frac{\lambda}{p}\|x-y\|_p^p$.
For every $j\geq 0$, let
$x_j:=\mathrm{ODC}_p\bigl(
\widehat{\nabla}f,\psi,2^j\bigr)$, and define the corresponding approximate Moreau gradient 
\[g_j:=\lambda J_p(y-x_j),\qquad
J_p(z):=\bigl(|z_k|^{p-2}z_k\bigr)_{k=1}^d.\]
Let $\mathbb P(J=j)=2^{-j}$ for $j \in \{1,2,\dots\}$, and define
\begin{equation}\label{eq:morgrad-mlmc-sample}
g':=g_0+\begin{cases}
2^J\bigl(g_J-g_{J-1}\bigr),
& J\leq j_{\max},\\[4pt]
0,& J>j_{\max}.
\end{cases}
\end{equation} where $j_{\max} = \lfloor \log_2 T_{\max} \rfloor$. 
\begin{algorithm}[H] \caption{$\ell_p$-\textsc{MorGradEst}} \label{alg:lp-morgradest} \begin{algorithmic}[1] \State \textbf{Input:} $\widehat{\nabla}f,\psi,\lambda,\delta,\sigma^q,\mathcal X$ \State \textbf{Parameters:} $\delta$: desired bias, $\sigma^q$: desired $q$-th moment error \State $T_{\max} = \left \lceil (\frac{2C_1 G}{\min \{\delta, \sigma \}})^{p(p -1 )} \right \rceil $ (see \eqref{eq:C1defn})

and $N = \left \lceil 2^p C_3 \frac{G^p}{\sigma^p} T_{\max}^{\frac{p - 2}{p - 1}} \right \rceil$ (see \eqref{eq:C3defn})

\For{$i=1,\dots,N$}  \State $\hat{g}^{(i)} \gets$ draw of estimator \eqref{eq:morgrad-mlmc-sample}\EndFor \State \Return $\widehat{\nabla}f_{p,\lambda}(y) = \frac1N \sum_{i=1}^N \widehat g^{(i)}$ \end{algorithmic} \end{algorithm}

\begin{lemma}
\label{lem:morgrad_mlmc}
Let $f$ and $\widehat{\nabla}f$ satisfy \Cref{Assumption}, let
$y\in\mathcal X$, $p > 2$, and let $\lambda>0$.
Let $\psi(x)=\frac{\lambda}{p}\|x-y\|_p^p$, suppose
$F=f+\psi$ has minimizer $x_*$ over $\mathcal X$, and let
$T_{\max}\in\mathbb N$ with $T_{\max}\ge1$.
For any \odc algorithm on $\mathcal X$ with constant $c$,
let $g'$ be the estimator in \eqref{eq:morgrad-mlmc-sample},
where $J$ is independent of all oracle randomness.
The following statements hold, for some constant $C_p>0$
depending only on $p$ and $c$:
\[
\bigl\|\E g'-g_*\bigr\|_q
\le C_p G T_{\max}^{-\frac1{p(p-1)}}
\quad\text{and}\quad
\E\|g'-\E g'\|_q^q
\le C_p G^q T_{\max}^{\frac{p-2}{(p-1)^2}},
\]
where $g_*=\lambda J_p(y-x_*)=\nabla f_{p,\lambda}(y)$.
Moreover, the expected number of stochastic-gradient evaluations
required to compute $g'$ is $O(1+\log T_{\max})$.
\end{lemma}

\begin{proof}
Since $\psi(y)=0$ and $x_*$ is a minimizer over $\mathcal X$,
\begin{flalign*}
    &f(x_*)+\frac{\lambda}{p}\|x_*-y\|_p^p\le f(y).
\end{flalign*}
By Jensen's inequality, it holds that
\[
\|\mathbb E\widehat{\nabla}f(x)\|_q^q
\le
\mathbb E\|\widehat{\nabla}f(x)\|_q^q
\le G^q.
\]
Since $\mathbb E\widehat{\nabla}f(x)\in\partial f(x)$,
this implies that $f$ is $G$-Lipschitz with respect to
$\|\cdot\|_p$ on $\mathcal X$.
Moreover, by convexity and H\"older's inequality,
\begin{flalign*}
    &f(y)-f(x_*)
    \le
    \left\langle
    \mathbb E\widehat{\nabla}f(y),y-x_*
    \right\rangle
    \le
    \|\mathbb E\widehat{\nabla}f(y)\|_q\|y-x_*\|_p
    \le G\|y-x_*\|_p.
\end{flalign*}
Therefore,
\[
\frac{\lambda}{p}\|x_*-y\|_p^p
\le G\|x_*-y\|_p,
\]
which implies that
\[
\|y-x_*\|_p
\le
p^{1/(p-1)}
\left(\frac{G}{\lambda}\right)^{1/(p-1)}.
\]
Since $F$ is $2^{2-p}\lambda$-uniformly convex of order $p$,
the ODC guarantee gives
\begin{flalign*}
    &\mathbb E\|x_j-x_*\|_p^p
    \le
    c\left(\frac{G}{2^{2-p}\lambda}\right)^q
    2^{-\frac{j}{p-1}}
    =
    c\,2^{(p-2)q}
    \left(\frac{G}{\lambda}\right)^q
    2^{-\frac{j}{p-1}}.
\end{flalign*}
By Lyapunov's inequality, for $0<s\le p$, it holds that
\begin{flalign*}
    &\mathbb E\|x_j-x_*\|_p^s
    \le
    \left(\mathbb E\|x_j-x_*\|_p^p\right)^{s/p}
    \le
    C_p
    \left(\frac{G}{\lambda}\right)^{\frac{s}{p-1}}
    2^{-\frac{js}{p(p-1)}}.
\end{flalign*}
Moreover, it holds that
\begin{flalign*}
    &\|J_p(u)-J_p(v)\|_q
    \le
    (p-1)(\|u\|_p+\|v\|_p)^{p-2}\|u-v\|_p
\end{flalign*}
for all $u,v\in\mathbb R^d$.
We now apply this inequality with $u=y-x_j$ and $v=y-x_*$.
Then,
\begin{flalign*}
    \|y-x_*\|_p+\|y-x_j\|_p &\le
    2\|y-x_*\|_p+\|x_j-x_*\|_p \\
    &\le 2p^{1/(p-1)}\left(\frac{G}{\lambda}\right)^{1/(p-1)}
    +\|x_j-x_*\|_p.
\end{flalign*}
Using that $(s+t)^{p-2}\le C_p(s^{p-2}+t^{p-2})$
for nonnegative $s,t$ and some constant $C_p$, we have
\begin{flalign}
    \label{ineq:grad_diff}
    &\|g_j-g_*\|_q
    \le
    C_p\lambda
    \left[
    \left(\frac{G}{\lambda}\right)^{\frac{p-2}{p-1}}
    \|x_j-x_*\|_p
    +\|x_j-x_*\|_p^{p-1}
    \right].
\end{flalign}
Taking expectations in \eqref{ineq:grad_diff} and using
the Lyapunov bound with $s=1$ and $s=p-1$ yields
\begin{flalign}
    \label{eq:usefulbound}
    \|\mathbb E g_j-g_*\|_q
    &\le \mathbb E\|g_j-g_*\|_q \le
    C_pG
    \left(
    2^{-\frac{j}{p(p-1)}}+2^{-\frac{j}{p}}
    \right)
    \le
    C_pG\,2^{-\frac{j}{p(p-1)}},
\end{flalign}
where we use
$2^{-j/p}\le2^{-j/[p(p-1)]}$ for $j\ge0$ and $p>2$.
Moreover, it holds that
\begin{flalign*}
    &\|g_j-g_*\|_q^q
    \le
    C_p\lambda^q
    \left[
    \left(\frac{G}{\lambda}\right)^{\frac{(p-2)q}{p-1}}
    \|x_j-x_*\|_p^q
    +\|x_j-x_*\|_p^p
    \right],
\end{flalign*}
where we used $(p-1)q=p$.
Absorbing constants into $C_p$ and taking expectations yields
\begin{flalign*}
    &\mathbb E\|g_j-g_*\|_q^q
    \le
    C_pG^q
    \left(
    2^{-\frac{j}{(p-1)^2}}
    +2^{-\frac{j}{p-1}}
    \right)
    \le
    C_pG^q\,2^{-\frac{j}{(p-1)^2}},
\end{flalign*}
because $p>2$ implies
$2^{-j/(p-1)}\le2^{-j/(p-1)^2}$ for $j\ge0$.
Therefore, for $j\ge1$,
\begin{flalign*}
    \mathbb E\|g_j-g_{j-1}\|_q^q
    &\le
    2^{q-1}
    \left(
    \mathbb E\|g_j-g_*\|_q^q
    +\mathbb E\|g_{j-1}-g_*\|_q^q
    \right)\le
    C_pG^q\,2^{-\frac{j}{(p-1)^2}}.
\end{flalign*}
On the other hand,
\[
\mathbb E g'
=
\mathbb E g_0
+\sum_{j=1}^{j_{\max}}
2^{-j}2^j\mathbb E(g_j-g_{j-1})
=
\mathbb E g_{j_{\max}}.
\]
Since $2^{j_{\max}}>T_{\max}/2$, combining this with
\eqref{eq:usefulbound} implies that
\begin{flalign*}
    &\|\mathbb E g'-g_*\|_q
    \le
    C_pG T_{\max}^{-\frac1{p(p-1)}}.
\end{flalign*}
Since
\[g'-g_*=(g_0-g_*)+\boldsymbol{1}[J\le j_{\max}]\,2^J(g_J-g_{J-1}),\]
it holds that
\begin{flalign*}
    \mathbb E\|g'-g_*\|_q^q &\le
    2^{q-1}\mathbb E\|g_0-g_*\|_q^q+
    2^{q-1}\mathbb E \left\|
    \boldsymbol{1}[J\le j_{\max}]\,2^J(g_J-g_{J-1})\right\|_q^q\\
    &\le
    2^{q-1}\mathbb E\|g_0-g_*\|_q^q
    +2^{q-1}\sum_{j=1}^{j_{\max}}
    2^{-j}2^{qj}\mathbb E\|g_j-g_{j-1}\|_q^q\\
    &\le C_pG^q\left(1+\sum_{j=1}^{j_{\max}}
    2^{\left(\frac1{p-1}-\frac1{(p-1)^2}\right)j}\right)=C_pG^q
    \left(1+\sum_{j=1}^{j_{\max}}
    2^{\frac{(p-2)j}{(p-1)^2}}\right)
\end{flalign*}
where the second inequality follows by the law of total
expectation and the independence of $J$ from the oracle randomness.
Since $p>2$, the geometric sum gives
\[\mathbb E\|g'-g_*\|_q^q
\le C_pG^qT_{\max}^{\frac{p-2}{(p-1)^2}}.\]
Therefore,
\begin{flalign*}
    \mathbb E\|g'-\mathbb E g'\|_q^q
    &\le
    2^{q-1}\mathbb E\|g'-g_*\|_q^q
    +2^{q-1}\|\mathbb E g'-g_*\|_q^q\le
    C_pG^q
    \left(
    T_{\max}^{\frac{p-2}{(p-1)^2}}
    +T_{\max}^{-\frac{q}{p(p-1)}}
    \right)\\
    &\le
    C_pG^qT_{\max}^{\frac{p-2}{(p-1)^2}}.
\end{flalign*}
Moreover, the expected number of stochastic-gradient evaluations
is at most
\begin{flalign*}
    &1+\sum_{j=1}^{j_{\max}}
    2^{-j}(2^j+2^{j-1})
    =
    1+\frac32j_{\max}
    =
    O(1+\log T_{\max}),
\end{flalign*}
since level $j$ is drawn with probability $2^{-j}$ and requires
at most $2^j+2^{j-1}$ oracle calls, while the base term $g_0$
requires at most one oracle call.
This finishes the proof.
\end{proof}

\begin{theorem}\label{thm:MorGradEst}
Let $f$ and $\widehat{\nabla}f$ satisfy Assumption~1, let
$y\in\mathcal X$, and let $\lambda,\delta,\sigma>0$ and $p > 2$.
The procedure
\[
\ell_p\text{-}\mathrm{MorGradEst}
\bigl(\widehat{\nabla}f,y,\lambda,\delta,\sigma^q,\mathcal X\bigr)
\]
outputs an estimator $\widehat{\nabla}f_{p,\lambda}(y)$ satisfying
\[\bigl\|\E\,\widehat{\nabla}f_{p,\lambda}(y)-\nabla f_{p,\lambda}(y)\bigr\|_q
\le\delta\quad\text{and}\quad\E\bigl\|
\widehat{\nabla}f_{p,\lambda}(y)
-\nabla f_{p,\lambda}(y)\bigr\|_q^q\le\sigma^q,\]
using $N_{\widehat{\nabla}f}$ stochastic-gradient computations, where
\[\E N_{\widehat{\nabla}f}=O\!\left(
\left[1+\frac{G^p}{\sigma^p}\left(
\max\left\{1,\frac{G}{\min\{\delta,\sigma\}}\right\}\right)^{p(p-2)}
\right]\left[1+\log\max\left\{1,\frac{G}{\min\{\delta,\sigma\}}\right\}\right]
\right).\]
\end{theorem}

\begin{proof}
Write the algorithm's output as
$\widehat{\nabla}f_{p,\lambda}(y)=
\frac1N\sum_{i=1}^N g'^{(i)}$,
with $g'^{(i)}$ i.i.d.\ copies of the MLMC estimator.
Note that from Lemma~\ref{lem:morgrad_mlmc}, we have
\begin{flalign}
    &\|\mathbb E g'-\nabla f_{p,\lambda}(y)\|_q
    \le
    C_1G T_{\max}^{-\frac1{p(p-1)}}\quad,
    \quad\mathbb E\|g'-\mathbb E g'\|_q^q
    \le
    C_2G^q T_{\max}^{\frac{p-2}{(p-1)^2}}.
    \label{eq:C1defn}
\end{flalign}
Since $T_{\max}\ge\left(\frac{2C_1G}{\min\{\delta,\sigma\}}\right)^{p(p-1)}$,
it holds that
\begin{flalign*}
    &\|\mathbb E g'-\nabla f_{p,\lambda}(y)\|_q
    \le
    \frac{\min\{\delta,\sigma\}}2.
\end{flalign*}
Since
$\mathbb E\widehat{\nabla}f_{p,\lambda}(y)=\mathbb E g'$,
it holds that
\begin{flalign*}
    &\|
    \mathbb E\widehat{\nabla}f_{p,\lambda}(y)
    -\nabla f_{p,\lambda}(y)
    \|_q
    \le
    \frac{\min\{\delta,\sigma\}}2
    \le\frac{\delta}{2}
    \le\delta.
\end{flalign*}
The random vectors
$g'^{(i)}-\mathbb E g'^{(i)}$ are i.i.d.\ and centered.
Since $1<q\le2$, the von Bahr--Esseen inequality applied
coordinatewise, together with the lemma, gives, for some
constant $C_{\mathrm{vBE}}$,
\begin{flalign*}
    \mathbb E
    \left\|
    \frac1N\sum_{i=1}^N
    \bigl(g'^{(i)}-\mathbb E g'^{(i)}\bigr)
    \right\|_q^q
    &\le
    \frac{C_{\mathrm{vBE}}}{N^{q-1}}
    \mathbb E\|g'-\mathbb E g'\|_q^q\le
    \frac{
    C_{\mathrm{vBE}}C_2G^q
    T_{\max}^{\frac{p-2}{(p-1)^2}}
    }{N^{q-1}}.
\end{flalign*}
Since
\begin{flalign}
    N\ge 2^p(C_{\mathrm{vBE}}C_2)^{p-1}
    \frac{G^p}{\sigma^p}
    T_{\max}^{\frac{p-2}{p-1}}=
    2^pC_3\frac{G^p}{\sigma^p}
    T_{\max}^{\frac{p-2}{p-1}},
    \qquad
    C_3:=(C_{\mathrm{vBE}}C_2)^{p-1},
    \label{eq:C3defn}
\end{flalign}
it holds that
\[\frac{C_{\mathrm{vBE}}C_2G^q
T_{\max}^{\frac{p-2}{(p-1)^2}}}{N^{q-1}}
\le\left(\frac{\sigma}{2}\right)^q.\]
Therefore,
\begin{flalign*}
    &\mathbb E
    \left\|
    \widehat{\nabla}f_{p,\lambda}(y)
    -\mathbb E\widehat{\nabla}f_{p,\lambda}(y)
    \right\|_q^q
    \le
    \frac{\sigma^q}{2^q}.
\end{flalign*}
Thus,
\begin{flalign*}
    \mathbb E\left\|
    \widehat{\nabla}f_{p,\lambda}(y)
    -\nabla f_{p,\lambda}(y)\right\|_q^q&\le2^{q-1}\mathbb E
    \left\|
    \widehat{\nabla}f_{p,\lambda}(y)
    -\mathbb E\widehat{\nabla}f_{p,\lambda}(y)\right\|_q^q+2^{q-1}\left\|
    \mathbb E\widehat{\nabla}f_{p,\lambda}(y)-\nabla f_{p,\lambda}(y)
    \right\|_q^q\\&\le\frac{\sigma^q}{2}
    +2^{q-1}\left(\frac{\min\{\delta,\sigma\}}2\right)^q
    \le\sigma^q
\end{flalign*}
Using $\lceil x\rceil\le1+x$, it holds that
\begin{flalign*}
    T_{\max}&\le1+\left(\frac{2C_1G}{\min\{\delta,\sigma\}}\right)^{p(p-1)}\le
    \left(1+(2C_1)^{p(p-1)}\right)\left(
    \max\left\{1,\frac{G}{\min\{\delta,\sigma\}}\right\}
    \right)^{p(p-1)}\\
    &=C_p\left(\max\left\{1,\frac{G}{\min\{\delta,\sigma\}}\right\}
    \right)^{p(p-1)}
\end{flalign*}
for some constant $C_p\ge1$.
Consequently,
\begin{flalign*}
    1+\log T_{\max}&\le(1+\log C_p)
    +p(p-1)\log\max\left\{1,\frac{G}{\min\{\delta,\sigma\}}
    \right\}\\
    &\le C_p'\left(1+\log\max\left\{1,\frac{G}{\min\{\delta,\sigma\}}
    \right\}\right)
\end{flalign*}
for some constant $C_p'$.
Since $N$ is deterministic, the expected total number of
stochastic-gradient computations satisfies
\begin{flalign*}
    &\mathbb E N_{\widehat{\nabla}f}
    \le
    N\left(1+\frac32j_{\max}\right)
    =
    O\bigl(N(1+\log T_{\max})\bigr).
\end{flalign*}
Moreover,
\begin{flalign*}
    N&\le1+2^pC_3\frac{G^p}{\sigma^p}T_{\max}^{\frac{p-2}{p-1}}\le1+C_p''\frac{G^p}{\sigma^p}\left(\max\left\{1,\frac{G}{\min\{\delta,\sigma\}}\right\}
    \right)^{p(p-2)}
\end{flalign*}
for some constant $C_p''$.
Combining these bounds gives
\[\mathbb EN_{\widehat{\nabla}f}
=O\!\left(\left[1+\frac{G^p}{\sigma^p}\left(\max\left\{1,\frac{G}{\min\{\delta,\sigma\}}\right\}
\right)^{p(p-2)}\right]\left[
1+\log\max\left\{1,\frac{G}{\min\{\delta,\sigma\}}\right\}
\right]\right),\]
which finishes the proof.
\end{proof}

We now provide a concrete \odc algorithm. Our Algorithm \ref{alg:epochda_odc} extends Algorithm $3$ of \cite{nesterov2009primal} to a composite setting. Following \cite{asi2021stochastic}, to show that Algorithm \ref{alg:epochda_odc} is an ODC$_p$ algorithm, we first establish a convergence guarantee for the single-epoch Algorithm \ref{alg:dam_algorithm}. 

Throughout this section, we use the global prox function
\[d_{\bar{x},R}(x):=
\frac{1}{pR^p}\|x-\bar{x}\|_p^p,
\qquad x\in\mathcal{X},\]
where \(R>0\).
Let \(\mu_p>0\) denote the \(p\)-uniform convexity constant of
\(u\mapsto \frac{1}{p}\|u\|_p^p\) with respect to
\(\|\cdot\|_p\).
Consequently, \(d_{\bar{x},R}\) is
\(\mu_p/R^p\)-uniformly convex on all of \(\mathcal{X}\).

\begin{equation}
\label{eq:V}
V_{z,R,\beta,c}(s)= \max_{x\in\mathcal X}\left\{
\langle s,x-z\rangle-c\psi(x)
-\beta d_{z,R}(x)\right\}.
\end{equation}
\begin{equation}
\label{eq:V,Pi}
\pi_{z,R,\beta,c}(s)=
\arg\max_{x\in\mathcal X}\left\{
\langle s,x-z\rangle-c\psi(x)
-\beta d_{z,R}(x)\right\}.
\end{equation}

where \(\psi:\mathcal X\to\mathbb R\) is \(\mu\)-uniformly convex of order \(p\), and \(c,\beta>0\).
We also assume \(z\in\mathcal
X\).

\begin{algorithm}[H]
\caption{
Dual Averaging Method
\(DAM(\hat{\nabla}f,\psi,\mathcal X,
\bar{x},\gamma,R,N)\)
}
\label{alg:dam_algorithm}
\begin{algorithmic}[1]
\State Set $x_0 = \bar{x}, s_0=0 \in E^\star,c_0=0 \text{ and } \beta=\gamma N^{1/q}$
\For{$i=0,1,\dots,N-1$}
    \State $c_{i+1}=c_i+1$
    \State $x_i=
    \pi_{\bar{x},R,\beta,c_{i+1}}(-s_i)$
    \State Query a stochastic gradient $\hat{\nabla}f(x_i)$,
    and set $s_{i+1}=s_i+\hat{\nabla}f(x_i)$
\EndFor
\State \Return
$\bar{x}_N(\bar{x},R)=
\frac1N\sum_{i=0}^{N-1}x_i$
\end{algorithmic}
\end{algorithm}

\begin{algorithm}[H]
\caption{EpochDA$(\hat{\nabla}f,\psi,\mathcal X,x_0,R_0,N)$}
\label{alg:epochda_odc}
\begin{algorithmic}[1]
\State Set
$y_0=x_0,\tau=p-1$
\State For each \(k\ge1\), define
$N_k=\left\lceil\frac{2^p p^p G^p}
{\mu^p\mu_pR_{k-1}^{p(p-1)}}
\right\rceil$. 
Define $m(N)=\max\left\{k\ge0:
\sum_{j=1}^{k}N_j\le N\right\}$.
\For{$k=1,\dots,m(N)$}
    \State 
    $R_k^p=2^{-k}R_0^p$,
    $\gamma_k=\frac{GR_{k-1}}{\mu_p^{1/p}}$.
    \State 
    $y_k=DAM(\hat{\nabla}f,\psi,
    \mathcal X,y_{k-1},\gamma_k,
    R_{k-1},N_k)$
\EndFor
\State \Return
$\hat{x}_N=y_{m(N)}$.

\end{algorithmic}
\end{algorithm}
 
\begin{lemma} (Lemma 1 of \cite{nesterov2009primal}, Lemma~1 of \cite{iouditski2014primal})\label{lem:V,pi} Let $E^\star$ be the dual space of $E$. The following statements hold: 
\begin{itemize}
    \item If $\beta_1 \leq \beta_2$, then \[
    V_{z,R,\beta_1,c}(s) \geq V_{z,R,\beta_2,c}(s)\]
    \item The function \(V_{z,R,\beta,c}\) is convex and
    differentiable on \(E^\star\).
    Moreover,
    \[\left\|V'_{z,R,\beta,c}(s_1)
    -V'_{z,R,\beta,c}(s_2)\right\|_p
    \le\left(\frac{R^p}{\beta\mu_p}
    \right)^{q-1}\|s_1-s_2\|_q^{q-1},\qquad \forall s_1,s_2\in E^\star.\]
    \item For any $s\in E^*$,
    \[V'_{z,R,\beta,c}(s)+z
    =\pi_{z,R,\beta,c}(s)\in \mathcal{X} \]
\end{itemize}
\end{lemma}
One might notice that $V_{z,R,\beta,c}(s), \pi_{z,R,\beta,c}(s)$ are defined without the term $\psi(x)$ in $\cite{iouditski2014primal}$. However, the same result will be achieved by completely following the proof in $\cite{nesterov2009primal}$. 

\begin{lemma}
\label{lem:epoch convergence}
Let $f: \mathcal{X} \rightarrow \R$ and $\hat{\nabla}f$ satisfy \Cref{Assumption}. For any $k \geq 1, N \geq 1$ and $x \in \mathcal{X} $, \Cref{alg:dam_algorithm} guarantees that
\begin{align}
\sum_{i=0}^{N-1}\Big[\langle
\hat{\nabla}f(x_i),x_i-x\rangle
+\psi(x_i)-\psi(x)\Big]\le
\beta d_{\bar{x},R}(x)+\frac{R^q}
{q\mu_p^{q-1}\beta^{q-1}}
\sum_{i=0}^{N-1}
\|\hat{\nabla}f(x_i)\|_q^q.
\label{eq:epoch-convergence}
\end{align}
\end{lemma}
\begin{proof}
Because $d_{\bar{x},R}(x)=
\frac{1}{pR^p}\|x-\bar{x}\|_p^p$
is \(\mu_p/R^p\)-uniformly convex,
$c\psi(x)+\beta d_{\bar{x},R}(x)$
is at least \(\beta\mu_p/R^p\)-uniformly convex.
For each \(i=0,\ldots,N-1\), since
$s_{i+1}=s_i+\hat{\nabla}f(x_i)$,
the standard smoothness inequality for \(V_{\bar{x},R,\beta,c_{i+1}}\) gives
\begin{align}
V_{\bar{x},R,\beta,c_{i+1}}(-s_{i+1}) &\le V_{\bar{x},R,\beta,c_{i+1}}(-s_i)
-\left\langle \hat{\nabla}f(x_i),
V'_{\bar{x},R,\beta,c_{i+1}}(-s_i)
\right\rangle+\frac{R^q}
{q\mu_p^{q-1}\beta^{q-1}}
\|\hat{\nabla}f(x_i)\|_q^q.
\label{eq:V-step}
\end{align}
By the definition of \(x_i\),
$V'_{\bar{x},R,\beta,c_{i+1}}(-s_i)+\bar{x}=x_i$.
Therefore,
\begin{align}
V_{\bar{x},R,\beta,c_{i+1}}(-s_{i+1}) &\le V_{\bar{x},R,\beta,c_{i+1}}(-s_i)
-\langle\hat{\nabla}f(x_i),
x_i-\bar{x}\rangle +\frac{R^q}
{q\mu_p^{q-1}\beta^{q-1}}
\|\hat{\nabla}f(x_i)\|_q^q.
\end{align}
Since \(c_{i+1}=c_i+1\), and \(x_i\) maximizes
the objective defining
\(V_{\bar{x},R,\beta,c_{i+1}}(-s_i)\),
we also have
\[V_{\bar{x},R,\beta,c_{i+1}}(-s_i)
\le V_{\bar{x},R,\beta,c_i}(-s_i)
-\psi(x_i).\]
Combining the previous two inequalities yields
\begin{align}
V_{\bar{x},R,\beta,c_{i+1}}(-s_{i+1})&\le
V_{\bar{x},R,\beta,c_i}(-s_i)
-\langle\hat{\nabla}f(x_i),
x_i-\bar{x}\rangle-\psi(x_i)+
\frac{R^q}{q\mu_p^{q-1}
\beta^{q-1}}
\|\hat{\nabla}f(x_i)\|_q^q
\end{align}
Summing over \(i=0,\ldots,N-1\), we obtain
\begin{align}
\sum_{i=0}^{N-1}\left[\langle
\hat{\nabla}f(x_i),x_i-\bar{x}
\rangle+\psi(x_i)\right]&\le
V_{\bar{x},R,\beta,c_0}(-s_0)
-V_{\bar{x},R,\beta,c_N}(-s_N)
+\frac{R^q}
{q\mu_p^{q-1}\beta^{q-1}}
\sum_{i=0}^{N-1}
\|\hat{\nabla}f(x_i)\|_q^q.
\label{eq:telescope-V}
\end{align}
Since $s_0=0,c_0=0$ and
\[d_{\bar{x},R}(x)\ge0,
\qquad d_{\bar{x},R}(\bar{x})=0,
\qquad \bar{x}\in\mathcal X,\]
we have
\[V_{\bar{x},R,\beta,c_0}(-s_0)
=\max_{x\in\mathcal X}
\{-\beta d_{\bar{x},R}(x)\}=0\]
For any fixed \(x\in\mathcal X\),
the definition of \(V\) gives
\begin{align}
-V_{\bar{x},R,\beta,c_N}(-s_N)
&\le \langle s_N,x-\bar{x}\rangle
+N\psi(x)+\beta d_{\bar{x},R}(x)
\nonumber\\ &=\sum_{i=0}^{N-1}
\langle\hat{\nabla}f(x_i),
x-\bar{x}\rangle+N\psi(x)+
\beta d_{\bar{x},R}(x).
\end{align}
Substituting this into
\eqref{eq:telescope-V}
and rearranging proves
\eqref{eq:epoch-convergence}.
\end{proof}

\begin{lemma}
\label{lem:single-epoch-expected}
Let $F=f+\psi$, and let \(x^\star\in\arg\min_{x\in\mathcal X}F(x)\).
Assume that the stochastic oracle satisfies
\[\mathbb E[\hat{\nabla}f(x_i)
\mid\mathcal F_i]=\nabla f(x_i)
\quad \text{and} \quad
\mathbb E[\|\hat{\nabla}f(x_i)\|_q^q
\mid\mathcal F_i]\le G^q,\]
where \(\mathcal F_i\) denotes the history before the fresh oracle call at iteration \(i\).
Suppose that \(\bar{x}\) satisfies
$\mathbb E\|\bar{x}-x^\star\|_p^p
\le R^p$. Then the output of
\Cref{alg:dam_algorithm} satisfies
\begin{align}
\mathbb E[F(\bar{x}_N(\bar{x},R))-F(x^\star)] &\le \frac{1}{N^{1/p}}\left[\frac{\gamma}{p}+
\frac{R^qG^q}{q\gamma^{q-1}
\mu_p^{q-1}}\right],
\label{eq:single-epoch-function}
\\ \mathbb E
\|\bar{x}_N(\bar{x},R)-x^\star\|_p^p
&\le \frac{p}{\mu N^{1/p}}\left[
\frac{\gamma}{p}+\frac{R^qG^q}
{q\gamma^{q-1}\mu_p^{q-1}}\right].
\label{eq:single-epoch-distance}
\end{align}
\end{lemma}
\begin{proof}
Apply \Cref{lem:epoch convergence} with \(x=x^\star\).
Since $d_{\bar{x},R}(x^\star)=
\frac{\|\bar{x}-x^\star\|_p^p
}{pR^p}$, we have
$\mathbb E[d_{\bar{x},R}(x^\star)]
\le\frac1p$. Moreover, conditional unbiasedness of the stochastic oracle gives
\begin{align}
\mathbb E \left[\left\langle
\hat{\nabla}f(x_i),x_i-x^\star
\right\rangle\right]=\mathbb E
\left[\left\langle\nabla f(x_i),
x_i-x^\star\right\rangle\right]
\ge \mathbb E[f(x_i)-f(x^\star)],
\end{align}
where the last inequality follows from convexity of \(f\). Therefore,
\begin{align}
\mathbb E \Big[\langle \hat{\nabla}f(x_i),x_i-x^\star
\rangle+\psi(x_i)-\psi(x^\star)
\Big]
\ge\mathbb E[F(x_i)-F(x^\star)].
\end{align}
Taking expectations in
\eqref{eq:epoch-convergence},
using $\beta=\gamma N^{1/q}$,
and the moment bound $\mathbb E \|\hat{\nabla}f(x_i)\|_q^q
\le G^q$, gives
\begin{align}\label{eq:epoch-da convergence}
\frac1N \sum_{i=0}^{N-1}\mathbb E[
F(x_i)-F(x^\star)]&\le \frac{1}{N^{1/p}}\left[\frac{\gamma}{p}+
\frac{R^qG^q}
{q\gamma^{q-1}\mu_p^{q-1}}\right].
\end{align}
By convexity of \(F\),
\[F(\bar{x}_N(\bar{x},R)) = F\left(
\frac1N\sum_{i=0}^{N-1}x_i\right)
\le\frac1N\sum_{i=0}^{N-1}F(x_i).\]
Plugging this into \eqref{eq:epoch-da convergence} proves
\eqref{eq:single-epoch-function}.
Finally, since \(F=f+\psi\) is
\(\mu\)-uniformly convex of order \(p\), \[F(x)-F(x^\star)\ge
\frac{\mu}{p}\|x-x^\star\|_p^p.\]
Hence
\[\mathbb E\|\bar{x}_N(\bar{x},R)-x^\star\|_p^p\le \frac{p}{\mu}
\mathbb E[F(\bar{x}_N(\bar{x},R))-F(x^\star)],\] which proves
\eqref{eq:single-epoch-distance}.
\end{proof}

\begin{corollary}
\label{cor:ODC}
Under the assumptions of
\Cref{lem:single-epoch-expected},
choose $\gamma=\frac{GR}{\mu_p^{1/p}}$. Then
\begin{align}
\mathbb E[F(\bar{x}_N(\bar{x},R))-F(x^\star)]&\le\frac{GR}{\mu_p^{1/p}N^{1/p}},
\label{eq:DAM-function-rate}\\ \mathbb E
\|\bar{x}_N(\bar{x},R)-x^\star\|_p^p &\le \frac{
GRp}{\mu\mu_p^{1/p}N^{1/p}}.
\label{eq:DAM-distance-rate}
\end{align}
\end{corollary}
   
\begin{theorem}
\label{thm: convergence of ODC}
Assume $\|y_0-x^\star\|_p^p
\le R_0^p$. Then the sequence \(\{y_k\}_{k=0}^{m(N)}\)
generated by
\Cref{alg:epochda_odc}
satisfies 
\[\mathbb E \|y_k-x^\star\|_p^p\le R_k^p=
2^{-k}R_0^p, \qquad
k=0,1,\ldots,m(N).\]
Moreover, for every completed epoch \(k\ge1\),
\[\mathbb E[F(y_k)
-F(x^\star)]\le \frac{\mu}{p} R_k^p.\]
\end{theorem}
\begin{proof}
We first prove by induction that
$\mathbb E \|y_k-x^\star\|_p^p\le R_k^p$.
For \(k=0\), this follows from the assumption
$\|y_0-x^\star\|_p^p \le
R_0^p$. Now assume that for some \(k\ge1\),
$\mathbb E \|y_{k-1}-x^\star\|_p^p\le R_{k-1}^p$. By \Cref{cor:ODC}, this gives
\[\mathbb E[F(y_k)-F(x^\star)]\le\frac{
GR_{k-1}}{
\mu_p^{1/p}N_k^{1/p}}.\]
By the definition of \(N_k\), $N_k\ge\frac{2^p p^pG^p}{
\mu^p\mu_pR_{k-1}^{p(p-1)}}.$ Therefore,
$N_k^{1/p}\ge\frac{2pG}{
\mu\mu_p^{1/p}R_{k-1}^{p-1}
}$. Substituting this bound yields
\begin{align}
\mathbb E[F(y_k)-F(x^\star)
]&\le\frac{GR_{k-1}}{\mu_p^{1/p}}\cdot\frac{
\mu\mu_p^{1/p}R_{k-1}^{p-1}
}{2pG}=\frac{\mu}{2p}R_{k-1}^p=\frac{\mu}{p} R_k^p.
\end{align}
Since \(F\) is \(\mu\)-uniformly convex of order \(p\),
\[F(y_k)-F(x^\star)\ge
\frac{\mu}{p}\|y_k-x^\star\|_p^p.\]
Taking expectations gives
\[\mathbb E\|y_k-x^\star\|_p^p\le R_k^p.\]
which completes the induction. Since $\tau=p-1$, then $N_k=\left\lceil
\frac{2^p p^pG^p}{
\mu^p\mu_pR_0^{p(p-1)}}2^{(k-1)\tau}\right\rceil$.
Let $m=m(N)$. By maximality of \(m\),
$N<\sum_{k=1}^{m+1}N_k$. Then,
\begin{align}
N&<\frac{2^p p^pG^p}{
\mu^p\mu_pR_0^{p(p-1)}}\sum_{k=0}^{m}2^{k\tau}+(m+1)\le \frac{2^pp^pG^p}{
\mu^p\mu_pR_0^{p(p-1)}}
\frac{2^{(m+1)\tau}}{2^\tau-1}+(m+1).
\end{align}
Since \(\tau\ge1\) and \(m\ge1\),
$m+1 \le 2^m \le 2^{m\tau}$. Hence
\[N\le\left[\frac{2^p p^pG^p}{
\mu^p\mu_pR_0^{p(p-1)}}\frac{2^\tau}{
2^\tau-1}+1\right]2^{m\tau}.\]
Therefore,
\[2^{-m}\le\left[\frac{2^p p^pG^p}{
\mu^p\mu_pR_0^{p(p-1)}}\frac{2^\tau}{(2^\tau-1)N}+\frac1N\right]^{1/\tau}.\]
From the epoch guarantee,
\[\mathbb E[F(\hat{x}_N)-F(x^\star)]\le
\frac{\mu}{p} R_m^p= \frac{\mu}{p} R_0^p2^{-m}.\]
Therefore,
\[\mathbb{E}[F(\widehat{x}_N)-F(x^\star)]\le \frac{\mu}{p} R_0^p
\left[\frac{2^p p^p G^p}{\mu^p\mu_p^p R_0^{p(p-1)}}\frac{2^\tau}{(2^\tau-1)N}+\frac{1}{N}\right]^{1/\tau}.\]
Using $R_0\le p^{q-1}\left(\frac{G}{\mu}
\right)^{1/(p-1)}$, both terms are bounded by $C_p\frac{G^q}{\mu^{1/(p-1)}N^{1/(p-1)}}$, for a constant \(C_p\) depending only on \(p\).
Thus,
\[\mathbb E[F(\hat{x}_N)-F(x^\star)]\le
C_p\frac{G^q}{\mu^{1/(p-1)}N^{1/(p-1)}}.\]
Finally, uniform convexity gives
\[\mathbb E\|\hat{x}_N-x^\star\|_p^p
\le\frac1\mu\mathbb E[F(\hat{x}_N)-F(x^\star)],\]and hence 
\[\mathbb E\|\hat{x}_N-x^\star\|_p^p
\le C_p\frac{G^q}{\mu^qN^{1/(p-1)}}.\]
\end{proof}

\begin{corollary}
\Cref{alg:epochda_odc} is an
\(\mathrm{ODC}_p\) algorithm with a constant \(C_p\) depending only on \(p\).
\end{corollary}
\subsection{Projection-Efficient $\ell_p$-Lipschitz Optimization}
\label{subsection:projection-efficient}
We combine our $\ell_p$-Moreau envelope gradient estimator
with an accelerated gradient method to reduce the number
of expensive constrained updates over $\mathcal X$.
We assume that $f:\R^d\to\R$ and its stochastic subgradient
oracle satisfy Assumption~\ref{Assumption} on $\R^d$,
and that the initial point $x_0\in\mathcal X$ satisfies
$\|x_0-x^\star\|_p\le R$.
The Moreau envelope is defined over $\R^d$, and its gradient
is estimated using $ODC_p$ over $\ell_p$ balls. Only the outer updates involve the complicated set
$\mathcal X$.
\newline \newline 
Our approach is inspired by that of Algorithm 3 in \cite{asi2021stochastic}, which involves applying a variant of Nesterov’s accelerated gradient descent method on the Moreau envelope
$f_{\lambda}(y) = \min_x\{f(x) + \frac{\lambda}{2}\|x-y\|_2^2\}$ in the Euclidean setting. Importantly, it may be seen that projections onto \(\mathcal X\) are not needed for computing the Euclidean Moreau envelope itself, so that accurate approximations of \(\nabla f_\lambda\) allow the
outer accelerated method to find an \(O(\epsilon)\)-suboptimal point of $f_{\lambda}$ using $T = O\left(\sqrt{\lambda R^2/\epsilon}\right)$
projections onto \(\mathcal X\). To ensure that this point is also \(O(\epsilon)\)-suboptimal
for the original objective \(f\), one must choose $\lambda = G^2/\epsilon$,
which yields the projection complexity $O(GR/\epsilon)$.

\begin{algorithm}[H]
\caption{Stochastic AGD with $\ell_p$ Moreau envelope}
\label{alg:proximal_point}
\begin{algorithmic}[1]
\State Parameters $\lambda,M>0$ and $T\ge1$
\State $r_\lambda=2(G/\lambda)^{1/(p-1)}$
\State $A_0=0$, $x_0=z_0=y_0=x_{\mathrm{init}}$
\State Bias tolerances $\{\delta_t\}_{t=1}^T$ and
moment bounds $\{\sigma_t^q\}_{t=1}^T$
\For{$t=0,\ldots,T-1$}
    \State $A_{t+1}=2^{p+q-1}(p-1)a_{t+1}^q,
    \quad A_{t+1}=A_t+a_{t+1}$
    \State $y_{t+1}\gets
    \dfrac{A_t}{A_{t+1}}x_t
    +\dfrac{a_{t+1}}{A_{t+1}}z_t$
    \State $\begin{aligned}[t]
    \nabla f'_{p,\lambda}(y_{t+1})\gets
    \ell_p\textsc{-MorGradEst}(    &\widehat{\nabla}f,y_{t+1},\lambda,\delta_{t+1},\sigma_{t+1}^q,\mathbb B_{r_\lambda}(y_{t+1}))
    \end{aligned}$
    \State $x_{t+1}\gets    \displaystyle\arg\min_{x\in\mathcal X}
    \left\{
    \langle\nabla f'_{p,\lambda}(y_{t+1}),x-y_{t+1}\rangle
    +\frac{M}{p}\|x-y_{t+1}\|_p^p
    \right\}$
    \State $z_{t+1}\gets    \displaystyle\arg\min_{z\in\mathcal X}
    \left\{
    \frac{a_{t+1}}{M}
    \langle\nabla f'_{p,\lambda}(y_{t+1}),z\rangle
    +\omega_p(z,z_t)
    \right\}$
\EndFor
\State \Return $x_T$
\end{algorithmic}
\end{algorithm}
Throughout this section, we define
\[\psi(u)=\frac1p\|u-x_{\mathrm{init}}\|_p^p,
\qquad \omega_p(u,v)
=\psi(u)-\psi(v)-\langle\nabla\psi(v),u-v\rangle.\]
\begin{lemma}\label{lem:section3convergence}
Let $\mathcal X\subseteq\mathbb B_R(0)$ be a nonempty
closed convex set. Let $f:\R^d\to\R$ and
$\widehat{\nabla}f$ satisfy Assumption~\ref{Assumption}. Let $x_0=x_{\mathrm{init}}\in\mathcal X$ satisfy
$\|x_0-x^\star\|_p\le R$ for some $x^\star\in\mathcal X$.
For any $\epsilon>0$ and integer $T\ge1$, choose
$\lambda,M>0$ such that
\[L_\lambda=(p-1)\lambda^{q-1}G^{2-q}, \qquad M\ge
2L_\lambda^{p/2}\left(\frac{T}{\epsilon}\right)^{(p-2)/2}.\]
For $t=0,\ldots,T-1$, set
\begin{equation}\label{eq:accuracy-parameters}
\delta_{t+1}=\frac{\epsilon}{2R},
\qquad
\sigma_{t+1}^q=
\frac{
q\epsilon M^{1/(p-1)}a_{t+1}
}{
2^{q-1}A_{t+1}
}.
\end{equation}
Then the iterates of \Cref{alg:proximal_point} satisfy,
for $t=0,\ldots,T-1$,
\begin{equation}\label{eq:one-step-convergence}
\begin{aligned}
&\frac{A_{t+1}}{M}
\E\!\left[
f_{p,\lambda}(x_{t+1})-f_{p,\lambda}(x^\star)
\right]
+\E\!\left[\omega_p(x^\star,z_{t+1})\right]
\\&\qquad\le\frac{A_t}{M}
\E\!\left[
f_{p,\lambda}(x_t)-f_{p,\lambda}(x^\star)
\right]
+\E\!\left[\omega_p(x^\star,z_t)\right]
+\frac{3a_{t+1}}{M}\epsilon.
\end{aligned}
\end{equation}
\end{lemma}

\begin{lemma}[{\cite{zalinescu1983uniformly}}]\label{lem:uniform-convexity}
For any $x,y\in\R^n$ and $p\ge2$, 
$\omega_p(y,x)\ge\frac1{p2^{p+1}}\|y-x\|_p^p$.
\end{lemma}

\begin{lemma}\cite{adil2026convex}\label{lem:A_T}
With $A_0=0$, suppose that $a_{t+1}>0$,
$A_{t+1}=A_t+a_{t+1}$, and
$A_{t+1}=2^{p+q-1}(p-1)a_{t+1}^q$.
Then, for all $t\ge0$,
\[
A_{t+1}^{1/p}\ge A_t^{1/p}
+\frac1{p[2^{p+q-1}(p-1)]^{1/q}}.
\]
Consequently, for every $T\ge1$,
\[
A_T^{1/p}\ge A_0^{1/p}
+\frac{T}{p[2^{p+q-1}(p-1)]^{1/q}}.
\]
\end{lemma}

\subsection{\texorpdfstring{Proof of \Cref{lem:section3convergence}}{Proof of the convergence lemma}}
\begin{proof}
The proof is generally inspired by those of \cite{adil2026convex,asi2021stochastic}.
First, $x_t,z_t,y_{t+1}\in\mathcal X$ by construction
and the convexity of $\mathcal X$.
By the convexity of $f_{p,\lambda}$ and the identity
$a_{t+1}(y_{t+1}-z_t)=A_t(x_t-y_{t+1})$, we have
\begin{equation}\label{eq:convexity-coupling}
\begin{aligned}
\frac{A_{t+1}}{M}
\bigl(f_{p,\lambda}(y_{t+1})-f_{p,\lambda}(x^\star)\bigr)\le
\frac{A_t}{M}
\bigl(f_{p,\lambda}(x_t)-f_{p,\lambda}(x^\star)\bigr)
+\frac{a_{t+1}}{M}
\left\langle
\nabla f_{p,\lambda}(y_{t+1}),z_t-x^\star
\right\rangle.
\end{aligned}
\end{equation}
The optimality condition for $z_{t+1}$ and the
three-point identity give
\begin{equation}\label{eq:mirror-bound}
\begin{aligned}
\frac{a_{t+1}}{M}
\left\langle
\nabla f'_{p,\lambda}(y_{t+1}),z_{t+1}-x^\star
\right\rangle\le
\omega_p(x^\star,z_t)
-\omega_p(x^\star,z_{t+1})
-\omega_p(z_{t+1},z_t).
\end{aligned}
\end{equation}
Proximal optimality bounds the minimizing displacement by
$r_\lambda/2$, so the smoothness argument of
\cite{diakonikolas2024complementary} gives $L_\lambda$.
By the radius bound in the proof of Lemma~3 with $\mathcal X=\R^d$,
the ball in \Cref{alg:proximal_point} contains the unconstrained
proximal minimizer. Hence,
\begin{equation}\label{eq:section3-p-model}
\begin{aligned}
f_{p,\lambda}(u)
&\le f_{p,\lambda}(v)
+\langle\nabla f_{p,\lambda}(v),u-v\rangle
+\frac{L_\lambda}{2}\|u-v\|_p^2\\
&\le f_{p,\lambda}(v)
+\langle\nabla f_{p,\lambda}(v),u-v\rangle
+\frac{M}{2p}\|u-v\|_p^p+\frac{\epsilon}{T}.
\end{aligned}
\end{equation}
On the other hand, by Young's inequality we have,
\begin{align*}
&\left\langle
\nabla f_{p,\lambda}(y_{t+1})
-\nabla f'_{p,\lambda}(y_{t+1}),
x_{t+1}-y_{t+1}
\right\rangle\\
&\qquad\le
\frac{M}{2p}\|x_{t+1}-y_{t+1}\|_p^p
+\frac{2^{q-1}}{qM^{q-1}}
\left\|
\nabla f_{p,\lambda}(y_{t+1})
-\nabla f'_{p,\lambda}(y_{t+1})
\right\|_q^q.
\end{align*}
By the convexity of $\mathcal X$, \[\frac{A_t}{A_{t+1}}x_t
+\frac{a_{t+1}}{A_{t+1}}z_{t+1}\in\mathcal X\]
With the optimality of
$x_{t+1}$ and \eqref{eq:section3-p-model},
we obtain
\begin{equation}\label{eq:section3-primal-progress}
\begin{aligned}
\frac{A_{t+1}}{M}
\bigl(f_{p,\lambda}(x_{t+1})
-f_{p,\lambda}(y_{t+1})\bigr)
&\le
\frac{a_{t+1}}{M}
\left\langle
\nabla f'_{p,\lambda}(y_{t+1}),z_{t+1}-z_t
\right\rangle
+\frac{a_{t+1}^p}{pA_{t+1}^{p-1}}
\|z_{t+1}-z_t\|_p^p\\
&\qquad\quad+
\frac{A_{t+1}\epsilon}{MT}
+\frac{2^{q-1}A_{t+1}}{qM^q}
\left\|
\nabla f_{p,\lambda}(y_{t+1})
-\nabla f'_{p,\lambda}(y_{t+1})
\right\|_q^q.
\end{aligned}
\end{equation}
By Lemma~\ref{lem:uniform-convexity},
\[
\omega_p(z_{t+1},z_t)
\ge \frac{1}{p\,2^{p+1}}\|z_{t+1}-z_t\|_p^p.
\]
Using $A_{t+1}=2^{p+q-1}(p-1)a_{t+1}^q$ and $q(p-1)=p$, we obtain
\[\begin{aligned}\frac{a_{t+1}^p}{pA_{t+1}^{p-1}}=\frac{1}{p[2^{p+q-1}(p-1)]^{p-1}}=\frac{1}{p(p-1)^{p-1}2^{p(p-1)+1}}\le\frac{1}{p\,2^{p+1}},
\end{aligned}\]
where the last inequality follows from $p\ge 2$, which gives
$(p-1)^{p-1}\ge 1$ and $p(p-1)+1\ge p+1$. Consequently,
\[
\frac{a_{t+1}^p}{pA_{t+1}^{p-1}}\|z_{t+1}-z_t\|_p^p
\le \omega_p(z_{t+1},z_t).
\]
Combining \eqref{eq:convexity-coupling},
\eqref{eq:mirror-bound}, and
\eqref{eq:section3-primal-progress}, the terms involving
$\langle\nabla f'_{p,\lambda}(y_{t+1}),z_{t+1}-z_t\rangle$
cancel.
Taking conditional expectations, applying Theorem~\ref{thm:MorGradEst}, and using
$\|z_t-x^\star\|_p\le2R$, we obtain
\begin{equation}\label{eq:potential-before-parameters}
\begin{aligned}
\frac{A_{t+1}}{M}
\E\!\left[
f_{p,\lambda}(x_{t+1})-f_{p,\lambda}(x^\star)
\right]
+\E[\omega_p(x^\star,z_{t+1})]
&\le
\frac{A_t}{M}
\E\!\left[
f_{p,\lambda}(x_t)-f_{p,\lambda}(x^\star)
\right]
+\E[\omega_p(x^\star,z_t)]\\
&\qquad\quad+
\frac{2Ra_{t+1}}{M}\delta_{t+1}
+\frac{A_{t+1}\epsilon}{MT}
+\frac{2^{q-1}A_{t+1}}{qM^q}\sigma_{t+1}^q.
\end{aligned}
\end{equation}
To bound the last three terms, observe that
\[A_k^{1/p}-A_{k-1}^{1/p}
\le\frac{a_k}{A_k^{1/q}}
=\frac{1}{[2^{p+q-1}(p-1)]^{1/q}}.\]
Since $A_0=0$, summing this inequality gives
\[\frac{A_{t+1}}{a_{t+1}}
=[2^{p+q-1}(p-1)]^{1/q}A_{t+1}^{1/p}
\le t+1\le T.\]
Substituting \eqref{eq:accuracy-parameters}, we therefore have
\begin{align*}
\frac{2Ra_{t+1}}{M}\delta_{t+1}
+\frac{A_{t+1}\epsilon}{MT}
+\frac{2^{q-1}A_{t+1}}{qM^q}\sigma_{t+1}^q=
\frac{2a_{t+1}\epsilon}{M}
+\frac{A_{t+1}\epsilon}{MT}
\le\frac{3a_{t+1}\epsilon}{M}.
\end{align*}
Combining this with \eqref{eq:potential-before-parameters}
proves \eqref{eq:one-step-convergence}.
\end{proof}

\begin{theorem}\label{thm:section3}
Under the conditions of \Cref{lem:section3convergence},
let $x^\star\in\operatorname*{argmin}_{x\in\mathcal X}f(x)$.
For any $0<\varepsilon\le GR$, with the parameter choices
specified below, \Cref{alg:proximal_point} computes
$x_T\in\mathcal X$ such that
$\E[f(x_T)]\le f(x^\star)+\varepsilon$, using
\[
T=O\!\left(
\left(\frac{GR}{\varepsilon}\right)^{2p/(p+2)}
\right)
\]
outer iterations and $N_{\widehat{\nabla}f}$ stochastic
subgradient computations, where
\[
\E[N_{\widehat{\nabla}f}]
=
O\!\left(
\left(\frac{GR}{\varepsilon}\right)^{p(p-1)}
\left(1+\log\frac{GR}{\varepsilon}\right)
\right).
\]
\end{theorem}

\begin{proof}
Set $\epsilon=\varepsilon/5$.
Since $A_0=0$, summing
\Cref{lem:section3convergence} from $t=0$ to $T-1$ gives
\begin{align*}
&\frac{A_T}{M}\E\!\left[
f_{p,\lambda}(x_T)-f_{p,\lambda}(x^\star)\right]+\E[\omega_p(x^\star,z_T)]\le\omega_p(x^\star,z_0)
+\sum_{t=0}^{T-1}\frac{3a_{t+1}}{M}\epsilon=\omega_p(x^\star,z_0)+\frac{3A_T}{M}\epsilon.
\end{align*}
Using $\omega_p(x^\star,z_T)\ge0$, we obtain
\[\E\!\left[f_{p,\lambda}(x_T)-f_{p,\lambda}(x^\star)\right]
\le\frac{M}{A_T}\omega_p
(x^\star,z_0)+3\epsilon.\]
By the $G$-Lipschitz property of $f$, for any $y\in\R^d$,
\begin{align*}
0\le f(y)-f_{p,\lambda}(y)\le
\sup_{r\ge0}\left\{Gr-\frac{\lambda}{p}r^p\right\}
=\frac{G^q}{q\lambda^{q-1}}.
\end{align*}
Since $f_{p,\lambda}(x^\star)\le f(x^\star)$,
\[\E[f(x_T)]-f(x^\star)\le
\frac{M}{A_T}
\omega_p(x^\star,z_0)
+\frac{G^q}{q\lambda^{q-1}}
+3\epsilon.\]
Choose $\lambda=\left(\frac{G^q}{q\epsilon}\right)^{p-1}$.
Then
\[\frac{G^q}{q\lambda^{q-1}}
=\epsilon,\qquad L_\lambda
=(p-1)\lambda^{q-1}G^{2-q}
=\frac{p-1}{q}\frac{G^2}{\epsilon}.\]
By \Cref{lem:A_T}, we have
\[A_T\ge\frac{T^p}{
p^p[2^{p+q-1}(p-1)]^{p-1}}.\]
Moreover, since $z_0=x_{\mathrm{init}}$,
\[\omega_p(x^\star,z_0)
=\frac1p\|x^\star-x_{\mathrm{init}}\|_p^p
\le\frac{R^p}{p}.\]
Taking
\[M=2L_\lambda^{p/2}
\left(\frac{T}{\epsilon}\right)^{(p-2)/2},\]
we therefore have
\[\frac{M}{A_T}
\omega_p(x^\star,z_0)\le
\frac{2p^{p-1}[2^{p+q-1}(p-1)]^{p-1}L_\lambda^{p/2}R^p
}{\epsilon^{(p-2)/2}T^{(p+2)/2}
}.\]
Thus it suffices to choose
\[T=\left\lceil\left(
2p^{p-1}[2^{p+q-1}(p-1)]^{p-1}
\right)^{2/(p+2)}
\left(\frac{L_\lambda R^2}{\epsilon}\right)^{p/(p+2)}
\right\rceil.\]
Since $L_\lambda=\Theta(G^2/\epsilon)$, this gives
\[T=\Theta\!\left(
\left(\frac{GR}{\epsilon}\right)^{2p/(p+2)}
\right),\qquad\frac{M}{A_T}\omega_p(x^\star,z_0)\le\epsilon.\]
Consequently,
\[\E[f(x_T)]-f(x^\star)
\le\epsilon+\epsilon+3\epsilon
=5\epsilon=\varepsilon.\]
To finish, we bound the expected total number of stochastic
subgradient oracle calls.
Recall that
\[
\delta_t=\frac{\epsilon}{2R},
\qquad
\sigma_t^q=
\frac{q\epsilon M^{1/(p-1)}a_t}
{2^{q-1}A_t},
\qquad t=1,\ldots,T.
\]
The choices of $T$ and $M$ above imply $
M=\Theta\!\left(\frac{\epsilon T^p}{R^p}\right)$.
Since $p/q=p-1$ and $(q-1)(p-1)=1$,
\[\sigma_t^p=
\frac{q^{p-1}}{2}\epsilon^{p-1}M
\left(\frac{a_t}{A_t}\right)^{p-1}.\]
As shown in the proof of
\Cref{lem:section3convergence}, $A_t/a_t\le t\le T$.
Hence
\[\sigma_t^p\ge\frac{q^{p-1}\epsilon^{p-1}M}{2T^{p-1}}
=\Theta\!\left(\frac{\epsilon^pT}{R^p}\right).\]
Together with $\delta_t=\epsilon/(2R)$ and $T\ge1$,
this implies, uniformly over $t=1,\ldots,T$,
\[\min\{\delta_t,\sigma_t\}
=\Theta(\epsilon/R),\qquad
1+\log\frac{G}{\min\{\delta_t,\sigma_t\}}=
O\!\left(1+\log\frac{GR}{\epsilon}\right).\]
Also,
\begin{align*}
\sum_{t=1}^T\frac{G^p}{\sigma_t^p}&=
\frac{2G^p}{q^{p-1}\epsilon^{p-1}M}\sum_{t=1}^T
\left(\frac{A_t}{a_t}\right)^{p-1}\le\frac{2G^pT^p}{q^{p-1}\epsilon^{p-1}M}=O\!\left(
\left(\frac{GR}{\epsilon}\right)^p\right).
\end{align*}
Thus, by Theorem~\ref{thm:MorGradEst}, the summation of stochastic gradient oracle is
\begin{align*}
\E[N_{\widehat{\nabla}f}]
&=
O\!\left(
T+
\sum_{t=1}^T
\frac{G^p}{\sigma_t^p}
\left(
\frac{G}{\min\{\delta_t,\sigma_t\}}
\right)^{p(p-2)}
\left(
1+\log\frac{G}{\min\{\delta_t,\sigma_t\}}
\right)
\right)\\
&=
O\!\left(
T+
\left(\frac{GR}{\epsilon}\right)^{p(p-1)}
\left(1+\log\frac{GR}{\epsilon}\right)
\right)\\
&=
O\!\left(
\left(\frac{GR}{\varepsilon}\right)^{p(p-1)}
\left(1+\log\frac{GR}{\varepsilon}\right)
\right),
\end{align*}
where the last equality uses $\epsilon=\varepsilon/5$
and the bound on $T$.
\end{proof}

\subsection{Convex Composite Minimization}
\label{subsection:convex_composite}
Consider the convex composite optimization problem
\begin{equation}\label{problem:section 4}
\min_{x \in \mathcal X} \Psi(x) := \Lambda(x) + f(x),
\end{equation}
where $\Lambda$ is $(L,\kappa)$-weakly smooth with respect to
$\|\cdot\|_p$ and $f$ satisfies \Cref{Assumption}.
We assume that there exists
$x^\star \in \arg\min_{x\in\mathcal X}\Psi(x)$ such that
$\|x_0-x^\star\|_p \le R$.
A closely related setting was studied in
\cite{diakonikolas2024complementary}, who considered complementary
composite minimization with $\Lambda$ being $(L,\kappa)$-weakly
smooth and $f$ being uniformly convex.
Their generalized AGD+ framework admits an
$\epsilon$-additive optimality gap in
$O\left(\left(LR^{\kappa}/\epsilon\right)^{
\frac{p}{\kappa p-p+\kappa}}\right)$
iterations when specialized to $\ell_p$ geometry.

In this section, we extend the stochastic composite AGD framework of
\cite{asi2021stochastic} to the $\ell_p$ setting using Bregman
proximal steps and our reduced-bias gradient estimator.
Our method handles a nonsmooth stochastic component $f$
without requiring uniform convexity, while preserving
the accelerated outer-loop structure.
\begin{lemma} 
\label{lemma:composite_minimization_helper}
(Lemma 1 of \cite{diakonikolas2024complementary})
    Suppose $\Lambda$ is $(L, \kappa)$-weakly smooth with respect to $||\cdot||_p$ for some $\kappa \in (1, 2]$. Then, for $\delta > 0$, if $M \ge (\frac{2(p - \kappa)}{p \kappa \delta})^{\frac{p - \kappa}{\kappa}} L^{\frac{p}{\kappa}}$, then for all $x, y \in \mathcal{X}$, it holds that \begin{flalign*}
        &\Lambda(y) \le \Lambda(x) + \langle \nabla \Lambda(x), y - x \rangle + \frac{M}{p} ||y - x||_p^p + \frac{\delta}{2}.
    \end{flalign*}
\end{lemma}

\begin{definition}
Let $D_h(z,y):=h(z)-h(y)-\langle\nabla h(y),z-y\rangle.$
For $F(z):=f(z)+\lambda D_h(z,y)$, define
\begin{equation}\label{eq:inexact-prox}
P_{f,\lambda}^{\epsilon,h}(y)
:=
\text{any }x\in\overline{\mathcal X}
\text{ such that }
\E[F(x)]
\le
\min_{z\in\overline{\mathcal X}}F(z)+\epsilon.
\end{equation}
where $\bar X=X\cap B_p(v_0,R)$
\end{definition}

\begin{algorithm}[H]
\caption{Stochastic composite accelerated gradient descent}
\label{alg:composite_convex_opt_alg_label}
\begin{algorithmic}[1]
\State \textbf{Input:} A problem of the form \eqref{problem:section 4} with
$\Lambda, f, \nabla\Lambda, \widehat{\nabla}f$.
\State \textbf{Parameters:} Step size parameters
$\beta_k=M\gamma_k^{p-1}$ and $\gamma_k=\frac{p}{k+p-1}$,
iteration number $N$, approximation parameters
$\{\epsilon_k,\delta_k,\sigma_k^q\}$, and $x_0=v_0$
where $x_0\in X$ satisfying $\|x_0-x^\ast\|_p\le R$.
\State Set
$h(z)=\frac{2^{p-2}}p\|z-v_0\|_p^p$.

\For{$k=1,2,\cdots,N$}

    \State $y_k=(1-\gamma_k)x_{k-1}+\gamma_k v_{k-1}$

    \State $\bar v_k=
    P^{\epsilon_k,h}_{\bar\Lambda_k+f,\beta_k}(v_{k-1})$
    for $\bar\Lambda_k(v):=
    \Lambda(y_k)+\langle\nabla\Lambda(y_k),v-y_k\rangle$
    \State $\psi_k(z):=
    \beta_kD_h(z,v_{k-1})+\bar\Lambda_k(z)$
    \State $\widehat g_k=
    -\frac{1}{\beta_k}\,
    \ell_p\textsc{-MorGradEst}\!\left(
    \widehat{\nabla}f,\psi_k,\beta_k, \beta_k\delta_k,\beta_k^q\sigma_k^q,\bar X\right)$
    \State $v_k=\displaystyle\arg\min_{v\in\bar X}
    \left\{h(v)-\langle\widehat g_k,v\rangle\right\}$
    \State $x_k=(1-\gamma_k)x_{k-1}+\gamma_k\bar v_k$
\EndFor
\State \Return $x_N$
\end{algorithmic}
\end{algorithm}

\begin{lemma}
\label{lem:composite_morgradest}
Let $f$ and $\widehat{\nabla}f$ satisfy Assumption \ref{Assumption}, $p > 2$, and let
\[v_k^\star=\arg\min_{z\in\bar X}\{f(z)+\psi_k(z)\}.\]
The procedure
\[\widehat g_k= -\frac1{\beta_k}\,
\ell_p\text{-}\mathrm{MorGradEst}\!\left(\widehat{\nabla}f,\psi_k,2^{p-2}\beta_k,
\beta_k\delta_k,\beta_k^q\sigma_k^q,\bar X\right)\]
outputs an estimator satisfying
\[\bigl\|\mathbb E\widehat g_k-\nabla h(v_k^\star)\bigr\|_q
\le\delta_k \quad\text{and}\quad
\mathbb E\bigl\|
\widehat g_k-\nabla h(v_k^\star)
\bigr\|_q^q\le\sigma_k^q,\]
using $N_{\widehat{\nabla}f}$ stochastic-gradient
computations, where
\[\mathbb E N_{\widehat{\nabla}f}
\le N\left( 1+
\frac32\lfloor\log_2T_{\max}\rfloor\right)\]
The parameters are
\[T_{\max}=\left\lceil
2\left(\frac{2^{p-1}(p-1)c^{1/p}
R^{p-2}(G/\beta_k)^{1/(p-1)}
}{\min\{\delta_k,\sigma_k\}}
\right)^{p(p-1)}\right\rceil\]
and
\[N=\left\lceil\frac{
2^{p^2+p+1}(p-1)^p c}{\sigma_k^p
}\left(\frac{1+2^{1/(p-1)^2}
}{1-2^{-(p-2)/(p-1)^2}}
\right)^{p-1}R^{p(p-2)}
\left(\frac{G}{\beta_k}\right)^q
T_{\max}^{(p-2)/(p-1)}
\right\rceil,\]
where $c$ is the constant in the
$\mathrm{ODC}_p$ guarantee.
\end{lemma}

\begin{proof}
The proof follows the MLMC argument in
Lemma~\ref{lem:morgrad_mlmc} and the averaging argument
in Theorem~\ref{thm:MorGradEst}.
Since $f+\psi_k$ is $\beta_k$-uniformly convex,
the ODC iterates satisfy
\[\mathbb E\|x_j-v_k^\star\|_p^p
\le c(G/\beta_k)^q2^{-j/(p-1)}.\]
Moreover, on $B_p(v_0,R)$,
\[\|\nabla h(u)-\nabla h(v)\|_q
\le 2^{p-2}(p-1)R^{p-2}\|u-v\|_p.\]
Consequently, the original MLMC estimator $g'$,
after normalization $Z=-g'/\beta_k$, satisfies
\[\|\mathbb E Z-\nabla h(v_k^\star)\|_q \le
2^{p-2+1/[p(p-1)]}(p-1)c^{1/p}
R^{p-2}(G/\beta_k)^{1/(p-1)}
T_{\max}^{-1/[p(p-1)]}.\]
Evaluating the geometric series in the same
moment calculation gives
\[\mathbb E\|Z-\mathbb E Z\|_q^q
\le\frac{2^{3q-2}\bigl(1+2^{1/(p-1)^2}\bigr)
}{1-2^{-(p-2)/(p-1)^2}}\left[
2^{p-2}(p-1)c^{1/p}
R^{p-2}(G/\beta_k)^{1/(p-1)}
\right]^qT_{\max}^{(p-2)/(p-1)^2}.\]
The choice of $T_{\max}$ makes the bias at most
$\min\{\delta_k,\sigma_k\}/2$, and the choice of $N$
makes the centered $q$th moment of the averaged
estimator at most $(\sigma_k/2)^q$.
Combining these bounds as in
Theorem~\ref{thm:MorGradEst} proves the guarantees. Finally, each MLMC sample uses at most
\[1+\sum_{j=1}
^{\lfloor\log_2T_{\max}\rfloor}
2^{-j}(2^j+2^{j-1})= 1+\frac32\lfloor\log_2T_{\max}\rfloor\]
oracle calls in expectation, yielding the stated cost.
\end{proof}

\begin{proposition}
\label{proposition:convex_composite}
Suppose $x^*$ minimizes $\Psi$ such that
$\|x_0-x^*\|_p\le R$, and let the iterates be generated
by Algorithm~\ref{alg:composite_convex_opt_alg_label}.
Suppose $M$ satisfies
Lemma~\ref{lemma:composite_minimization_helper}
with approximation error $\eta>0$.
Choose
\[
\sigma_k^q=\frac{qR^p}{2N},
\qquad
\delta_k=\frac{R^{p-1}}{4N},
\qquad
\epsilon_k=
\frac{p^{p-2}MR^p}{(k+p-1)^{p-1}N},
\]
where also
$\gamma_k=p/(k+p-1)$ and
$\beta_k=M\gamma_k^{p-1}$.
Further, for all $k\ge2$, let
$\Gamma_k=(1-\gamma_k)\Gamma_{k-1}$,
where $\Gamma_1=1$.
Then
\begin{align*}
\mathbb E[\Psi(x_N)-\Psi(x^*)]
&\le
\left(2^{p-1}+2+\frac1p\right)
\frac{Mp^pR^p}{(N+p-1)^p}
+\Gamma_N\sum_{k=1}^N\frac{\eta}{2\Gamma_k}\\
&=
\left(2^{p-1}+2+\frac1p\right)
\frac{Mp^pR^p}{(N+p-1)^p}
+\frac{\eta(N+p)}{2(p+1)}.
\end{align*}
\end{proposition}

\begin{proof}
Let
\[
\Phi_k(v)
=\bar\Lambda_k(v)+f(v)+\beta_kD_h(v,v_{k-1}),
\qquad
v_k^\star=\arg\min_{v\in\bar X}\Phi_k(v).
\]
Recall that
\[
D_h(u,v)\ge\frac1p\|u-v\|_p^p.
\]
By Lemma~\ref{lemma:composite_minimization_helper},
the definitions of $x_k,y_k$, and
$\beta_k=M\gamma_k^{p-1}$, we obtain
\begin{align*}
\Lambda(x_k)
&\le
(1-\gamma_k)\Lambda(x_{k-1})
+\gamma_k\bar\Lambda_k(\bar v_k)
+\frac{M\gamma_k^p}{p}
\|\bar v_k-v_{k-1}\|_p^p+\frac{\eta}{2}\\
&\le
(1-\gamma_k)\Lambda(x_{k-1})
+\gamma_k\left[
\bar\Lambda_k(\bar v_k)
+\beta_kD_h(\bar v_k,v_{k-1})
\right]+\frac{\eta}{2}.
\end{align*}
Adding the convexity inequality for $f$ gives
\[\Psi(x_k)\le
(1-\gamma_k)\Psi(x_{k-1})
+\gamma_k\Phi_k(\bar v_k)+\frac{\eta}{2}.\]
For any fixed $u\in\bar X$, convexity of $\Lambda$
and optimality of $v_k^\star$ give
\[\Psi(u)\ge\Phi_k(u)-\beta_kD_h(u,v_{k-1}),\qquad
\Phi_k(u)\ge \Phi_k(v_k^\star)+\beta_kD_h(u,v_k^\star).\]
Consequently,
\begin{align*}
\Psi(x_k)-\Psi(u) &\le
(1-\gamma_k)[\Psi(x_{k-1})-\Psi(u)]+\gamma_k[\Phi_k(\bar v_k)-\Phi_k(v_k^\star)]\\
&\quad+\gamma_k\beta_k
[D_h(u,v_{k-1})-D_h(u,v_k^\star)]
+\frac{\eta}{2}.
\end{align*}
Now, define
\[D_0=D_h(u,v_0), \qquad
D_k= \mathbb E\left[
h(u)-h(v_k)-\langle\widehat g_k,u-v_k\rangle\right],
\quad k \ge 1.\]
Optimality of $v_k$ in the mirror update implies
\[h(u)-h(v_k)-\langle\widehat g_k,u-v_k\rangle
\ge D_h(u,v_k)\ge0.\]
Moreover, Young's inequality and
$D_h(v_k,v_k^\star)\ge
\|v_k-v_k^\star\|_p^p/p$ imply
\begin{align*}
&h(u)-h(v_k)-\langle\widehat g_k,u-v_k\rangle\\
&\quad\le
D_h(u,v_k^\star)
+\langle
\widehat g_k-\nabla h(v_k^\star),v_k^\star-u
\rangle
+\frac1q
\|\widehat g_k-\nabla h(v_k^\star)\|_q^q.
\end{align*}
Conditioning on $\mathcal F_{k-1}$ and applying
Lemma~\ref{lem:composite_morgradest}, together with
$\|v_k^\star-u\|_p\le2R$, therefore gives
\begin{align*}
D_k \le \mathbb E D_h(u,v_k^\star)+2R\delta_k+\frac{\sigma_k^q}{q}= \mathbb E
D_h(u,v_k^\star)+\frac{R^p}{N}
\end{align*}
Also,
$\mathbb E D_h(u,v_{k-1})\le D_{k-1}$.
By the definition of the approximate proximal step,
\[\mathbb E[\Phi_k(\bar v_k)-\Phi_k(v_k^\star)] \le
\epsilon_k.\]
Taking expectations in the preceding recursion yields
\begin{align*}
\mathbb E[\Psi(x_k)-\Psi(u)]
&\le
(1-\gamma_k)\mathbb E[\Psi(x_{k-1})-\Psi(u)]
+\gamma_k\epsilon_k\\
&\quad+
\gamma_k\beta_k(D_{k-1}-D_k)
+\frac{\gamma_k\beta_kR^p}{N}
+\frac{\eta}{2}.
\end{align*}
For $u,v\in\bar X$, H\"older's inequality gives
\[D_h(u,v)\le \frac{2^{p-2}}pR^p
+\frac{2^{p-2}(p-1)}pR^p
+2^{p-2}R^p=2^{p-1}R^p.\]
Hence, including $k=0$,
\[0\le D_k\le\left(2^{p-1}+1\right)R^p.\]
Since $\gamma_1=1$, dividing the recursion by
$\Gamma_k$ and summing gives, with
$a_k=\gamma_k\beta_k/\Gamma_k$,
\begin{align*}
\mathbb E[\Psi(x_N)-\Psi(u)]
&\le
\Gamma_N\sum_{k=1}^N a_k(D_{k-1}-D_k)\\
&\quad+
\frac{(1+1/p)R^p\Gamma_N}{N}
\sum_{k=1}^N a_k
+\Gamma_N\sum_{k=1}^N\frac{\eta}{2\Gamma_k},
\end{align*}
where we used
$\epsilon_k=\beta_kR^p/(pN)$.
By Bernoulli's inequality,
\[\frac{a_{k+1}}{a_k}=\frac{k+p}{k}\left(1-\frac1{k+p}\right)^p
\ge1.\]
Therefore,
\begin{align*}
\sum_{k=1}^N a_k(D_{k-1}-D_k)
&=a_1D_0+\sum_{k=1}^{N-1}(a_{k+1}-a_k)D_k-a_ND_N\\
&\le (2^{p-1}+1)R^pa_N,
\end{align*}
and $\sum_{k=1}^Na_k\le Na_N$.
Using
\[
\Gamma_Na_N
=\gamma_N\beta_N
=M\gamma_N^p
=\frac{Mp^p}{(N+p-1)^p},
\]
we conclude that
\[
\mathbb E[\Psi(x_N)-\Psi(u)]
\le
\left(2^{p-1}+2+\frac1p\right)
\frac{Mp^pR^p}{(N+p-1)^p}
+\Gamma_N\sum_{k=1}^N\frac{\eta}{2\Gamma_k}.
\]
Finally, the recurrence for $\Gamma_k$ gives
\[
\Gamma_N\sum_{k=1}^N\frac1{\Gamma_k}
=\frac{N+p}{p+1}.
\]
Setting $u=x^*$ finishes the proof.
\end{proof}

\begin{theorem}
Suppose the conditions of
Proposition~\ref{proposition:convex_composite} hold
and let $\Lambda$ be $(L,\kappa)$-weakly smooth
with $\kappa\in(1,2]$.
For a target accuracy $0<\eps\le LR^\kappa$, choose
\begin{flalign*}
    &\eta=\frac{\eps}{N},
    \qquad
    M=\left[
    \frac{2(p-\kappa)}{p\kappa\eta}
    \right]^{\frac{p-\kappa}{\kappa}}
    L^{\frac p\kappa}.
\end{flalign*}
With $N$ chosen sufficiently large such that
\begin{flalign*}
    &N=\Theta\left(
    \left(\frac{LR^\kappa}{\eps}\right)^{
    \frac{p}{\kappa p-p+\kappa}}
    \right),
\end{flalign*}
Algorithm~\ref{alg:composite_convex_opt_alg_label},
using EpochDA for the inner computations,
returns an iterate satisfying
$\mathbb E[\Psi(x_N)]\le\Psi(x^*)+\eps$.
The number of evaluations of $\nabla\Lambda$ is
\begin{flalign*}
N_{\nabla\Lambda}=N=O\left(
\left(\frac{LR^\kappa}{\eps}\right)^{\frac{p}{\kappa p-p+\kappa}}\right),
\end{flalign*}
and the expected number of stochastic-gradient
evaluations, including both inner computations, is
\begin{flalign*}
    &\mathbb E N_{\nabla f}
    =\widetilde O\left(
    N+\left(\frac{GR}{\eps}\right)^p
    N^{p(p-2)}
    \right).
\end{flalign*}
\end{theorem}

\begin{proof}
By Proposition~\ref{proposition:convex_composite},
\begin{flalign*}
    &\mathbb E[\Psi(x_N)-\Psi(x^*)]
    \le
    \left(2^{p-1}+2+\frac1p\right)
    \frac{Mp^pR^p}{(N+p-1)^p}
    +\frac{\eta(N+p)}{2(p+1)}.
\end{flalign*}
Choose $\eta=\eps/N$.
Since $N+p\le(p+1)N$, the second term is at most
$\eps/2$.
Substituting this choice into the definition of $M$
gives
\begin{flalign*}
    &M=
    \left[
    \frac{2(p-\kappa)N}{p\kappa\eps}
    \right]^{\frac{p-\kappa}{\kappa}}
    L^{\frac p\kappa}.
\end{flalign*}
Now choose $N$ sufficiently large such that
\begin{flalign*}
    &N\ge
    \left[
    2p^p\left(2^{p-1}+2+\frac1p\right)
    \left(\frac{2(p-\kappa)}{p\kappa}\right)^{
    \frac{p-\kappa}{\kappa}}
    \right]^{\frac{\kappa}{\kappa p-p+\kappa}}
    \left(\frac{LR^\kappa}{\eps}\right)^{
    \frac{p}{\kappa p-p+\kappa}}.
\end{flalign*}
Using $N+p-1\ge N$, this ensures
\begin{flalign*}
    &\left(2^{p-1}+2+\frac1p\right)
    \frac{Mp^pR^p}{(N+p-1)^p}
    \le\frac{\eps}{2}.
\end{flalign*}
Therefore,
$\mathbb E[\Psi(x_N)-\Psi(x^*)]\le\eps$.
Each outer iteration evaluates $\nabla\Lambda$
once, so $N_{\nabla\Lambda}=N$.
For the stochastic-gradient complexity, at each
outer iteration $k$, substitute
\begin{flalign*}
    &\delta_k=\frac{R^{p-1}}{4N},
    \qquad
    \sigma_k=
    R^{p-1}\left(\frac{q}{2N}\right)^{1/q}
\end{flalign*}
into Lemma~\ref{lem:composite_morgradest}.
For fixed $p>2$, $\delta_k$ and $\sigma_k$
scale as $N^{-1}$ and $N^{-1/q}$, respectively,
and $\min\{\delta_k,\sigma_k\}=\delta_k$.
Consequently, the MLMC cost at iteration $k$ is
\begin{flalign*}
    &\widetilde O\left(
    1+
    \left(\frac{G}{\beta_kR^{p-1}}\right)^p
    N^{p^2-p-1}
    \right).
\end{flalign*}
To compute $\bar v_k$, we choose
$z_0\in\arg\min_{z\in\bar X}\psi_k(z)$
and set $R_0=(pG/\beta_k)^{1/(p-1)}$. An expected suboptimality of at most $\epsilon_k$
can then be achieved using
$O(1+G^p/(\beta_k\epsilon_k^{p-1}))$ stochastic subgradient calls, where
$\epsilon_k=\beta_kR^p/(pN)$.
The total cost at iteration $k$ therefore satisfies
\begin{flalign*}
    &\mathbb E N_k
    =\widetilde O\left(
    1+
    \left(\frac{G}{\beta_kR^{p-1}}\right)^p
    N^{p^2-p-1}
    \right).
\end{flalign*}
Since $\beta_k=M\gamma_k^{p-1}$ and
$\gamma_k\ge1/N$, summing over the $N$ outer
iterations gives
\begin{flalign*}
    &\mathbb E N_{\nabla f}
    =\widetilde O\left(
    N+
    \left(\frac{G}{MR^{p-1}}\right)^p
    N^{2p(p-1)}
    \right).
\end{flalign*}
Finally, the choices of $M$ and $N$ imply
$M=\Theta(\eps N^p/R^p)$.
Substituting this relation yields
\begin{flalign*}
    &\mathbb E N_{\nabla f}
    =\widetilde O\left(
    N+\left(\frac{GR}{\eps}\right)^p
    N^{p(p-2)}
    \right),
\end{flalign*}
which finishes the proof.
\end{proof}

\subsection{$\ell_p$ Accelerated Proximal Point Methods}
\label{subsection:lp_proximal_point}
In this section, we establish the guarantees of our stochastic accelerated $\ell_p$ proximal point method. We first analyze the NextLambda procedure, showing that it controls the proximal point movement while requiring only logarithmically many BROO calls. We then establish a potential decrease bound using the $\ell_p$ Bregman divergence and our reduced-bias Moreau gradient estimator. Combining these results, we obtain convergence and stochastic gradient query complexity guarantees for the accelerated proximal point method.
We define \begin{flalign*}
        &\tilde{P}_{f, \lambda}^{\varphi}(y) \in \{x \in X: \mathbb{E}[F(x)] \le F(P_{f, \lambda}(y)) + \varphi, F(x) = f(x) + \frac{\lambda}{p}||x - y||_p^p \}
    \end{flalign*}
\begin{algorithm}[H]
\caption{Stochastic accelerated proximal point method}\label{alg:stochastic_accelerated_prox_point}
\begin{algorithmic}[1]
\State Initialize $x_0 = v_0 = 0$, $A_0 \ge 0$, gradient estimator $\hat{\nabla} f$, function NextLambda
\For{$k = 0, 1, \dots,$}
    \State $\lambda_{k + 1} = $ NextLambda$(x_k, v_k, A_k)$
    \State Solve $(A_k + a_{k + 1})^{p - 1} = 2^{p - 1} \lambda_{k + 1} a_{k + 1}^p$ for $a_{k + 1}$
    \State $A_{k + 1} = A_k + a_{k + 1}$
    \State $y_k = \frac{A_k}{A_{k + 1}} x_k + \frac{a_{k + 1}}{A_{k + 1}} v_k$
    \State $x_{k + 1} = \tilde{P}_{f, \lambda_{k + 1}}^{\varphi_{k + 1}}(y_k)$
    \State $g_{k + 1} = $ MORGRADEST$(\hat{\nabla}f, y_k, \lambda_{k + 1}, \delta_{k + 1}, \sigma_{k + 1}^q, X \cap \mathbb{B}_r(y_k))$
    \State $v_{k + 1} = \argmin_{v \in X}\{\frac{1}{2} a_{k + 1} \langle g_{k + 1}, v \rangle + D_\Phi(v, v_k) \}$, $\Phi(x) = \frac{1}{p}||x||_p^p$
    \State \Return
$x_{k + 1}$ if $A_{k + 1} \ge A_{\max}$ or $k + 1 = K_{\max}$
\EndFor
\end{algorithmic}
\end{algorithm}

\subsubsection{NextLambda Algorithm}

\begin{algorithm}[H]
\caption{NextLambda$(x,v,A)$}
\label{alg:next_lambda_alg}
\begin{algorithmic}[1]

\State Input:
$x,v\in\mathbb{R}^d$, $A\ge0$,
$\lambda_{\min},\lambda_{\max}$,
$G,R,r$
\State For $\lambda'>0$, let $a_{\lambda'}>0$ satisfy
$(A+a_{\lambda'})^{p-1}=
2^{p-1}\lambda'a_{\lambda'}^p$, where $A_{\lambda'} = A + a_{\lambda'}$ and $y_{\lambda'} = \frac{A}{A_{\lambda'}}x + \frac{a_{\lambda'}}{A_{\lambda'}}v$

\State Define
$\Delta(\lambda)=\left\|\mathcal{O}_{\lambda,r/17}(y_\lambda)-y_\lambda
\right\|_p.$
\State $\lambda\gets\lambda_{\max}$
\While{
$\lambda\ge\lambda_{\min}$
and
$\Delta(\lambda)\le 13r/16$
} $\lambda\gets\lambda/2$
\EndWhile
\If{$\lambda<\lambda_{\min}$} \Return $2\lambda$
\EndIf
\State
$\lambda_u=2\lambda$,
$\lambda_\ell=\lambda$,
$\lambda_m=\sqrt{\lambda_u\lambda_\ell}$
\If{$\Delta(\lambda_\ell)\le15r/16$}
\Return $\lambda_\ell$
\EndIf
\While{
$\Delta(\lambda_m)\notin[13r/16,15r/16]$ and $\log\frac{\lambda_u}{\lambda_\ell}
>\frac{r}{8\left(R+\frac{1}{p-1}
\left(\frac{G}{\lambda_\ell}
\right)^{1/(p-1)}\right)}$
}
    \If{$\Delta(\lambda_m)<13r/16$}
    $\lambda_u=\lambda_m$
    \Else\ $\lambda_\ell=\lambda_m$
    \EndIf
    \State $\lambda_m=\sqrt{\lambda_u\lambda_\ell}$
\EndWhile
\State \Return $\lambda_m$

\end{algorithmic}
\end{algorithm}
Here, we discuss our implementation of the subroutine NextLambda
(inspired by \cite{asi2021stochastic} and \cite{carmon2021thinking}),
such that for $\hat{x}_{k+1}=P_{f,\lambda}(y_k)$, it holds that
$\|\hat{x}_{k+1}-y_k\|_p\le r$ and either
$\|\hat{x}_{k+1}-y_k\|_p\ge 3r/4$ or
$\lambda<2\lambda_{\min}$. We let $J_p(z) = \nabla (\frac{1}{p} ||z||_p^p) = (|z_i|^{p - 2} z_i)_i \quad \text{and} \quad \hat{x}_\lambda = \argmin_{x} \{f(x) + \frac{\lambda}{p}||x - y_\lambda||_p^p \}$. Assume that $f$ is $G$-Lipschitz with respect to the $\ell_p$ norm. 

\begin{definition}
\label{definition:broo_oracle}
We say that a mapping
$\mathcal{O}_{\lambda,\rho}(\cdot):\mathcal{X}\to\mathbb{R}^d$
is a Ball Regularized Optimization Oracle of radius $r$ for $f$ if,
for every query point $\overline{x}$, regularization parameter
$\lambda>0$, and desired accuracy $\rho>0$, it returns
$\tilde{x}=\mathcal{O}_{\lambda,\rho}(\overline{x})\in \mathcal{X} \cap \mathbb{B}_r(\overline{x})$
such that
\begin{flalign*}
f(\tilde{x})+\frac{\lambda}{p}
\|\tilde{x}-\overline{x}\|_p^p
\le\min_{x\in\mathbb{B}_r(\overline{x})}\left\{f(x)+\frac{\lambda}{p}
\|x-\overline{x}\|_p^p\right\}
+\frac{\lambda 2^{2-p}}{p}\rho^p.
\end{flalign*}
\end{definition}

\begin{lemma}
\label{lemma:x_lambda_y_lambda_diff_bound}
Let $\hat{x}_\lambda = \argmin_{z} \{f(z) + \frac{\lambda}{p} ||z - y_\lambda||_p^p \}$. For all $\lambda > 0$, it holds that \begin{flalign*}
            &||\hat{x}_\lambda - y_\lambda||_p \le (\frac{G}{\lambda})^{1/(p - 1)}
        \end{flalign*}
\end{lemma}

\begin{proof}
    Let $\delta_\lambda = \hat{x}_\lambda - y_\lambda$. We know that for some $g_\lambda \in \partial f(\hat{x}_\lambda)$, we have that $g_\lambda + \lambda J_p(\delta_\lambda) = 0$. Since $f$ is $G$-Lipschitz, it holds that $||g_\lambda||_q \le G$, which implies that $\lambda ||J_p(\delta_\lambda)||_q \le G$. Since $||J_p(z)||_q = ||z||_p^{p - 1}$, it holds that $\lambda \|\delta_\lambda\|_p^{p - 1} \leq G$ from which the stated claim follows. 
\end{proof}

\begin{lemma}
\label{lemma:y_lambda_derivative}
Let $A\ge 0$ and $x,v\in\mathbb R^d$
satisfy $\|x-v\|_p\le 2R$.
For every $\lambda>0$, let $a_\lambda>0$
be the unique solution of
$(A+a_\lambda)^{p-1}
=2^{p-1}\lambda a_\lambda^p$,
and define
$y_\lambda=\frac{Ax+a_\lambda v}{A+a_\lambda}$.
Then
\[\left\|\frac{d}{d\lambda}y_\lambda\right\|_p
\le\frac{R}{\lambda}.\]
\end{lemma}

\begin{proof}
 Denote $\alpha_\lambda = a_\lambda/(A + a_\lambda)$. Taking logarithmic derivatives on both sides of $(A + a_\lambda)^{p - 1} = 2^{p - 1}\lambda a_\lambda^p$, we get that \begin{flalign*}
&(p - 1) \frac{a_\lambda'}{A + a_\lambda} = \frac{1}{\lambda} + p \frac{a_\lambda'}{a_\lambda} \implies a_\lambda' = -\frac{a_\lambda(A + a_\lambda)}{\lambda(a_\lambda + (pA))}
\end{flalign*} Since $a_\lambda + pA \ge A + a_\lambda$, this means that $\lvert a_\lambda^\prime\rvert\le a_\lambda/\lambda$. Now, notice that \begin{flalign*}
&|\alpha_\lambda'| = |\frac{A a_\lambda'}{(A + a_\lambda)^2} | \le \frac{Aa_\lambda}{\lambda(A + a_\lambda)^2} \le \frac{1}{2 \lambda}\end{flalign*} and since $y_\lambda = (1 - \alpha_\lambda)x + \alpha_\lambda v$, this implies that \begin{flalign*}
&||\frac{d}{d \lambda} y_\lambda||_p = |\alpha_\lambda'| \cdot ||v - x||_p \le \frac{R}{\lambda}\end{flalign*} as $\|v-x\|_p \leq 2R$.
\end{proof}

\begin{lemma}
\label{lemma:delta_norm_continuity}
Let $f:\mathbb{R}^d\to\mathbb{R}$ be convex and $G$-Lipschitz,
let $y_\lambda$ be continuously differentiable for $\lambda>0$.
For every $0<\lambda_1\le\lambda_2$, it holds that
\[\left|\|\delta_{\lambda_2}\|_p
-\|\delta_{\lambda_1}\|_p\right|
\leq \int_{\lambda_1}^{\lambda_2}
\left\|\frac{d}{d\lambda}y_\lambda
\right\|_p\,d\lambda+\frac{1}{p-1}
\int_{\lambda_1}^{\lambda_2}
\frac{\|\delta_\lambda\|_p}{\lambda}
\,d\lambda.\]
\end{lemma}

\begin{proof}
We first smooth $f$ and use a nondegenerate approximation of the
$p$-power regularizer. Let
\[f_\sigma(x)=\mathbb{E}[f(x+\nu_\sigma)],\qquad\nu_\sigma\sim\mathcal{N}(0,\sigma^2I),\]
and, for $\epsilon>0$, define
$\phi_\epsilon(z)=\frac{1}{p}
\sum_{i=1}^d\left[(z_i^2+\epsilon^2)^{p/2}-\epsilon^p\right].$
Let \[
\hat{x}_\lambda^{\sigma,\epsilon}
=\argmin_x\left\{f_\sigma(x)+\lambda
\phi_\epsilon(x-y_\lambda)\right\},
\qquad\delta_\lambda^{\sigma,\epsilon}=\hat{x}_\lambda^{\sigma,\epsilon}
-y_\lambda,\]
and define $\rho_\epsilon(z)=
\left(\sum_{i=1}^d
(z_i^2+\epsilon^2)^{p/2}
\right)^{1/p}$.
Notice that
$\phi_\epsilon(z)=\frac{1}{p}
\left(\rho_\epsilon(z)^p-d\epsilon^p
\right)$
and $\rho_\epsilon(z) \rightarrow\|z\|_p$ as $\epsilon\rightarrow 0.$
Since $f_\sigma$ is smooth and convex and
$\phi_\epsilon$ has a positive definite Hessian,
$\hat{x}_\lambda^{\sigma,\epsilon}$ is uniquely defined and
differentiable in $\lambda$.
Let $J_{p,\epsilon}(z)=\nabla\phi_\epsilon(z)$,
so that
\[[J_{p,\epsilon}(z)]_i
=z_i(z_i^2+\epsilon^2)^{p/2-1}.\]
The Hessian of $\phi_\epsilon$ is diagonal, with
$[\nabla^2\phi_\epsilon(z)]_{ii}
=(z_i^2+\epsilon^2)^{p/2-2}\left[
\epsilon^2+(p-1)z_i^2\right]>0.$
Now, the first-order optimality gives
\[\nabla f_\sigma
(\hat{x}_\lambda^{\sigma,\epsilon})
+\lambda J_{p,\epsilon}
(\delta_\lambda^{\sigma,\epsilon})
=0.\]
Differentiating with respect to $\lambda$ gives
\[\left[\nabla^2 f_\sigma
(\hat{x}_\lambda^{\sigma,\epsilon})
+\lambda\nabla^2\phi_\epsilon
(\delta_\lambda^{\sigma,\epsilon})
\right]\frac{d}{d\lambda}
\delta_\lambda^{\sigma,\epsilon}
=-J_{p,\epsilon}
(\delta_\lambda^{\sigma,\epsilon})
-\nabla^2 f_\sigma
(\hat{x}_\lambda^{\sigma,\epsilon})
\frac{d}{d\lambda}y_\lambda.\]
For brevity, write
\[H=\nabla^2 f_\sigma
(\hat{x}_\lambda^{\sigma,\epsilon}),
\qquad D=\nabla^2\phi_\epsilon
(\delta_\lambda^{\sigma,\epsilon}).\]
 By construction we have $H\succeq0$ and $D\succ0$, and
\[(H+\lambda D)\frac{d}{d\lambda}
\delta_\lambda^{\sigma,\epsilon}
=-J_{p,\epsilon}(\delta_\lambda^{\sigma,\epsilon})
-H\frac{d}{d\lambda}y_\lambda.\]
Since $\nabla\rho_\epsilon(z)
=\frac{J_{p,\epsilon}(z)}{
\rho_\epsilon(z)^{p-1}}$, we have
$\frac{d}{d\lambda}\rho_\epsilon
(\delta_\lambda^{\sigma,\epsilon})
=\frac{J_{p,\epsilon}(\delta_\lambda^{\sigma,\epsilon})^\top \frac{d}{d\lambda}
\delta_\lambda^{\sigma,\epsilon}
}{\rho_\epsilon
(\delta_\lambda^{\sigma,\epsilon})^{p-1}}.$
Therefore,
\[\left|\frac{d}{d\lambda}
\rho_\epsilon
(\delta_\lambda^{\sigma,\epsilon})
\right|\le\frac{J_{p,\epsilon}
(\delta_\lambda^{\sigma,\epsilon})^\top(H+\lambda D)^{-1}
J_{p,\epsilon}(\delta_\lambda^{\sigma,\epsilon})
}{\rho_\epsilon
(\delta_\lambda^{\sigma,\epsilon})^{p-1}}+\frac{\left|J_{p,\epsilon}
(\delta_\lambda^{\sigma,\epsilon})^\top(H+\lambda D)^{-1}H\frac{d}{d\lambda}y_\lambda\right|}{\rho_\epsilon(\delta_\lambda^{\sigma,\epsilon})^{p-1}}.\]
We first bound the first term.
Since $H+\lambda D\succeq\lambda D$
and $D\succ0$, we have
$(H+\lambda D)^{-1}\preceq\frac{1}{\lambda}D^{-1}$.
Hence
\[J_{p,\epsilon}
(\delta_\lambda^{\sigma,\epsilon})^\top(H+\lambda D)^{-1}J_{p,\epsilon}
(\delta_\lambda^{\sigma,\epsilon})
\le \frac{1}{\lambda}J_{p,\epsilon}
(\delta_\lambda^{\sigma,\epsilon})^\top D^{-1}J_{p,\epsilon}
(\delta_\lambda^{\sigma,\epsilon}).
\]
For every coordinate,
\[\frac{[J_{p,\epsilon}
(\delta_\lambda^{\sigma,\epsilon})]_i^2}{D_{ii}}=
\frac{|\delta_{\lambda,i}
^{\sigma,\epsilon}|^2\left(
|\delta_{\lambda,i}^{\sigma,\epsilon}|^2+\epsilon^2\right)^{p/2}}{
\epsilon^2+(p-1)
|\delta_{\lambda,i}^{\sigma,\epsilon}|^2}\le\frac{1}{p-1}\left(
|\delta_{\lambda,i}^{\sigma,\epsilon}|^2+\epsilon^2\right)^{p/2}.\]
Therefore,
\[J_{p,\epsilon}
(\delta_\lambda^{\sigma,\epsilon})^\top D^{-1}J_{p,\epsilon}
(\delta_\lambda^{\sigma,\epsilon})
\le\frac{\rho_\epsilon
(\delta_\lambda^{\sigma,\epsilon})^p
}{p-1},\]
and thus
\[J_{p,\epsilon}
(\delta_\lambda^{\sigma,\epsilon})^\top(H+\lambda D)^{-1}J_{p,\epsilon}
(\delta_\lambda^{\sigma,\epsilon})
\le\frac{\rho_\epsilon
(\delta_\lambda^{\sigma,\epsilon})^p
}{(p-1)\lambda}.\]
We next bound the second term.
Let $S=D^{-1/2}HD^{-1/2}$.
Since $H\succeq0$, we have $S\succeq0$, and
\[D^{-1/2}H(H+\lambda D)^{-1}
J_{p,\epsilon}(\delta_\lambda^{\sigma,\epsilon})=
S(S+\lambda I)^{-1}D^{-1/2}
J_{p,\epsilon}(\delta_
\lambda^{\sigma,\epsilon}).\]
Every eigenvalue of
$S(S+\lambda I)^{-1}$ is in $[0,1]$. Therefore,
\[
\left\|
D^{-1/2}
H
(H+\lambda D)^{-1}
J_{p,\epsilon}
(\delta_\lambda^{\sigma,\epsilon})
\right\|_2
\le
\left\|
D^{-1/2}
J_{p,\epsilon}
(\delta_\lambda^{\sigma,\epsilon})
\right\|_2.
\]
It follows from Cauchy--Schwarz that
\[\left|J_{p,\epsilon}
(\delta_\lambda^{\sigma,\epsilon})^\top(H+\lambda D)^{-1}H\frac{d}{d\lambda}y_\lambda\right|\le
\left[J_{p,\epsilon}
(\delta_\lambda^{\sigma,\epsilon})^\top D^{-1}J_{p,\epsilon}
(\delta_\lambda^{\sigma,\epsilon})
\right]^{1/2}\left[\left(
\frac{d}{d\lambda}y_\lambda\right)^\top D
\left(\frac{d}{d\lambda}y_\lambda
\right)\right]^{1/2}.\]
We already showed that
$J_{p,\epsilon}
(\delta_\lambda^{\sigma,\epsilon})^\top D^{-1}J_{p,\epsilon}
(\delta_\lambda^{\sigma,\epsilon})
\le\frac{\rho_\epsilon
(\delta_\lambda^{\sigma,\epsilon})^p
}{p-1}.$
Moreover,
\[D_{ii}=\left(
|\delta_{\lambda,i}^{\sigma,\epsilon}|^2+\epsilon^2\right)^{p/2-2}
\left[\epsilon^2+(p-1)
|\delta_{\lambda,i}^{\sigma,\epsilon}|^2\right]\le(p-1)\left(
|\delta_{\lambda,i}^{\sigma,\epsilon}|^2+\epsilon^2\right)^{(p-2)/2}.\]
Hence, with H\"older's inequality, we have
\[\left(\frac{d}{d\lambda}y_\lambda
\right)^\top D\left(
\frac{d}{d\lambda}y_\lambda
\right)\le(p-1)\sum_i\left(
|\delta_{\lambda,i}^{\sigma,\epsilon}|^2+\epsilon^2\right)^{(p-2)/2}
\left|\frac{d}{d\lambda}y_{\lambda,i}
\right|^2\le (p-1)\rho_\epsilon
(\delta_\lambda^{\sigma,\epsilon})^{p-2}\left\|\frac{d}{d\lambda}y_\lambda\right\|_p^2\]
Combining the preceding two inequalities gives
\[\left| J_{p,\epsilon}
(\delta_\lambda^{\sigma,\epsilon})^\top(H+\lambda D)^{-1}H\frac{d}{d\lambda}y_\lambda\right|\le
\rho_\epsilon
(\delta_\lambda^{\sigma,\epsilon})^{p-1}\left\|\frac{d}{d\lambda}y_\lambda\right\|_p.\]
Consequently,
\[\left|\frac{d}{d\lambda}
\rho_\epsilon
(\delta_\lambda^{\sigma,\epsilon})
\right|\le\left\|\frac{d}{d\lambda}y_\lambda\right\|_p+
\frac{\rho_\epsilon
(\delta_\lambda^{\sigma,\epsilon})
}{(p-1)\lambda}.\]
Integrating from $\lambda_1$ to $\lambda_2$ gives
\[\left|\rho_\epsilon
(\delta_{\lambda_2}^{\sigma,\epsilon})-\rho_\epsilon
(\delta_{\lambda_1}^{\sigma,\epsilon})\right|\le\int_{\lambda_1}^{\lambda_2}\left\|\frac{d}{d\lambda}y_\lambda\right\|_p\,d\lambda+\frac{1}{p-1}
\int_{\lambda_1}^{\lambda_2}
\frac{\rho_\epsilon
(\delta_\lambda^{\sigma,\epsilon})
}{\lambda}\,d\lambda.\]
It remains to let $\sigma,\epsilon\rightarrow0$.
Since $f$ is $G$-Lipschitz,
$f_\sigma$ is also $G$-Lipschitz and
\[|f_\sigma(x)-f(x)|\le
G\mathbb{E}\|\nu_\sigma\|_p
\longrightarrow0\]uniformly in $x$.
Now, by the definition of $\hat{x}_\lambda^{\sigma,\epsilon}$,\[
f_\sigma(\hat{x}_\lambda^{\sigma,\epsilon})+\lambda\phi_\epsilon
(\delta_\lambda^{\sigma,\epsilon}) =f_\sigma(\hat{x}_\lambda^{\sigma,\epsilon})+\lambda\phi_\epsilon
(\hat{x}_\lambda^{\sigma,\epsilon}-y_\lambda)\le f_\sigma(y_\lambda)+
\lambda\phi_\epsilon(y_\lambda-y_\lambda) = f_\sigma(y_\lambda)+
\lambda\phi_\epsilon(0)\]
We now compute $\phi_\epsilon(0)$. Since $\rho_\epsilon(0)^p=\sum_{i=1}^d
(0^2+\epsilon^2)^{p/2}=\sum_{i=1}^d
\epsilon^p=d\epsilon^p$,
we have
$\phi_\epsilon(0)=\frac{1}{p}\left[
\rho_\epsilon(0)^p-d
\epsilon^p\right]=\frac{1}{p}\left[
d\epsilon^p-d\epsilon^p\right]=0$.
Therefore, \[f_\sigma
(\hat{x}_\lambda^{\sigma,\epsilon})
+\lambda\phi_\epsilon
(\delta_\lambda^{\sigma,\epsilon})
\le f_\sigma(y_\lambda)\]
Using $\phi_\epsilon(z)=\frac{1}{p}
\left[\rho_\epsilon(z)^p-d\epsilon^p
\right]$ again, we obtain
\[f_\sigma
(\hat{x}_\lambda^{\sigma,\epsilon})
+\frac{\lambda}{p}\left[\rho_\epsilon
(\delta_\lambda^{\sigma,\epsilon})^p
-d\epsilon^p\right]\le
f_\sigma(y_\lambda) \Rightarrow \frac{\lambda}{p}\left[\rho_\epsilon
(\delta_\lambda^{\sigma,\epsilon})^p
-d\epsilon^p\right]
\le f_\sigma(y_\lambda)-f_\sigma
(\hat{x}_\lambda^{\sigma,\epsilon})\]
Since $f_\sigma$ is $G$-Lipschitz,
$f_\sigma(y_\lambda)-f_\sigma
(\hat{x}_\lambda^{\sigma,\epsilon})
\le G\left\|y_\lambda-
\hat{x}_\lambda^{\sigma,\epsilon}
\right\|_p \leq G
\|\delta_\lambda^{\sigma,\epsilon}\|_p$. Moreover,
\[\rho_\epsilon
(\delta_\lambda^{\sigma,\epsilon})^p
=\sum_{i=1}^d\left(
|\delta_{\lambda,i}^{\sigma,\epsilon}|^2+\epsilon^2\right)^{p/2}\ge
\sum_{i=1}^d|\delta_{\lambda,i}^{\sigma,\epsilon}|^p=\|\delta_\lambda^{\sigma,\epsilon}\|_p^p.\]
Hence$\|\delta_\lambda^{\sigma,\epsilon}\|_p\le\rho_\epsilon(\delta_\lambda^{\sigma,\epsilon})$.
Combining the preceding inequalities gives
\[\frac{\lambda}{p}\left[
\rho_\epsilon(\delta_\lambda
^{\sigma,\epsilon})^p-d\epsilon^p
\right]\le G\rho_\epsilon(\delta_
\lambda^{\sigma,\epsilon}) \Rightarrow \frac{\lambda}{p}
\rho_\epsilon(\delta_\lambda
^{\sigma,\epsilon})^p\le G
\rho_\epsilon(\delta_\lambda
^{\sigma,\epsilon})+
\frac{\lambda d}{p}\epsilon^p\]
Since $\lambda\in[\lambda_1,\lambda_2]$ with $\lambda_1>0$,
this implies that $\rho_\epsilon
(\delta_\lambda^{\sigma,\epsilon})$
is uniformly bounded over
$\lambda\in[\lambda_1,\lambda_2]$ for sufficiently small
$\sigma,\epsilon$. Therefore the corresponding minimizers remain in a common bounded
set. On every bounded set,
\[f_\sigma(x)+
\lambda\phi_\epsilon(x-y_\lambda)
\longrightarrow f(x)+\frac{\lambda}{p}\|x-y_\lambda\|_p^p\]
uniformly as $\sigma,\epsilon\rightarrow0$.
Since the limiting objective is strictly convex, its minimizer is
unique. Also,
$\hat{x}_\lambda^{\sigma,\epsilon}
\longrightarrow\hat{x}_\lambda$ and $\delta_\lambda^{\sigma,\epsilon}\longrightarrow\delta_\lambda$
for every $\lambda\in[\lambda_1,\lambda_2]$.
Moreover,
$\rho_\epsilon
(\delta_\lambda^{\sigma,\epsilon})
\longrightarrow\|\delta_\lambda\|_p.$
The uniform bound above and $\lambda\ge\lambda_1>0$ allow us to
apply dominated convergence to the integral. Taking
$\sigma,\epsilon\rightarrow0$ yields
\[\left|\|\delta_{\lambda_2}\|_p
-\|\delta_{\lambda_1}\|_p\right|
\le\int_{\lambda_1}^{\lambda_2}
\left\|\frac{d}{d\lambda}y_\lambda
\right\|_p\,d\lambda+\frac{1}{p-1}
\int_{\lambda_1}^{\lambda_2}
\frac{\|\delta_\lambda\|_p}{\lambda}
\,d\lambda,\]
which proves the claim.
\end{proof}

\begin{lemma}
\label{lemma:delta_hat_diff_bound}
Under the same context as
Lemma \ref{lemma:delta_norm_continuity}, let $\hat{\Delta}(\lambda) =
\|\hat{x}_\lambda-y_\lambda\|_p$.
If $\|x-v\|_p\le2R$, then for every
$0<\lambda_1\le\lambda_2$,
\[\left|\hat{\Delta}(\lambda_1)
-\hat{\Delta}(\lambda_2)\right|
\le\left[R+\frac{1}{p-1}\left(
\frac{G}{\lambda_1}\right)^{1/(p-1)}
\right]\log\frac{\lambda_2}{\lambda_1}.\]
\end{lemma}

\begin{proof}
By Lemma \ref{lemma:delta_norm_continuity},
\[\left|\hat{\Delta}(\lambda_1)-
\hat{\Delta}(\lambda_2)\right|
\le\int_{\lambda_1}^{\lambda_2}\left\|\frac{d}{d\lambda}y_\lambda
\right\|_p\,d\lambda+\frac{1}{p-1}
\int_{\lambda_1}^{\lambda_2}
\frac{\hat{\Delta}(\lambda)}{\lambda}\,d\lambda.\]
Lemma \ref{lemma:y_lambda_derivative} gives
$\left\|\frac{d}{d\lambda}y_\lambda
\right\|_p \le \frac{R}{\lambda}$,
and Lemma \ref{lemma:x_lambda_y_lambda_diff_bound} gives
\[\hat{\Delta}(\lambda)\le \left(\frac{G}{\lambda}\right)^{1/(p-1)} \le
\left(\frac{G}{\lambda_1}
\right)^{1/(p-1)}\] for every $\lambda\in[\lambda_1,\lambda_2]$.
Therefore,
\[\left|\hat{\Delta}(\lambda_1)-
\hat{\Delta}(\lambda_2)\right|\le
\left[R+\frac{1}{p-1}\left(
\frac{G}{\lambda_1}\right)^{1/(p-1)}
\right]\int_{\lambda_1}^{\lambda_2}
\frac{d\lambda}{\lambda},\]
which proves the claim.
\end{proof}

\begin{lemma}
\label{lemma:nextlambda_helper_lemma}
Let $\hat{x}_\lambda=\argmin_x\left\{f(x)+\frac{\lambda}{p}\|x-y_\lambda\|_p^p\right\}$ and $
\hat{\Delta}(\lambda)
=\|\hat{x}_\lambda-y_\lambda\|_p$. Define
\[\Delta(\lambda)=\left\|
\mathcal{O}_{\lambda,\rho}(y_\lambda)-y_\lambda\right\|_p.\]
If either $\hat{\Delta}(\lambda)\le r-\rho$ or $\Delta(\lambda)<r-\rho$,
then
\[\bprox_{\lambda,r}^f(y_\lambda)
=\prox_\lambda^f(y_\lambda)\quad 
\text{and}\quad|\Delta(\lambda)-\hat{\Delta}(\lambda)|\le\rho.\]
where $\bprox_{\lambda,r}^f(y_\lambda) = x^{\star}_{\lambda,r}$ and $\prox_\lambda^f(y_\lambda) = \hat{x}_{\lambda}$
\end{lemma}
\begin{proof}
Let
\[F_\lambda(x)=f(x)+\frac{\lambda}{p}\|x-y_\lambda\|_p^p,\qquad
x_{\lambda,r}^*=\argmin_{x\in\mathbb B_r(y_\lambda)}F_\lambda(x).\]
By Definition~\ref{definition:broo_oracle},
\begin{equation}
\label{eq:nextlambda_broo_accuracy}
F_\lambda(
\mathcal O_{\lambda,\rho}(y_\lambda))
\le
F_\lambda(x_{\lambda,r}^*)
+\frac{\lambda 2^{2-p}}{p}\rho^p.
\end{equation}
\paragraph{Case 1:
$\hat{\Delta}(\lambda)\le r-\rho$.}
Since $\hat{x}_\lambda$ belongs to
$\mathbb B_r(y_\lambda)$ and minimizes
$F_\lambda$ globally, we have
$x_{\lambda,r}^*=\hat{x}_\lambda$.
By the $p$-uniform convexity of
$\frac1p\|\cdot\|_p^p$ and the global
optimality of $\hat{x}_\lambda$,
\[F_\lambda(\mathcal O_{\lambda,\rho}(y_\lambda))
-F_\lambda(\hat{x}_\lambda)\ge
\frac{\lambda 2^{2-p}}{p}
\|\mathcal O_{\lambda,\rho}(y_\lambda)-\hat{x}_\lambda\|_p^p.\]
Combining this with
\eqref{eq:nextlambda_broo_accuracy}
gives
\begin{equation}
\label{eq:nextlambda_distance_bound}
\|\mathcal O_{\lambda,\rho}(y_\lambda)
-\hat{x}_\lambda\|_p
\le\rho.
\end{equation}
\paragraph{Case 2:
$\Delta(\lambda)<r-\rho$.}
By the definition of $\Delta(\lambda)$,
$\mathcal O_{\lambda,\rho}(y_\lambda)
\in\mathbb B_r(y_\lambda)$.
Since $x_{\lambda,r}^*$ minimizes
$F_\lambda$ over this convex ball,
the constrained first-order optimality
condition and $p$-uniform convexity imply
\[
F_\lambda(
\mathcal O_{\lambda,\rho}(y_\lambda))
-F_\lambda(x_{\lambda,r}^*)
\ge
\frac{\lambda 2^{2-p}}{p}
\|\mathcal O_{\lambda,\rho}(y_\lambda)
-x_{\lambda,r}^*\|_p^p.
\]
Together with
\eqref{eq:nextlambda_broo_accuracy},
this yields
\[
\|\mathcal O_{\lambda,\rho}(y_\lambda)
-x_{\lambda,r}^*\|_p
\le\rho.
\]
Consequently,
\begin{align*}
\|x_{\lambda,r}^*-y_\lambda\|_p
&\le
\|x_{\lambda,r}^*
-\mathcal O_{\lambda,\rho}(y_\lambda)
\|_p+\Delta(\lambda)\le\rho+\Delta(\lambda)<r.
\end{align*}
Thus $x_{\lambda,r}^*$ lies in the
interior of $\mathbb B_r(y_\lambda)$.
By convexity, an interior minimizer
over the ball is also a global
minimizer of $F_\lambda$. Hence
$x_{\lambda,r}^*=\hat{x}_\lambda$,
and \eqref{eq:nextlambda_distance_bound}
holds in this case as well.
Therefore, in both cases,
\[\bprox_{\lambda,r}^f(y_\lambda)
=\prox_\lambda^f(y_\lambda).\]
Finally, the reverse triangle
inequality gives
\begin{align*}
|\Delta(\lambda)-\hat{\Delta}(\lambda)|&\le
\|\mathcal O_{\lambda,\rho}(y_\lambda)
-\hat{x}_\lambda\|_p\le\rho.
\end{align*}
\end{proof}
\begin{theorem}
\label{theorem:nextlambda_guarantee}
Let $f:\mathbb{R}^d\to\mathbb{R}$ be convex and $G$-Lipschitz,
let $x,v\in\mathbb{R}^d$ satisfy
$\|x-v\|_p\le2R$, $0<r\le R$,
and let
\[\hat{x}_\lambda=\argmin_x\left\{
f(x)+\frac{\lambda}{p}\|x-y_\lambda\|_p^p\right\},\qquad
\hat{\Delta}(\lambda)=
\|\hat{x}_\lambda-y_\lambda\|_p.\]
Given that $\lambda_{\max}\ge\frac{2^{p-1}G}{r^{p-1}}$ and $\lambda_{\min}
\in(0,\lambda_{\max})$,
Algorithm \ref{alg:next_lambda_alg} returns
$\lambda\in[\lambda_{\min},\lambda_{\max}]$ such that $\hat{\Delta}(\lambda)\le r$. Moreover, either
\begin{itemize}
\item $\lambda <2\lambda_{\min}$,
\item $\lambda \in [2\lambda_{\min}, \lambda_{\max}] \text{ and }\hat{\Delta}(\lambda) \in
\left[\frac{3r}{4},r\right)$.
\end{itemize}
Furthermore, the overall BROO complexity is
\[O\left(\log\frac{\lambda_{\max}}{\lambda_{\min}}+\log\left(1+\frac{
R+\frac{1}{p-1}\left(\frac{G}{\lambda_{\min}}\right)^{1/(p-1)}}{
r}\right)\right).\]
\end{theorem}

\begin{proof}
First, by
Lemma \ref{lemma:x_lambda_y_lambda_diff_bound}, \begin{flalign*}
&\hat{\Delta}(\lambda_{\max}) \le (\frac{G}{\lambda_{\max}})^{1/(p - 1)} \le  (\frac{G}{2^{p - 1}G/r^{p - 1}})^{1/(p - 1)} = \frac{r}{2}. 
\end{flalign*}
Since $\frac{r}{2}<r-\frac{r}{17}$,
Lemma \ref{lemma:nextlambda_helper_lemma} with
$\rho=r/17$ gives
\[\Delta(\lambda_{\max})\le
\hat{\Delta}(\lambda_{\max})+
\frac{r}{17}\le\frac{r}{2}+
\frac{r}{17}<\frac{13r}{16}.\]
Therefore, the first halving loop executes at least once. Suppose first that the halving loop terminates because
$\lambda<\lambda_{\min}$.
Let $\lambda_k$ denote this value and
$\lambda_{k-1}=2\lambda_k$ the preceding iterate.
Since the loop body executed at $\lambda_{k-1}$,
\[\lambda_{k-1}\ge\lambda_{\min},
\qquad \Delta(\lambda_{k-1}) \le
\frac{13r}{16}.\]
Because $\frac{13r}{16}<r-\frac{r}{17}$, Lemma \ref{lemma:nextlambda_helper_lemma} gives $\hat{\Delta}(\lambda_{k-1})
\le \frac{13r}{16}+\frac{r}{17}<r$.
The algorithm returns
\[\lambda_{k-1} =
2\lambda_k<2\lambda_{\min},\]
so the first case holds. For the remainder of the proof, the halving loop terminates because
$\Delta(\lambda_\ell)>\frac{13r}{16}$ while $\lambda_\ell\ge\lambda_{\min}$.
Since the loop executed at the previous iterate
$\lambda_u=2\lambda_\ell$,
$\Delta(\lambda_u)\le\frac{13r}{16}$. If $\Delta(\lambda_\ell)\le
\frac{15r}{16}$, then $\Delta(\lambda_\ell)\in\left(
\frac{13r}{16},\frac{15r}{16}\right]$. Since $\frac{15r}{16}<r-\frac{r}{17}$,
Lemma \ref{lemma:nextlambda_helper_lemma} gives
\[\hat{\Delta}(\lambda_\ell)
\in \left(\frac{13r}{16}-\frac{r}{17},\frac{15r}{16}+\frac{r}{17}
\right]\subset\left[\frac{3r}{4},
r\right).\]
Thus $\hat{\Delta}(\lambda_\ell)<r$.
It remains to consider the binary-search phase.
At the beginning of this phase,
\[\Delta(\lambda_\ell)>\frac{15r}{16},\qquad\Delta(\lambda_u)\le
\frac{13r}{16}.\]
The update rule preserves these two inequalities throughout the
binary search. If the binary search terminates because
\[\Delta(\lambda_m)\in\left[
\frac{13r}{16},\frac{15r}{16}
\right],\]
then, exactly as above, Lemma \ref{lemma:nextlambda_helper_lemma} gives $\hat{\Delta}(\lambda_m)\in\left[
\frac{3r}{4},r\right)$.
Finally, suppose the binary search terminates because
\[\log\frac{\lambda_u}{\lambda_\ell}
\le\frac{r}{8\left(R+\frac{1}{p-1}
\left(\frac{G}{\lambda_\ell}
\right)^{1/(p-1)}\right)}.\]
Then, Lemma \ref{lemma:nextlambda_helper_lemma} gives $\hat{\Delta}(\lambda_u)
\le\frac{13r}{16}+\frac{r}{17}
<\frac{7r}{8}$. We next show that
$\hat{\Delta}(\lambda_\ell)>
\frac{7r}{8}$: Suppose instead that
$\hat{\Delta}(\lambda_\ell)\le
\frac{7r}{8}$. Since $\frac{7r}{8}<r-\frac{r}{17}$, Lemma \ref{lemma:nextlambda_helper_lemma} would imply
\[\Delta(\lambda_\ell)\le\frac{7r}{8}+\frac{r}{17}<\frac{15r}{16},\]
contradicting the binary-search invariant. Therefore,
$\hat{\Delta}(\lambda_u)<\frac{7r}{8}<\hat{\Delta}(\lambda_\ell)$.
By Lemma \ref{lemma:delta_hat_diff_bound},
$\hat{\Delta}$ is continuous. It follows by the Intermediate Value that there exists
$\lambda^*\in(\lambda_\ell,\lambda_u)$ such that $\hat{\Delta}(\lambda^*)=\frac{7r}{8}$. Since
$\lambda_m=\sqrt{\lambda_u\lambda_\ell}$, we have $\left|\log\frac{\lambda_m}{\lambda^*}\right|\le\frac{1}{2}
\log\frac{\lambda_u}{\lambda_\ell}$.
Applying
Lemma \ref{lemma:delta_hat_diff_bound} to
$\lambda_m$ and $\lambda^*$, and using
$\min\{\lambda_m,\lambda^*\}\ge\lambda_\ell$, gives
\begin{align*}
\left|\hat{\Delta}(\lambda_m)-
\hat{\Delta}(\lambda^*)\right|&\le
\left[R+\frac{1}{p-1}\left(\frac{G}{\lambda_\ell}\right)^{1/(p-1)}
\right]\left|
\log\frac{\lambda_m}{\lambda^*}
\right|\\ &\le\frac{1}{2}\left[R+
\frac{1}{p-1}\left(\frac{G}{\lambda_\ell}\right)^{1/(p-1)}
\right]\log\frac{\lambda_u}{\lambda_\ell}\le\frac{r}{16}.
\end{align*}
Since $\hat{\Delta}(\lambda^*)=
\frac{7r}{8}$, it follows that
$\hat{\Delta}(\lambda_m)\in
\left[\frac{13r}{16},\frac{15r}{16}
\right]\subset\left[\frac{3r}{4},r
\right)$. It remains to bound the number of BROO calls. The halving phase clearly requires at most
$O\left(\log\frac{\lambda_{\max}}{\lambda_{\min}}\right)$ calls.
For the binary-search phase, initially $\log\frac{\lambda_u}{\lambda_\ell}=\log 2$. Since
$\lambda_m=\sqrt{\lambda_u\lambda_\ell}$, each bisection halves the logarithmic interval width. Therefore, after $k$ iterations,
$\log\frac{\lambda_u}{\lambda_\ell}
=\frac{\log 2}{2^k}$. The binary search terminates once $$\log\frac{\lambda_u}{\lambda_\ell}
\le \frac{r}{8\left[R+\frac{1}{p-1}
\left(\frac{G}{\lambda_\ell}
\right)^{1/(p-1)}\right]}$$
Throughout the binary search, $\lambda_\ell\ge\lambda_{\min}$.
Therefore,
$\left(\frac{G}{\lambda_\ell}
\right)^{1/(p-1)}\le\left(
\frac{G}{\lambda_{\min}}
\right)^{1/(p-1)}$. Consequently
\[\frac{r}{8\left[R+\frac{1}{p-1}
\left(\frac{G}{\lambda_\ell}\right)^{1/(p-1)}
\right]}\ge\frac{r}{8\left[R+
\frac{1}{p-1}\left(
\frac{G}{\lambda_{\min}}
\right)^{1/(p-1)}\right]}\]
Thus, it is sufficient to make the logarithmic interval width smaller
than the right-hand side above. In particular, the binary search is
guaranteed to have terminated once
\[\frac{\log 2}{2^k}\le\frac{r}{
8\left[R+\frac{1}{p-1}\left(
\frac{G}{\lambda_{\min}}
\right)^{1/(p-1)}\right]} \Rightarrow 2^k\ge\frac{8\log 2}{r}
\left[R+\frac{1}{p-1}\left(
\frac{G}{\lambda_{\min}}
\right)^{1/(p-1)}\right]\]
Therefore, the number of binary-search iterations is at most
\[O\left(\log\left(1+\frac{R+
\frac{1}{p-1}\left(\frac{G}{\lambda_{\min}}\right)^{1/(p-1)}}{
r}\right)\right)\] Combining the halving and binary-search phases, the overall BROO complexity is
\[O\left(\log\frac{\lambda_{\max}}{\lambda_{\min}}+\log\left(1+\frac{
R+\frac{1}{p-1}\left(\frac{G}{\lambda_{\min}}\right)^{1/(p-1)}}{
r}\right)\right).\]
\end{proof}

\subsubsection{$\ell_p$ Proximal Point Method Guarantees}
In the analysis, below, we define the filtration $\mathcal{F}_k = \sigma(x_1, v_1, A_1, \zeta_1, \ldots, x_k, v_k, A_k, \zeta_k)$, where $\zeta_i$ represents the internal randomness in NextLambda($x_i, v_i, A_i$). Furthermore, we let  $\hat{x}_k = P_{f, \lambda_k}(y_{k - 1})$, i.e the exact proximal mapping for which $x_k$ approximates.
\begin{lemma}
\label{lemma:biasedreduced_potentialdecreasebound}
Let $f:\mathcal{X} \rightarrow \mathbb{R}$ be convex. If $\mathcal{X}\subseteq\mathbb{B}_R(x_0)$, it holds that \begin{flalign*}
    &\mathbb{E}[A_{k + 1} (f(x_{k + 1}) - f(x^*)) + 2D_\Phi(x^*, v_{k + 1})] \\ 
        &\le A_k(f(x_k) - f(x^*)) + 2D_\Phi(x^*, v_k) + A_{k + 1} \varphi_{k + 1} - \frac{1 - 2^{1 - q}}{q} \lambda_{k + 1} A_{k + 1} ||\hat{x}_{k + 1} - y_k||_p^p \\ 
        &+ \frac{2^{2 - q} a_{k + 1}^q}{q} \sigma_{k + 1}^q + 2 a_{k + 1} \delta_{k + 1} R 
\end{flalign*}
\end{lemma}

\begin{proof} 
 Following a similar approach to \cite{asi2021stochastic}, we let $\hat{g}_k = \nabla f_{\lambda_k}(y_{k - 1})$ and we bound from above and below the quantity $a_{k + 1} \langle \hat{g}_{k + 1}, v_k - x^* \rangle$. First, observe that \begin{flalign*}
        &v_k - x^* = \hat{x}_{k + 1} - x^* + \frac{A_k}{a_{k + 1}}(\hat{x}_{k + 1} - x_k) - \frac{A_{k + 1}}{a_{k + 1}}(\hat{x}_{k + 1} - y_k)
\end{flalign*} which implies that \begin{flalign*}
    \langle \hat{g}_{k + 1}, v_k - x^* \rangle &= \langle \hat{g}_{k + 1}, \hat{x}_{k + 1} - x^* \rangle + \frac{A_k}{a_{k + 1}} \langle \hat{g}_{k + 1}, \hat{x}_{k + 1} - x_k \rangle -\frac{A_{k + 1}}{a_{k + 1}} \langle \hat{g}_{k + 1}, \hat{x}_{k + 1} - y_k \rangle \\ 
    &\ge f(\hat{x}_{k + 1}) - f(x^*) + \frac{A_k}{a_{k + 1}} (f(\hat{x}_{k + 1})- f(x_k)) + \frac{\lambda_{k + 1}A_{k + 1}}{a_{k + 1}} ||\hat{x}_{k + 1} - y_k||_p^p \\ 
    &= f(\hat{x}_{k + 1}) - f(x^*)  + \frac{A_k}{a_{k + 1}}(f(\hat{x}_{k + 1}) - f(x^*) - (f(x_k) - f(x^*))) + \frac{\lambda_{k + 1}A_{k + 1}}{a_{k + 1}} ||\hat{x}_{k + 1} - y_k||_p^p \\ 
    &= \frac{A_{k + 1}}{a_{k + 1}}(f(\hat{x}_{k + 1}) - f(x^*)) - \frac{A_k}{a_{k + 1}}(f(x_k) - f(x^*)) + \frac{\lambda_{k + 1}A_{k + 1}}{a_{k + 1}} ||\hat{x}_{k + 1} - y_k||_p^p
\end{flalign*}  where we used the fact that $f$ is convex and that $\langle\widehat{g}_{k+1},\widehat{x}_{k+1}-y_k\rangle
=-\lambda_{k+1}\|\widehat{x}_{k+1}-y_k\|_p^p$. We also have that \begin{flalign*}
        &\mathbb{E}[f(x_{k + 1})|\mathcal{F}_k] \le \mathbb{E}[f(x_{k + 1}) + \frac{\lambda_{k + 1}}{p} ||x_{k + 1} - y_k||_p^p] \le f(\hat{x}_{k + 1}) + \frac{\lambda_{k + 1}}{p} ||\hat{x}_{k + 1} - y_k||_p^p + \varphi_{k + 1} 
\end{flalign*} which implies that \begin{flalign*}
    &a_{k + 1}\langle \hat{g}_{k + 1}, v_k - x^* \rangle \\ 
    &\ge A_{k + 1} (\mathbb{E}[f(x_{k + 1})| \mathcal{F}_k] - \frac{\lambda_{k + 1}}{p}||\hat{x}_{k + 1} - y_k||_p^p - \varphi_{k + 1} - f(x^*)) - A_k(f(x_k) - f(x^*)) \\ 
    &+ \lambda_{k + 1} A_{k + 1}||\hat{x}_{k + 1} - y_k||_p^p \\ 
    &= A_{k + 1} (\mathbb{E}[f(x_{k + 1}) | \mathcal{F}_k] - f(x^*)) - A_k(f(x_k) - f(x^*)) + (1 - \frac{1}{p}) \lambda_{k + 1} A_{k + 1} ||\hat{x}_{k + 1} - y_k||_p^p - A_{k + 1} \varphi_{k + 1}
\end{flalign*} We now upper bound $a_{k + 1}\langle \hat{g}_{k + 1}, v_k - x^* \rangle$. We have that \begin{flalign*}
    &v_{k + 1} = \argmin_{v \in X} \{\frac{1}{2}a_{k + 1}\langle g_{k + 1}, v \rangle  + D_\Phi(v, v_k) \}  = \argmin_{v \in X} \{\frac{1}{2}a_{k + 1} \langle g_{k + 1}, v \rangle + \Phi(v) - \langle \nabla \Phi(v_k), v \rangle \}
\end{flalign*} By first-order optimality conditions, it holds that $0\in\tfrac12 a_{k+1}g_{k+1}
+\nabla\Phi(v_{k+1})-\nabla\Phi(v_k)+N_{\mathcal X}(v_{k+1})$, where $N_X$ denotes the normal cone, which implies that for all $u \in X$, \begin{flalign*}
    &\langle \frac{1}{2} a_{k + 1} g_{k + 1} + \nabla \Phi(v_{k + 1}) - \nabla \Phi(v_k), u - v_{k + 1} \rangle \ge 0. 
\end{flalign*} Therefore, it holds that \begin{flalign*}
    \frac{1}{2}a_{k + 1} \langle g_{k + 1}, v_{k + 1} - x^* \rangle \le \langle \nabla \Phi(v_{k + 1}) - \nabla \Phi(v_{k}), x^* - v_{k + 1} \rangle 
\end{flalign*} Since $\langle \nabla \Phi(v_{k + 1}) - \nabla \Phi(v_{k}), x^* - v_{k + 1} \rangle = D_\Phi(x^*, v_{k}) - D_\Phi(x^*, v_{k + 1}) - D_\Phi(v_{k + 1}, v_k)$, it holds that \begin{flalign*}
    \frac{1}{2}a_{k + 1} \langle g_{k + 1}, v_{k + 1} - x^* \rangle  \le D_\Phi(x^*, v_{k}) - D_\Phi(x^*, v_{k + 1}) - D_\Phi(v_{k + 1}, v_k)
\end{flalign*} 
Noting that $v_k - x^* = (v_k - v_{k + 1}) + (v_{k + 1} - x^*)$, we have that
\begin{flalign*}
    a_{k + 1}\langle \hat{g}_{k + 1}, v_k - x^* \rangle
    &= a_{k + 1}\langle g_{k + 1}, v_{k + 1} - x^* \rangle \\
    &\quad + a_{k + 1}\langle
    \hat{g}_{k + 1} + g_{k + 1}
    - \mathbb{E}[g_{k + 1}|\mathcal{F}_k],
    v_k - v_{k + 1} \rangle \\
    &\quad + a_{k + 1}\langle
    \mathbb{E}[g_{k + 1}|\mathcal{F}_k]
    - \hat{g}_{k + 1},
    x^* - v_{k + 1} \rangle \\
    &\quad - a_{k + 1}\langle
    g_{k + 1}
    - \mathbb{E}[g_{k + 1}|\mathcal{F}_k],
    v_k - x^* \rangle.
\end{flalign*}
By H\"older's inequality, we have that for all
$a,b\in\mathbb{R}^d$ and $\eta>0$,
$\langle a,b\rangle
\le \frac{\eta^{1-q}}{q}\|a\|_q^q
+\frac{\eta}{p}\|b\|_p^p$.
By choosing $\eta=2^{3-p}$, it holds that
\begin{flalign*}
    a_{k + 1}\langle
    \hat{g}_{k + 1} + g_{k + 1}
    - \mathbb{E}[g_{k + 1}|\mathcal{F}_k],
    v_k - v_{k + 1}\rangle
    &\le \frac{2^{3-2q}a_{k + 1}^q}{q}\|\hat{g}_{k + 1} + g_{k + 1}
    - \mathbb{E}[g_{k + 1}|\mathcal{F}_k]\|_q^q \\
    &\quad + \frac{2^{3-p}}{p}
    \|v_k-v_{k + 1}\|_p^p.
\end{flalign*}
By Fact \ref{fact:bregmanlowerbound},
we have that
$-2D_\Phi(v_{k + 1},v_k)
+\frac{2^{3-p}}{p}\|v_k-v_{k + 1}\|_p^p
\le 0$, which implies that
\begin{flalign*}
    a_{k + 1}\langle
    \hat{g}_{k + 1},v_k-x^*\rangle &\le 2D_\Phi(x^*,v_k)
    -2\mathbb{E}[D_\Phi(x^*,v_{k + 1})|\mathcal{F}_k] 
    \\
    &\quad+\frac{2^{3-2q}a_{k + 1}^q}{q}
    \mathbb{E}\big[
    \|\hat{g}_{k + 1}+g_{k + 1}
    -\mathbb{E}[g_{k + 1}|\mathcal{F}_k]\|_q^q
    \big|\mathcal{F}_k\big] \\
    &\quad +a_{k + 1}\mathbb{E}\big[
    \langle \mathbb{E}[g_{k + 1}|\mathcal{F}_k]
    -\hat{g}_{k + 1},
    x^*-v_{k + 1}\rangle
    \big|\mathcal{F}_k\big],
\end{flalign*}
where we used the fact that $v_k$ and $x^*$
are $\mathcal{F}_k$-measurable and therefore
\begin{flalign*}
    &\mathbb{E}\big[
    \langle g_{k + 1}
    -\mathbb{E}[g_{k + 1}|\mathcal{F}_k],
    v_k-x^*\rangle
    \big|\mathcal{F}_k\big]=0.
\end{flalign*}
It also holds that
\begin{flalign*}
    \mathbb{E}\big[
    \|\hat{g}_{k + 1}+g_{k + 1}
    -\mathbb{E}[g_{k + 1}|\mathcal{F}_k]\|_q^q
    \big|\mathcal{F}_k\big]
    &\le 2^{q-1}\|\hat{g}_{k + 1}\|_q^q  +2^{q-1}\mathbb{E}\big[
    \|g_{k + 1}
    -\mathbb{E}[g_{k + 1}|\mathcal{F}_k]\|_q^q
    \big|\mathcal{F}_k\big] \\
    &\le 2^{q-1}\sigma_{k + 1}^q
    +2^{q-1}\|\hat{g}_{k + 1}\|_q^q.
\end{flalign*}
Moreover, by the bias guarantee and
the fact that $\mathcal{X}\subseteq
\mathbb{B}_R(x_0)$, we have
\begin{flalign*}
    \mathbb{E}\big[
    \langle
    \mathbb{E}[g_{k + 1}|\mathcal{F}_k]
    -\hat{g}_{k + 1},
    x^*-v_{k + 1}\rangle
    \big|\mathcal{F}_k\big]
    &\le
    \|\mathbb{E}[g_{k + 1}|\mathcal{F}_k]
    -\hat{g}_{k + 1}\|_q
    \mathbb{E}[\|x^*-v_{k + 1}\|_p
    |\mathcal{F}_k] \\
    &\le 2\delta_{k + 1}R.
\end{flalign*}
Therefore, we have that
\begin{flalign*}
    a_{k + 1}\langle
    \hat{g}_{k + 1},v_k-x^*\rangle 
    &\le 2D_\Phi(x^*,v_k)
    -2\mathbb{E}[D_\Phi(x^*,v_{k + 1})
    |\mathcal{F}_k] \\
    &\quad +\frac{2^{2-q}a_{k + 1}^q}{q}
    \big(\sigma_{k + 1}^q
    +\|\hat{g}_{k + 1}\|_q^q\big)
    +2a_{k + 1}\delta_{k + 1}R \\
    &=2D_\Phi(x^*,v_k)
    -2\mathbb{E}[D_\Phi(x^*,v_{k + 1})
    |\mathcal{F}_k] \\
    &\quad +\frac{2^{2-q}a_{k + 1}^q}{q}
    \big(\sigma_{k + 1}^q
    +\lambda_{k + 1}^q
    \|\hat{x}_{k + 1}-y_k\|_p^p\big)
    +2a_{k + 1}\delta_{k + 1}R.
\end{flalign*}
This implies that
\begin{flalign*}
    &A_{k+1}(\mathbb E[f(x_{k+1})\mid\mathcal F_k]-f(x^*))-A_k(f(x_k)-f(x^*))+\frac1q\lambda_{k+1}A_{k+1}||\hat x_{k+1}-y_k||_p^p-A_{k+1}\varphi_{k+1}\\
    &\le 2D_\Phi(x^*,v_k)-2\mathbb E[D_\Phi(x^*,v_{k+1})\mid\mathcal F_k]
    +\frac{2^{2-q}a_{k+1}^q}{q}(\sigma_{k+1}^q+\lambda_{k+1}^q||\hat x_{k+1}-y_k||_p^p)
    +2a_{k+1}\delta_{k+1}R.
\end{flalign*}
This implies that
\begin{flalign*}
    &\mathbb E[A_{k+1}(f(x_{k+1})-f(x^*))+2D_\Phi(x^*,v_{k+1})\mid\mathcal F_k]\\
    &\le A_k(f(x_k)-f(x^*))+2D_\Phi(x^*,v_k)-\frac1q\lambda_{k+1}A_{k+1}||\hat x_{k+1}-y_k||_p^p\\
    &\quad+A_{k+1}\varphi_{k+1}+\frac{2^{2-q}a_{k+1}^q}{q}(\sigma_{k+1}^q+\lambda_{k+1}^q||\hat x_{k+1}-y_k||_p^p)
    +2a_{k+1}\delta_{k+1}R\\
    &\le A_k(f(x_k)-f(x^*))+2D_\Phi(x^*,v_k)
    -\frac{1-2^{1-q}}q\lambda_{k+1}A_{k+1}||\hat x_{k+1}-y_k||_p^p
    +A_{k+1}\varphi_{k+1}\\
    &\quad+\frac{2^{2-q}a_{k+1}^q}{q}\sigma_{k+1}^q+2a_{k+1}\delta_{k+1}R,
\end{flalign*}
where we used $A_{k+1}^{p-1}=2^{p-1}\lambda_{k+1}a_{k+1}^p$ in the last inequality.
\end{proof}

\begin{lemma}
\label{lemma:biasedreducedK_maxbound}
 Let $K \le K_{\max}$ be the iteration in which the algorithm stops and let \begin{flalign*}
    &\overline{\varepsilon} \geq \max_{k < K_{\max}} \{\frac{A_{k + 1} \varphi_{k + 1}}{a_{k + 1}} + \frac{2^{2-q}a_{k+1}^{q-1} \sigma_{k + 1}^q}{q} + 2 \delta_{k + 1}R \}
\end{flalign*}  We then have that \begin{flalign*}
    &\mathbb{E}\left[A_K(f(x_K) - f(x^*) - \overline{\varepsilon}) + 2D_\Phi(x^*, v_K) + \frac{1 - 2^{1 - q}}{q} \sum_{i = 1}^K \lambda_i A_i ||\hat{x}_i - y_{i - 1} ||_p^p\right] \\ 
    &\le A_0(f(x_0) - f(x^*)) + 3^{p - 2}(p - 1)R^p
\end{flalign*} 
\end{lemma}

\begin{proof}
    Let \begin{flalign*}
        &T_k := -\frac{1 - 2^{1 - q}}{q} \lambda_{k + 1} A_{k + 1} || \hat{x}_{k + 1} - y_k||_p^p + A_{k + 1} \varphi_{k + 1} + \frac{2^{2 - q} a_{k + 1}^q}{q} \sigma_{k + 1}^q + 2a_{k + 1} \delta_{k + 1} R.
    \end{flalign*}
    First, we have from Lemma \ref{lemma:biasedreduced_potentialdecreasebound} that 
    \begin{flalign*}
        &\mathbb{E}[A_{k + 1}(f(x_{k + 1} - f(x^*)) + 2D_\Phi(x^*, v_{k + 1})] \le A_k(f(x_k)- f(x^*)) + 2D_\Phi(x^*, v_k) + T_k. 
    \end{flalign*}  
    Define \begin{flalign*}
        M_k = A_k(f(x_k) - f(x^*) - \overline{\varepsilon}) + 2D_\Phi(x^*, v_k) + \frac{1 - 2^{1 - q}}{q} \sum_{i = 1}^k \lambda_i A_i ||\hat{x}_i - y_{i - 1} ||_p^p.
    \end{flalign*} We will now prove that $M_k$ is a supermartingale. It holds that \begin{flalign*}
        \mathbb{E}[M_{k + 1} | \mathcal{F}_k]  
        &= \mathbb{E}[A_{k + 1}(f(x_{k + 1}) - f(x^*) - \overline{\eps}) + 2D_\Phi(x^*, v_{k + 1})] + \frac{1 - 2^{1 - q}}{q} \sum_{i = 1}^k \lambda_i A_i ||\hat{x}_i - y_{i - 1}||_p^p \\ 
        &+ \frac{1 - 2^{1 - q}}{q} \lambda_{k + 1} A_{k + 1} ||\hat{x}_{k + 1} - y_k||_p^p \\ 
        &\le A_k(f(x_k) - f(x^*)) + 2D_\Phi(x^*, v_k) + T_k - A_{k + 1} \overline{\eps} + \frac{1 - 2^{1 - q}}{q} \sum_{i = 1}^k \lambda_i A_i ||\hat{x}_i - y_{i - 1}||_p^p \\ 
        &+ \frac{1 - 2^{1 - q}}{q} \lambda_{k + 1} A_{k + 1} ||\hat{x}_{k + 1} - y_k||_p^p \\ 
        &\le A_k(f(x_k) - f(x^*)) + 2D_\Phi(x^*, v_k) + T_k - (A_{k } + a_{k + 1}) \overline{\eps} + \frac{1 - 2^{1 - q}}{q} \sum_{i = 1}^k \lambda_i A_i ||\hat{x}_i - y_{i - 1}||_p^p \\ 
        &+ \frac{1 - 2^{1 - q}}{q} \lambda_{k + 1} A_{k + 1} ||\hat{x}_{k + 1} - y_k||_p^p \\ 
        &\le M_k + T_k + \frac{1 - 2^{1 - q}}{q} \lambda_{k + 1} A_{k + 1} ||\hat{x}_{k + 1} - y_k||_p^p - a_{k + 1} \overline{\eps} \\ 
        &= M_k + A_{k + 1} \varphi_{k + 1} + \frac{2^{2 - q} a_{k + 1}^q \sigma_{k + 1}^q}{q} + 2a_{k + 1} \delta_{k + 1} R - a_{k + 1} \overline{\eps}
    \end{flalign*} By our setting of $\overline{\varepsilon}$, we have that \begin{flalign*}
        &M_k + A_{k + 1} \varphi_{k + 1} + \frac{2^{2 - q} a_{k + 1}^q}{q} \sigma_{k + 1}^q + 2 a_{k + 1} \delta_{k + 1} R - a_{k + 1} \overline{\varepsilon} \le M_k
    \end{flalign*} so $M_k$ is a supermartingale. Therefore, by the optional stopping theorem, it holds that \begin{flalign*}
        \mathbb{E}[M_K] \le M_0 &= A_0(f(x_0) - f(x^*) - \overline\varepsilon) + 2D_\Phi(x^*, v_0)
        \end{flalign*}  Using that $v_0 = x_0$, we have that \begin{flalign*}
        D_\Phi(x^*, v_0) = D_\Phi(x^*, x_0) &\le \frac{p - 1}{2} (||x^*||_p + 2||x_0||_p)^{p - 2} ||x^* - x_0||_p^2 \\ 
        &\le \frac{p - 1}{2} (3||x_0||_p + ||x^* - x_0||_p)^{p - 2} R^2 \\ 
        &\le \frac{3^{p - 2}(p - 1)}{2} R^p
    \end{flalign*} which finishes the proof. 
\end{proof}

\begin{fact}(Lemma 32 of \cite{carmon2020acceleration})
\label{fact:growthlemma}
Let $\{B_k \}_{k \in \mathbb{N}}$ be a nonnegative, nondecreasing sequence such that $B_k \ge \sum_{i \in [k]} \alpha B_i$ for some $\alpha \in [0, 1)$ and all $k$. Then for all $k$, $B_k \ge \exp(\alpha(k - 1)) B_1$. 
    
\end{fact}

\begin{lemma}
Fix any $k>0$. Define
\[T_\lambda:=\sum_{i\le k}\mathbf{1}_{\{\lambda_i<2\lambda_{\min}\}},
\quad T_r:=\sum_{i\le k}\mathbf{1}_{\{\|\hat{x}_i-y_{i-1}\|\geq 3r/4\}}\]
Then, with probability one, for $\eps < \frac{2^{q - 1} - 1}{3\cdot 2^{q - 1} - 1}$, 
\[\frac{1}{A_k}\left(3C_pR^p-\frac{(1 - 2^{1 - q})}{q}\sum_{i\le k}\lambda_i A_i\|\hat{x}_i-y_{i-1}\|^p_p\right)\leq O\left(\min\left\{\frac{\lambda_{\min}R^p}{T_\lambda^p},\,\frac{R^p}{A_0}
\exp\left(-\Omega(1) \frac{r^{p/(p+1)}}{R^{p/(p+1)}}T_r\right)\right\}\right)\] See \eqref{defn:cplabel} for a definition of $C_p$. 
\end{lemma} 

\begin{proof}
    If $3C_pR^p- \frac{(1 - 2^{1 - q})}{q}\sum_{i\le k}\lambda_i A_i\|\hat{x}_i-y_{i-1}\|^p_p \le 0$, then the desired inequality trivially holds. Therefore, we only consider the case where $\frac{(1 - 2^{1 - q})}{q} \sum_{i\le k}\lambda_i A_i\|\hat{x}_i-y_{i-1}\|^p_p  \le 3C_pR^p$. We first prove that \begin{flalign}
        &A_t^{1/p} \ge \frac{1}{2p} \sum_{i = 1}^{t} \lambda_i^{-1/p}. \label{ineq:At_recurrence}
    \end{flalign} Since $A_i^{p - 1} = 2^{p - 1} \lambda_i a_i^p$, it holds that $a_i = \frac{A_i^{(p - 1)/p}}{2^{(p - 1)/p} \lambda_i^{1/p}}$. Therefore, \begin{flalign*}
        &a_i \ge \frac{1}{2^{(p - 1)/p}}A_i^{(p - 1)/p} \lambda_i^{-1/p} \ge  \frac{1}{2} A_i^{(p - 1)/p} \lambda_i^{-1/p}
    \end{flalign*} Since $f(x) = x^{1/p}$ is concave for $p \ge 2$, it holds that \begin{flalign}
        &A_i^{1/p} - A_{i - 1}^{1/p} \ge \frac{1}{p}A_i^{\frac{1}{p} - 1}(A_i - A_{i - 1}) = \frac{a_i}{p A_i^{(p - 1)/p}} \ge \frac{\lambda_i^{-1/p}}{2p} \label{ineq:onestep}
    \end{flalign} Telescoping \eqref{ineq:onestep} from $i = 1$ to $t$ and using that $A_0 \ge 0$ yields \eqref{ineq:At_recurrence}. For $t \le k$, define \begin{flalign*}
        &I_{r, t} = \left\{ i \le t: ||\hat{x}_i - y_{i - 1}||_p \ge \frac{3r}{4} \right\}
    \end{flalign*} It holds that \begin{flalign*}
        &\frac{1 - 2^{1 - q}}{q} \sum_{i \in I_{r, t}} \lambda_i A_i(\frac{3r}{4})^p \le \frac{1 - 2^{1 - q}}{q} \sum_{i \le t} \lambda_i A_i ||\hat{x}_i - y_{i - 1}||_p^p \le 3C_pR^p \implies \sum_{i \in I_{r, t}} \lambda_i A_i \le 3C_p(\frac{R}{r})^p (\frac{4}{3})^p \cdot \frac{q}{1 - 2^{1 - q}} 
    \end{flalign*} Now, by a reverse Hölder inequality, it holds that for all $\theta \in (0, 1)$, \begin{flalign*}
        &\sum_i u_i v_i \ge (\sum_i u_i^\theta)^{1/\theta} \cdot (\sum_i v_i^{-\theta/(1 - \theta)})^{-(1 - \theta)/\theta}
    \end{flalign*} Taking $u_i = A_i^{1/p} , v_i = \frac{1}{\lambda_i^{1/p} A_i^{1/p}}, \theta = p/(p + 1)$, it holds that \begin{flalign*}
        \sum_{i \in I_{r, t}} \frac{1}{\lambda_i^{1/p}} \ge (\sum_{i \in I_{r, t}} A_i^{1/(p + 1)})^{(p + 1)/p} \cdot (\sum_{i \in I_{r, t}} \lambda_i A_i)^{-1/p} &\ge (3C_p \cdot \frac{q}{1 - 2^{1 - q}} \cdot (\frac{4}{3})^p)^{-1/p} \cdot \frac{r}{R} \cdot (\sum_{i \in I_{r, t}} A_i^{1/(p + 1)})^{(p + 1)/p} \\ 
        &\ge C_p' \cdot \frac{r}{R} \cdot (\sum_{i \in I_{r, t}} A_i^{1/(p + 1)})^{(p + 1)/p}
    \end{flalign*} where we let $C_p' = (3C_p \cdot \frac{q}{1 - 2^{1 - q}} \cdot (\frac{4}{3})^p)^{-1/p}$. Combining with \eqref{ineq:At_recurrence} gives that \begin{flalign*}
        &A_t^{1/(p + 1)} \ge (\frac{C_p'}{2p})^{p/(p + 1)} \cdot (\frac{r}{R})^{\frac{p}{p + 1}} \sum_{i \in I_{r, t}} A_i^{1/(p + 1)} 
    \end{flalign*} By Fact \ref{fact:growthlemma}, it holds that \begin{flalign}
        &A_k^{1/(p + 1)} \ge A_0^{1/(p + 1)} \exp(\Omega(1) (\frac{r}{R})^{\frac{p}{p + 1}} T_r) \implies \frac{1}{A_k} \le \frac{1}{A_0} \exp(-\Omega(1) (\frac{r}{R})^{\frac{p}{p + 1}}T_r) \label{ineq:growthlemma_consequence}
    \end{flalign} Therefore, \begin{flalign}
        &\frac{1}{A_k}\left(3C_pR^p-\frac{1 - 2^{1 - q}}{q}\sum_{i\le k}\lambda_i A_i\|\hat{x}_i-y_{i-1}\|^p_p\right) \le \frac{3C_pR^p}{A_k} \le O(\frac{R^p}{A_0} \exp(-\Omega(1) (\frac{r}{R})^{\frac{p}{p + 1}}T_r)) \label{bound1}
    \end{flalign} Let $I_{\lambda} = \{i \le k: \lambda_i < 2\lambda_{\min} \}$. Combining \eqref{ineq:At_recurrence} with $\lambda_i < 2\lambda_{\min}$ yields \begin{flalign}
        &A_k^{1/p} \ge \frac{1}{2p} \cdot \frac{T_{\lambda}}{(2\lambda_{\min})^{1/p}} \implies \frac{1}{A_k} \le O(\frac{\lambda_{\min}}{T_{\lambda}^p}) \implies \frac{3C_pR^p}{A_k} \le O(\frac{R^p\lambda_{\min}}{T_{\lambda}^p}) \label{bound2}
    \end{flalign} Combining \eqref{bound1} and \eqref{bound2} finishes the proof. 
\end{proof}

\begin{proposition}
\label{proposition:stochastic_prox_method_guarantee}
Let $f : \mathcal{X} \to \mathbb{R}$ and $\widehat{\nabla} f$ satisfy Assumption \ref{Assumption}, and let $\mathcal{X} \subseteq B_R(x_0)$. For a target
accuracy $\varepsilon \le GR$ let
\[\varphi_k =
\frac{\varepsilon}{48C_p \lambda_k^{q - 1} a_k^{q - 1}},
\quad \delta_k =
\frac{\varepsilon}{96C_p R},
\quad \sigma_k^q =
\frac{q\varepsilon}{48C_p2^{2-q}a_k^{q-1}}, \quad A_0 =
\frac{R^{p-1}}{G}, \quad A_{\max} =
3C_p\frac{R^p}{\varepsilon}\] see \eqref{defn:cplabel}.
If $\lambda_k \ge \lambda_{\min} \ge
\frac{1}{A_{\max}}
= \Omega\left(\frac{\varepsilon}{R^p}\right)$ for all $k \le K_{\max}$, then lines 7 and 9 of Algorithm \ref{alg:stochastic_accelerated_prox_point} have total complexity \[
\mathbb{E}N_{\widehat{\nabla} f}
=O\left(K_{\max}\log\left(\frac{GR}{\varepsilon}\right)+ \left(\frac{GR}{\varepsilon}\right)^{p(p-1)}
\log\!\left(\frac{GR}{\varepsilon}\right)\right)\] If in addition $\|P_{f,\lambda_k}(y_{k-1}) - y_{k-1}\|_p \ge 3r/4$ whenever $\lambda_k \ge 2\lambda_{\min}$ then for
\[K_{\max}=O\left(\left(\frac{R}{r}\right)^{\frac{p}{p+1}}
\log\left(\frac{GR}{\varepsilon}\right)+\bigg(\frac{\lambda_{\min}R^p}{\varepsilon}\bigg)^{\frac{1}{p}}\right)\]
the algorithm’s output $x_K$ satisfies $f(x_K) - f(x^\star) \le \varepsilon$ with probability at least $2/3$.
\end{proposition}

\begin{proof}
    First, let us prove correctness of the algorithm. The settings of $\varphi_k$, $\delta_k$ and $\sigma_k$ in the proposition
guarantee that \[\max_{k < K_{\max}} \Bigg\{ \frac{A_{k + 1} \varphi_{k + 1}}{a_{k + 1}} + \frac{2^{2-q}a_{k + 1}^{q - 1} \sigma_{k + 1}^q}{q} + 2\delta_{k + 1}R\Bigg\} = \frac{\eps}{12C_p}\] Therefore, applying Lemma \ref{lemma:biasedreducedK_maxbound} with $\overline{\eps} = \eps/(12C_p) \le \eps/(6C_p) \le R^p/(2A_{\max})$ yields \[\mathbb{E}\left[A_K(f(x_K) - f(x^\star)) + \frac{1 - 2^{1 - q}}{q}
\sum_{i\le K}\lambda_i
A_i\|\hat{x}_i - y_{i-1}\|^p_p\right]
\leq 3^{p - 2} R^p(p -1 ) + \bar{\varepsilon} \cdot \mathbb{E}A_K + A_0(f(x_0) - f(x^\star))\] We now upper bound $a_K$ in terms of $A_{K - 1}$ and $\lambda_K$. It holds that $(A_{K - 1} + a_K)^{p -1 } = 2^{p - 1} \lambda_K a_K^p$. If $a_K > A_{K - 1}$, then $A_{K - 1} + a_K < 2a_K$, meaning that $a_K < 1/\lambda_K$ because $2^{p - 1} \lambda_K a_K^p = (A_{K - 1} + a_K)^{p -1 } < (2a_K)^{p - 1}$. Combining this with the $a_K < A_{K - 1}$ case, it holds that \begin{flalign*}
    &a_K \le \max \{A_{K - 1}, \frac{1}{\lambda_K} \} \implies a_K \le A_{K - 1} + \frac{1}{\lambda_K} \implies A_K \le 2A_{K - 1} + \frac{1}{\lambda_K}
\end{flalign*} Since $A_{K - 1} \le A_{\max}$ and $\lambda_K \ge 1/A_{\max}$, it holds that $A_K \le 3 A_{\max}$ with probability 1. Since $A_0 = R^{p - 1}/G$ and $f$ is G-Lipschitz, it holds that $A_0(f(x_0) - f(x^*)) \le A_0 G \cdot ||x_0 - x^*||_p \le A_0 GR =  R^p$. Moreover, using the definition of $A_{\max}$, we have that $\overline{\eps} A_K \le \frac{\eps}{12C_p} \cdot 3(3C_p \frac{R^p}{\eps}) = \frac{3}{4} R^p$. 
Thus, \begin{flalign*}
    \mathbb{E}[A_K(f(x_K) - f(x^*)) + \frac{1 - 2^{1 - q}}{q} \sum_{i \le K} \lambda_i A_i ||\hat{x}_i - y_{i - 1}||_p^p] &\le A_0(f(x_0) - f(x^*)) + 3^{p - 2} R^p(p - 1) + \overline{\eps} \mathbb{E}[A_K] \\ 
    &\le [3^{p - 2}(p - 1) + 1 + \frac{3}{4}] R^p
\end{flalign*} Denote \begin{flalign}
    &C_p = 3^{p - 2}(p - 1) + \frac{7}{4} \label{defn:cplabel}
\end{flalign} Therefore, with probability at least $2/3$, \begin{flalign*}
    &\mathbb{E}[A_K(f(x_K) - f(x^*)) + \frac{1 - 2^{1 - q}}{q} \sum_{i \le K} \lambda_i A_i ||\hat{x}_i - y_{i - 1}||_p^p] \le 3C_pR^p \implies f(x_K) - f(x^*) \le \frac{3C_p R^p}{A_K}
\end{flalign*} If $A_K \ge A_{\max}$, then we are done. Otherwise, $K = K_{\max}$ and by the assumption on NextLambda, it holds that $T_{\lambda} + T_r \ge K_{\max}$, which means that either $T_{\lambda} \ge K_{\max}/2$ or $T_r \ge K_{\max} /2$. If $T_{\lambda} \ge K_{\max}/2$, using \eqref{ineq:At_recurrence} for the $T_{\lambda}$ iterations where $\lambda_i < 2\lambda_{\min}$, we get that \begin{flalign*}
    &A_K^{1/p} \ge \frac{T_{\lambda}}{2p (2\lambda_{\min})^{1/p}} \implies \frac{3C_p R^p}{A_K} \le O\left(\frac{\lambda_{\min} R^p}{T_{\lambda}^p}\right) \le \eps
\end{flalign*} as long as $T_{\lambda} \ge \Omega((\frac{\lambda_{\min}R^p}{\eps})^{1/p})$. This implies that it suffices to let $K_{\max} \ge \Omega((\frac{\lambda_{\min}R^p}{\eps})^{1/p})$. If $T_r \ge K_{\max}/2$, then we can use that $A_K \ge A_0 \exp(\Omega(1)  (\frac{r}{R})^{p/(p + 1)} T_r)$ (from \eqref{ineq:growthlemma_consequence}) and that $R^p/A_0 = GR$ to say that \begin{flalign*}
    &\frac{3C_p R^p}{A_K} \le 3C_p GR\exp(-\Omega(1) (\frac{r}{R})^{p/(p + 1)}T_r)
\end{flalign*} which is at most $\eps$ provided that $T_r = \Omega((\frac{R}{r})^{p/(p + 1)} \log(\frac{GR}{\eps}))$. Therefore, taking \[K_{\max} = O\left(\left(\frac{R}{r}\right)^{p/(p + 1)} \log\left(\frac{GR}{\eps}\right) + \left(\frac{\lambda_{\min}R^p}{\eps}\right)^{1/p}\right)\] suffices to conclude that $f(x_K) - f(x^*) \le \eps$ with probability at least $2/3$. Now, we prove the complexity guarantee. Each step of computing $x_{k + 1}$ on line 7 of Algorithm \ref{alg:stochastic_accelerated_prox_point} requires $O(\frac{G^q}{\lambda_k \phi_k^{p-1}})$ queries by Definition \ref{defn:odc_alg} and the definition of an approximate proximal mapping. Moreover, by Theorem \ref{thm:MorGradEst} it requires \begin{flalign*}
    &O\left(\frac{G^p}{\sigma^p}\left(\frac{G}{\min \{\delta, \sigma \}}\right)^{p(p - 2)} \left(1 + \log \frac{G}{\min \{\delta, \sigma \}}\right)\right)
\end{flalign*} 
in expectation to compute $g_{k+1}$ in line 8 of Algorithm \ref{alg:stochastic_accelerated_prox_point}. 
Observe that \begin{flalign*}
    &\sigma_k^p = (\sigma_k^q)^{p - 1} = O(\frac{\varepsilon^{p - 1}}{a_k^{(q - 1)(p - 1)}}) = O(\frac{\varepsilon^{p - 1}}{a_k}) \implies \frac{G^p}{\sigma_k^p} = O(\frac{G^p a_k}{\varepsilon^{p - 1}}). 
\end{flalign*} Therefore, plugging in $\phi_k$,$\delta_k$ and $\sigma_k$, we have 
\begin{flalign*}
    \mathbb{E}N_{\widehat{\nabla} f} &= \sum_{k\in[K]}O\left(\frac{G^p}{\lambda_k\varphi_k^{p-1}}\right)+\sum_{k\in[K]} O([1 + \frac{G^p}{\sigma_k^p}(\frac{G}{\min \{\delta, \sigma_k \}})^{p(p - 2)}] \cdot \log \frac{G}{\min \{\delta, \sigma_k \}}) \\
    &\le O(\frac{G^p}{\varepsilon^{p - 1}} \sum_{k \in [K]} a_k) + O(\sum_{k \in [K]} \log \frac{G}{\min \{\delta, \sigma_k \}} + \frac{G^p}{\sigma_k^p} (\frac{G}{\min \{\delta, \sigma_k \}})^{p(p - 2)} \log \frac{G}{\min \{\delta, \sigma_k \}}) \\ 
    &\le O(\frac{G^p}{\varepsilon^{p - 1}} \sum_{k \in [K]} a_k) + O(\sum_{k \in [K]} \log \frac{GR}{\varepsilon} + \frac{G^p a_k}{\varepsilon^{p - 1}} (\frac{GR}{\varepsilon})^{p(p - 2)} \log \frac{GR}{\varepsilon}) \\ 
    &\le O(\frac{G^p}{\varepsilon^{p - 1}} A_{\max}) + O(K_{\max} \log \frac{GR}{\varepsilon}) + O(\frac{G^p A_{\max}}{\varepsilon^{p - 1}} (\frac{GR}{\varepsilon})^{p(p - 2)} \log \frac{GR}{\varepsilon}) \\ 
    &\le O((\frac{GR}{\varepsilon})^p ) + O(K_{\max} \log \frac{GR}{\varepsilon}) + O((\frac{GR}{\varepsilon})^{p(p - 1)} \log \frac{GR}{\varepsilon})
\end{flalign*}

since $a_k \le A_K \le 3A_{\max} \le O(R^p/\varepsilon)$,  $\min \{\delta_k, \sigma_k \} = \Omega(\varepsilon/R)$, and $\sum_{k \in [K]}a_k^{p-1} \leq (\sum_{k \in [K]}a_k)^{p-1} \le O( A_K^{p-1})$ and $A_k = O(A_{\max})$. This finishes the proof. 

\end{proof}

\section{$\ell_p$ Parallel Stochastic Convex Optimization: Upper Bounds}
\label{section: parallel_upper_bounds}
\subsection{Preliminaries}
We study the following parallel SCO problem, where access to the objective is only allowed through the stochastic gradient oracle. 
\begin{problem}
\label{problem:parallel_problem}
    Let $f: \mathbb{R}^d \rightarrow \mathbb{R}$ be convex. Assume that there exists a stochastic gradient oracle $g$ such that for all $x \in \mathbb{R}^d$, $\mathbb{E}[g(x)] \in \partial f(x)$ and $\mathbb{E}||g(x)||_q^q \le L^q$. Our goal is to produce $x \in \mathbb{R}^d$ such that \begin{flalign*}
        &\mathbb{E}[f(x)] \le \min_{x^* \in \mathbb{B}(R)}f(x^*) + \epsilon_{\mathrm{opt}}
    \end{flalign*} for some $\epsilon_{\mathrm{opt}} > 0$. 
\end{problem} We will use the query depth, computational depth, and the total number of stochastic gradient oracle queries as metrics for characterizing the efficiency of our algorithms. Following \cite{carmon2023resqueing}, we present a few useful oracle definitions which will be useful in our analysis. 

\begin{definition}[Line search oracle]
We say $\mathcal{O}_{\rm ls}$ is a $(\Delta,\lambda)$-line search oracle for
$F:\mathbb{R}^d \to \mathbb{R}$ if given $\bar{x}\in\mathbb{R}^d$,
$\mathcal{O}_{\rm ls}$ returns $x\in\mathbb{R}^d$ with
\[\left\|x -x^\star_{\bar{x},\lambda}
\right\|_p \le \Delta .\]
\end{definition}

\begin{definition}[Ball optimization oracle]
We say $\mathcal{O}_{\rm bo}$ is a $(\phi,\lambda)$-ball optimization oracle for
$F:\mathbb{R}^d \to \mathbb{R}$ if given $\bar{x}\in\mathbb{R}^d$,
$\mathcal{O}_{\rm bo}$ returns $x\in\mathbb{R}^d$ with
\[\mathbb{E}\left[F(x) + \frac{\lambda}{p}\|x-\bar{x}\|^p_p
\right]\leq F\left(x^\star_{\bar{x},\lambda}\right)+\frac{\lambda}{p}
\left\|x^\star_{\bar{x},\lambda}-\bar{x}\right\|^p_p+\phi .\]
\end{definition}

\begin{definition}[Stochastic proximal oracle]
We say $\mathcal{O}_{\rm sp}$ is a $(\Delta,\sigma,\lambda)$-stochastic
proximal oracle for $F:\mathbb{R}^d\to\mathbb{R}$ if, given
$\bar{x}\in\mathbb{R}^d$, $\mathcal{O}_{\rm sp}$ returns a random vector
$g\in\mathbb{R}^d$ such that
\[\left\|\mathbb{E}g-g^\star_{\bar{x},\lambda}\right\|_q
\leq \Delta,\qquad\mathbb{E}
\left\|g-\mathbb{E}g\right\|_q^q
\leq \sigma^q,\]
\end{definition}

\begin{proposition}
\label{proposition:init_prop_guarantee}
Let $F:\mathbb{R}^d \to \mathbb{R}$ be $L$-Lipschitz and convex, and let
$x^\star \in \mathbb{B}_p(R)$. There is an algorithm \textsc{BallAccel} (Algorithm \ref{alg:stochastic_accelerated_prox_point}) taking parameters $r \in [0,R]$ and $\epsilon_{\rm opt} \in (0,LR]$ with the following guarantee. Define
\[\kappa := \frac{LR}{\epsilon_{\rm opt}},\qquad K := \left(\frac{R}{r}\right)^{\frac{p}{p+1}},
\qquad \lambda^\star :=
\frac{\epsilon_{\rm opt}K^p}{R^p}\log^p\kappa.\]
For a universal constant $C_{\rm ba}>0$, \textsc{BallAccel} runs in at most
$C_{\rm ba}K\log\kappa$ iterations and produces $x\in\mathbb{R}^d$ such that
\[
\mathbb{E}F(x) \le F(x^\star) + \epsilon_{\rm opt}.
\]
Moreover, in each iteration \textsc{BallAccel} requires the following oracle calls
(all for $F$).

\begin{enumerate}
    \item At most $C_{\rm ba}\log\left(\frac{R\kappa}{r}\right)$ calls to a
    $\left(\frac{r}{C_{\rm ba}},\lambda\right)$-line search oracle with values of
    $\lambda \in \left[\frac{\lambda^\star}{C_{\rm ba}},
    \frac{C_{\rm ba}L}{r^{p-1}}\right]$.

    \item A single call to a
    $\left(\frac{\lambda r^p}{C_{\rm ba}\log^{p+1}\kappa},\lambda\right)$-ball
    optimization oracle with
    $\lambda \in \left[\frac{\lambda^\star}{C_{\rm ba}},
    \frac{C_{\rm ba}L}{r^{p-1}}\right]$.

    \item A single call to a
    $\left(\frac{\epsilon_{\rm opt}}{C_{\rm ba}R},
    \frac{\epsilon_{\rm opt}K^{1/p}}{C_{\rm ba}R},\lambda\right)$-stochastic
    proximal oracle with
    $\lambda \in \left[\frac{\lambda^\star}{C_{\rm ba}},
    \frac{C_{\rm ba}L}{r^{p-1}}\right]$.
\end{enumerate}
\end{proposition}

\begin{proof} \,
    \begin{itemize}
    \item \textbf{Iteration Count}:
    Having $\lambda_{\min} = O(\lambda^\star)$ and $\epsilon = \epsilon_{\rm opt}$, our analogous of \textsc{BallAccel} runs in at most $C_{\rm ba}K\log\kappa$ iterations, where the value of $K_{\max}$ from Proposition \ref{proposition:stochastic_prox_method_guarantee} is $O(K \log \kappa)$. 
    \item \textbf{Line search oracle parameters}:
    From Theorem \ref{theorem:nextlambda_guarantee}, the line search oracle has complexity bounded by \begin{flalign*}
    O\left(\log\frac{\lambda_{\max}}{\lambda_{\min}}+\log\left(1+\frac{
    R+\frac{1}{p-1}\left(\frac{G}{\lambda_{\min}}\right)^{1/(p-1)}}{
    r}\right)\right)
    \end{flalign*}
    With $\lambda_{\max} = \frac{C_{\rm ba}L}{r^{p-1}}$ and $\lambda_{\min} = \frac{\lambda^\star}{C_{\rm ba}}$, we have, as $K \ge 1$,
    \[\frac{\lambda_{\max}}{\lambda_{\min}} = \frac{C_{\rm ba}L}{r^{p-1}} \times \frac{C_{\rm ba}}{\lambda^\star} = \frac{C_{\rm ba}^2 LR^p}{r^{p-1} \epsilon K^p\log^p \kappa}  = \frac{C_{\rm ba}^2 \kappa R^{p-1}}{r^{p-1}  K^p\log^p \kappa} = \frac{C_{\rm ba}^2 \kappa }{ K^{1/p}\log^p \kappa} \le \frac{C_{\rm ba}^2 \kappa }{\log^p \kappa}\]
    Now, note that
    \begin{align*}
    \frac{1}{r}\left(\frac{G}{\lambda_{\min}}\right)^{\frac{1}{p-1}}&=
    \frac{R}{r}\frac{\kappa^{1/(p-1)}}{K^{p/(p-1)}
    \log^{p/(p-1)}\kappa}=\frac{\kappa^{1/(p-1)}}{
    (R/r)^{1/(p^2-1)}\log^{p/(p-1)}\kappa}\le
    \kappa^{1/(p-1)}\le \kappa,\end{align*}
    Therefore, \begin{align*}\log\left(1+\frac{
    R+\frac{1}{p-1}\left(\frac{G}{\lambda_{\min}}\right)^{1/(p-1)}}{r}\right)\le
    \log\left(1+\frac{R}{r}+\kappa\right)=
    O\left(\log\left(\frac{R\kappa}{r}\right)\right),
    \end{align*}
    Together,
    \[O\left(\log\frac{\lambda_{\max}}{\lambda_{\min}}+\log\left(1+\frac{
    R+\frac{1}{p-1}\left(\frac{G}{\lambda_{\min}}\right)^{1/(p-1)}}{
    r}\right)\right) \leq O\bigg(\log\frac{R\kappa}{r}\bigg)\]
    \item \textbf{Ball optimization oracle parameters}: By Lemma \ref{lemma:biasedreduced_potentialdecreasebound}, we have that \begin{flalign*}
    &\mathbb{E}[A_{k + 1} (f(x_{k + 1}) - f(x^*)) + 2D_\Phi(x^*, v_{k + 1})] \\ 
        &\le A_k(f(x_k) - f(x^*)) + 2D_\Phi(x^*, v_k) + A_{k + 1} \varphi_{k + 1} - \frac{1 - 2^{1 - q}}{q} \lambda_{k + 1} A_{k + 1} ||\hat{x}_{k + 1} - y_k||_p^p \\ 
        &+ \frac{2^{2 - q} a_{k + 1}^q}{q} \sigma_{k + 1}^q + 2 a_{k + 1} \delta_{k + 1} R 
         \\&\le A_k(f(x_k) - f(x^*))
        + 2D_\Phi(x^*, v_k)+ A_{k + 1} \varphi_{k + 1} - \frac{1}{q}\left(1-2^{1 - q}\right)
        \lambda_{k + 1} A_{k + 1}\|\hat{x}_{k + 1} - y_k\|_p^p \\
        &\quad + \frac{2^{2 - q} a_{k + 1}^q}{q}\sigma_{k + 1}^q+ 2 a_{k + 1} \delta_{k + 1} R .
\end{flalign*}
    If our algorithm does not terminate, one of the following occurs:
    \begin{itemize}
        \item [(1)] $||\hat{x}_{k + 1} - y_k||_p^p = \Omega(r^p)$
        \item  [(2)] $\lambda_{k+1} = O(\lambda^\star)$
    \end{itemize}
    If $(1)$ occurs, then $- \frac{1 - 2^{1-q}}{q} \lambda_{k + 1} A_{k + 1} ||\hat{x}_{k + 1} - y_k||_p^p = -\Omega(\lambda_{k+1}A_{k+1}r^p)$. If we choose $\varphi_{k+1}  \leq \frac{\lambda_{k+1}r^p}{C_{\rm{ba}}\log^{p+1}\kappa}$, then $A_{k+1}\varphi_{k+1}$ can be bounded by the first case. If $(2)$ occurs, since Proposition \ref{proposition:stochastic_prox_method_guarantee} requires $A_{\max} \leq O(R^p/\epsilon)$, with the fact that the total iteration is $K \log \kappa$ and $\varphi_{k+1}  \leq \frac{\lambda_{k+1}r^p}{C_{\rm{ba}}\log^{p+1}\kappa}$ by design, we have
    \[\sum_{k+1 \leq K}A_{k+1}\varphi_{k+1} \leq O\Bigg(K\log\kappa \cdot \frac{R^p}{\epsilon_{\rm{opt}}}\cdot \frac{\lambda^\star r^p}{C_{\rm{ba}}\log^{p+1}\kappa}\Bigg) \leq O(R^p)\]
    Here, the first term is the iteration count, the second term is due to an upper bound on $A_{k+1}$, and
    the third term is bounded since $\lambda_{k+1} = O(\lambda^\star)$.
    \item \textbf{Stochastic proximal oracle parameters}: In Proposition \ref{proposition:stochastic_prox_method_guarantee} we require $\sigma^q = O(\epsilon_{\rm opt}/a_{k+1}^{q-1})$, using the relation that $A_{k+1}^{p-1} = 2^{p - 1}\lambda_{k+1}a_{k+1}^p$, we have
    \[\sigma^q = \frac{\epsilon_{\rm opt}}{a_{k+1}^{q-1}} = \frac{\epsilon_{\rm opt}\lambda_{k+1}^{1/p(p-1)}2^{\frac{1}{p}}}{A_{k+1}^{1/p}} = \Omega\Bigg(\frac{\epsilon_{\rm opt}^{1 + 1/p} \cdot (\lambda^\star)^{1/p(p-1)}}{R}\Bigg) = \Omega\Bigg(\frac{\epsilon_{\rm opt}^{p/(p-1)} \cdot K^{1/(p-1)}\log^{1/(p-1)}\kappa}{R^{p/(p-1)}}\Bigg)\]
    where the second equality uses $A_{k+1} \le  O(R^p/\varepsilon_{\mathrm{opt}})$. Hence, $\sigma = \tilde{\Omega}(\epsilon_{\rm opt}K^{1/p}/R)$.
    \end{itemize}
\end{proof}

\begin{lemma}
\label{lem:bo-amplification}
Let $\bar{x}\in\mathcal X$, $\lambda>0$. Suppose a $\left(
\frac{2^{2-p}\lambda r^p}{p\,3^{p+1}C_{\rm ba}^p},
\lambda\right)$-ball optimization oracle returns $x\in\mathcal X\cap B_p(\bar x,r)$,
then \begin{equation}
\label{eq:bo-single-success}
\Pr\left(\|x-x^\star_{\bar{x},\lambda}\|_p
\le\frac{r}{3C_{\rm ba}}\right)\ge\frac{2}{3}.
\end{equation}
Moreover, let $x_1,\ldots,x_m$ be $m$ independent outputs of the ball optimization oracle. If
\[m = \lceil 20\log(1/\zeta)\rceil,\]
then, with probability at least $1-\zeta$, at least $0.51m$ of the points $x_1,\ldots,x_m$ satisfy
\[\|x_i-x^\star_{\bar{x},\lambda}\|_p
\le\frac{r}{3C_{\rm ba}}.\]
Consequently, applying $\mathrm{Aggregate}$ (Algorithm \ref{alg:aggregate}) to
$x_1,\ldots,x_m$ returns a point $\widehat{x}$ satisfying
\[\|\widehat{x}-x^\star_{\bar{x},\lambda}\|_p\le\frac{r}{C_{\rm ba}}\]
with probability at least $1-\zeta$.

\end{lemma}

\begin{proof}
By uniform convexity, the function $\Phi(x):=\frac{1}{p}\|x\|_p^p$ satisfies, for all $u,v\in\mathbb R^d$, \begin{equation} \label{eq:lp-uniform-convexity} 
\Phi(u) \ge \Phi(v) + \langle J_p(v),u-v\rangle + \frac{2^{2-p}}{p}\|u-v\|_p^p, \end{equation} 
where $ J_p(v) = \left( |v_i|^{p-2}v_i \right)_{i=1}^d$. Therefore, \begin{align*} 
\frac{\lambda}{p}\|x-\bar{x}\|_p^p &\ge \frac{\lambda}{p}\|y-\bar{x}\|_p^p + \lambda \left\langle J_p(y-\bar{x}),x-y \right\rangle + \frac{2^{2-p}\lambda}{p}\|x-y\|_p^p. \end{align*} Since $F$ is convex, for any $g\in\partial F(y)$, $ F(x) \ge F(y)+\langle g,x-y\rangle$.  Combining two inequalities gives \begin{align} 
F(x) + \frac{\lambda}{p}\|x-\bar{x}\|_p^p \ge F(y) + \frac{\lambda}{p}\|y-\bar{x}\|_p^p + \left\langle g+\lambda J_p(y-\bar{x}), x-y \right\rangle + \frac{2^{2-p}\lambda}{p} \|x-y\|_p^p. \label{eq:prox-uniform-convexity} \end{align} Taking $y=x^\star_{\bar{x},\lambda}$ and using the first-order optimality condition for $x^\star_{\bar{x},\lambda}$ gives \begin{equation} \label{eq:prox-gap-distance} F(x) + \frac{\lambda}{p}\|x-\bar{x}\|_p^p \ge F(y) + \frac{\lambda}{p}\|y-\bar{x}\|_p^p + \frac{2^{2-p}\lambda}{p} \|x-y\|_p^p.
\end{equation} 
Combining \eqref{eq:prox-gap-distance} with the definition of $\left(\frac{\lambda r^p}{C_{\rm ba}},\lambda\right)$-ball optimization oracle gives 
\begin{align*} \E\|x-x^\star_{\bar{x},\lambda}\|_p^p &\le \frac{p\,2^{p-2}}{\lambda} \E\left[ F(x) + \frac{\lambda}{p}\|x-\bar{x}\|_p^p - F(x^\star_{\bar{x},\lambda}) - \frac{\lambda}{p}
\|x^\star_{\bar{x},\lambda}-\bar{x}\|_p^p  \right] \\ &\le \frac{p\,2^{p-2}}{\lambda} \cdot \frac{2^{2-p}\lambda r^p} {p\,3^{p+1}C_{\rm ba}^p} = \frac{r^p}{3^{p+1}C_{\rm ba}^p}. \end{align*} Therefore, by Markov's inequality, \begin{align*} \Pr\left( \|x-x^\star_{\bar{x},\lambda}\|_p > \frac{r}{3C_{\rm ba}} \right) &= \Pr\left( \|x-x^\star_{\bar{x},\lambda}\|_p^p > \frac{r^p}{3^pC_{\rm ba}^p} \right) \le \frac{ \E\|x-x^\star_{\bar{x},
\lambda}\|_p^p}{r^p/(3^pC_{\rm ba}^p)} \le \frac{1}{3}. \end{align*} 
Hence \[\Pr\left(\|x-x^\star_{\bar{x},\lambda}\|_p \le \frac{r}{3C_{\rm ba}} \right) \ge \frac{2}{3}. \] For $i=1,\ldots,m$, define \[ Z_i := \mathbf{1}\left\{ \|x_i-x^\star_{\bar{x},\lambda}\|_p \le \frac{r}{3C_{\rm ba}} \right\}. \] Then $Z_1,\ldots,Z_m$ are independent Bernoulli random variables with $\E Z_i\ge\frac{2}{3}$. By a Chernoff bound, with probability at least $1-\zeta$, at least $0.51m$ of the points $x_1,\ldots,x_m$ satisfy
\[\|x_i-x^\star_{\bar{x},\lambda}\|_p
\le\frac{r}{3C_{\rm ba}}.\]
Therefore, the preconditions of Fact \ref{fact:aggregateguarantee} are satisfied, and
$\mathrm{Aggregate}$ (Algorithm \ref{alg:aggregate}) returns $\widehat{x}$ such that \[\|\widehat{x}-x^\star_{\bar{x},\lambda}\|_p\le\frac{r}{C_{\rm ba}}.\] 
\end{proof}

\begin{algorithm}[t]
\caption{$p$-MorGradEst}
\label{alg:p-morgrad}
\begin{algorithmic}[1]
\State \textbf{Input:} $r$-BROO $\mathcal O_{\lambda,\delta}(\cdot)$ of radius $r$,
query point $y$, regularization $\lambda>0$,
bias $\beta>0$, $\sigma>0$

\State Set
$T_{\max}=\max\left\{2,\frac{C_p\lambda^p r^{p(p-1)}}{\beta^p}\right\},
T_0=\max\left\{1,\frac{C_p\lambda^p r^{p(p-1)}\log^{p-1}(2T_{\max})}{\sigma^p}\right\}$
\State $j_{\max}=\lfloor\log_2T_{\max}\rfloor,\quad\delta_j=\frac{r}{(2^jT_0)^{1/p}}$
\State $x_0=\mathcal O_{\lambda,\delta_0}(y)$, $g_0=\lambda J_p(y-x_0)$
\State Sample $J\sim\mathrm{Geom}(1/2)$
\If{$J\le j_{\max}$}
    \State $x_J=\mathcal O_{\lambda,\delta_J}(y)$, $x_{J-1}=\mathcal O_{\lambda,\delta_{J-1}}(y)$
    \State $g_J=\lambda J_p(y-x_J)$, $g_{J-1}=\lambda J_p(y-x_{J-1})$
    \State $\widehat g=g_0+2^J(g_J-g_{J-1})$
\Else
    \State $\widehat g=g_0$
\EndIf
\State \Return $\widehat g$
\end{algorithmic}
\end{algorithm}

\begin{lemma}
\label{lem:p-morgrad-ball-feasible}
Let $F:\mathcal X\to\mathbb R$ satisfy Assumption~\ref{Assumption}.
For a query point $y\in\mathcal X$ and $\lambda,r>0$, define
\[
x^\star\in\arg\min_{x\in\mathcal X}
\left\{
F(x)+\frac{\lambda}{p}\|x-y\|_p^p
\right\},
\qquad
g=\lambda J_p(y-x^\star),
\]
where $J_p(z)=(|z_i|^{p-2}z_i)_{i=1}^d$.
Suppose that $\|x^\star-y\|_p\le r$.
Let $\widehat g$ be the geometric estimator in
Algorithm~\ref{alg:p-morgrad}, run with these values of
$T_0$ and $j_{\max}$ and with
\[
\delta_j=\frac{r}{(2^jT_0)^{1/p}},
\qquad 0\le j\le j_{\max},
\]
Then, for a constant $C_p>0$ depending only on $p$,
\begin{align}
\|\mathbb E\widehat g-g\|_q
&\le
\frac{C_p\lambda r^{p-1}}{\sqrt{2^{j_{\max}}T_0}},
\label{eq:morgrad-ball-bias}\\
\mathbb E\|\widehat g-\mathbb E\widehat g\|_q^q
&\le\frac{C_p\lambda^qr^p}{T_0^{q/2}}\left(
1+\sum_{j=1}^{j_{\max}}2^{-j(1-q/2)}\right).
\label{eq:morgrad-ball-q-moment}
\end{align}
\end{lemma}
\begin{proof}
Let $\Phi(u)=\|u\|_p^p/p$, so that
\[D_\Phi(u,v)=\Phi(u)-\Phi(v)-\langle J_p(v),u-v\rangle.\]
We first show that
\begin{equation}
\label{eq:morgrad-duality-bregman}
\|J_p(u)-J_p(v)\|_q^2
\le
C_p(\|u\|_p+\|v\|_p)^{p-2}D_\Phi(u,v).
\end{equation}
For $p=2$, this follows from
$\|u-v\|_2^2=2D_\Phi(u,v)$.
For $p>2$, integrating the second derivative of
$|a|^p/p$ gives
\[D_\Phi(u,v)=(p-1)
\sum_{i=1}^d(u_i-v_i)^2
\int_0^1(1-t)|v_i+t(u_i-v_i)|^{p-2}\,dt\ge
c_p\sum_{i=1}^d
(|u_i|+|v_i|)^{p-2}(u_i-v_i)^2\]
for a constant $c_p>0$ depending only on $p$.
Indeed, after rescaling so that $|a|+|b|=1$,
the integral
$\int_0^1(1-t)|b+t(a-b)|^{p-2}\,dt$
is continuous and strictly positive on this compact set,
and therefore has a positive minimum.
Also,
\[\bigl||u_i|^{p-2}u_i-|v_i|^{p-2}v_i\bigr|\le
(p-1)(|u_i|+|v_i|)^{p-2}|u_i-v_i|.\]
H\"older's inequality therefore yields
\begin{align*}
\|J_p(u)-J_p(v)\|_q^2
&\le
(p-1)^2
\left(\sum_{i=1}^d(|u_i|+|v_i|)^p\right)^{(p-2)/p}
\sum_{i=1}^d
(|u_i|+|v_i|)^{p-2}(u_i-v_i)^2\\
&\le
C_p(\|u\|_p+\|v\|_p)^{p-2}D_\Phi(u,v),
\end{align*}
which proves \eqref{eq:morgrad-duality-bregman}.
By the convexity of $F$ and optimality of $x^\star$,
for every $x\in\mathcal X$,
\[
F(x)+\frac{\lambda}{p}\|x-y\|_p^p
-F(x^\star)-\frac{\lambda}{p}\|x^\star-y\|_p^p
\ge
\lambda D_\Phi(x-y,x^\star-y).
\]
Since $\|x_j-y\|_p\le r$ and
$\|x^\star-y\|_p\le r$, applying
\eqref{eq:morgrad-duality-bregman} with
$u=x_j-y$ and $v=x^\star-y$ gives
\begin{align}
\mathbb E\|g_j-g\|_q^2
&\le
C_p\lambda r^{p-2}
\mathbb E\left[
F(x_j)+\frac{\lambda}{p}\|x_j-y\|_p^p
-F(x^\star)-\frac{\lambda}{p}\|x^\star-y\|_p^p
\right]\notag\\
&\le
\frac{C_p\lambda^2r^{2p-2}}{2^jT_0}.
\label{eq:morgrad-ball-level-second-moment}
\end{align}
Here we used $g_j=\lambda J_p(y-x_j)$,
$g=\lambda J_p(y-x^\star)$, and
$J_p(-u)=-J_p(u)$.
By the definition of $\widehat g$,
\[
\mathbb E\widehat g
=
\mathbb E g_0+
\sum_{j=1}^{j_{\max}}
\left(\mathbb E g_j-\mathbb E g_{j-1}\right)
=
\mathbb E g_{j_{\max}}.
\]
Consequently,
\[
\|\mathbb E\widehat g-g\|_q
\le
\left(\mathbb E\|g_{j_{\max}}-g\|_q^2\right)^{1/2}
\le
\frac{C_p\lambda r^{p-1}}
{\sqrt{2^{j_{\max}}T_0}},
\]
which proves \eqref{eq:morgrad-ball-bias}.
Since $1<q\le2$, Jensen's inequality and
$(p-1)q=p$ give
\[
\mathbb E\|g_j-g\|_q^q
\le
\left(\mathbb E\|g_j-g\|_q^2\right)^{q/2}
\le
\frac{C_p\lambda^qr^p}{(2^jT_0)^{q/2}}.
\]
It follows that
\[
\mathbb E\|g_j-g_{j-1}\|_q^q
\le
2^{q-1}\left(
\mathbb E\|g_j-g\|_q^q+
\mathbb E\|g_{j-1}-g\|_q^q
\right)
\le
\frac{C_p\lambda^qr^p}{(2^jT_0)^{q/2}}.
\]
Using the definition of $\widehat g$ once more,
\begin{align*}
\mathbb E\|\widehat g-\mathbb E\widehat g\|_q^q
&\le
2^q\mathbb E\|\widehat g-g\|_q^q\\
&\le
C_p\left(
\mathbb E\|g_0-g\|_q^q+
\sum_{j=1}^{j_{\max}}
2^{j(q-1)}\mathbb E\|g_j-g_{j-1}\|_q^q
\right)\\
&\le
\frac{C_p\lambda^qr^p}{T_0^{q/2}}
\left(
1+\sum_{j=1}^{j_{\max}}2^{-j(1-q/2)}
\right).
\end{align*}
This proves \eqref{eq:morgrad-ball-q-moment}.
The argument also applies after conditioning on the history
before the call.
\end{proof}
\begin{lemma}
\label{lem:p-lemma8}
Let $F$ satisfy Assumption~\ref{Assumption} with a minimizer $x^\star$.
Let $\ell=\log\left(\frac{GR^2}{r\epsilon_{\rm opt}}\right)$.
Let
\[
\varphi_k=\frac{\lambda_k r^p}{C_p\ell^{p+1}}
\quad,\quad
\sigma_k^q=\frac{\lambda_k^q r^p}{C_p\ell^{p+1}}
\quad,\quad
\beta_k=\frac{\epsilon_{\rm opt}}{C_pR}
\quad,\quad
A_0=\frac{R^{p-1}}{G}
\text{ and }A_{\max}=\frac{C_pR^p}{\epsilon_{\rm opt}}.
\]
Define
$\widehat{x}_k:=\argmin_{x\in\mathcal X}
\left\{F(x)+\frac{\lambda_k}{p}
\|x-y_{k-1}\|_p^p\right\}$,
and assume that for each $k$ we have
$\|\widehat{x}_k-y_{k-1}\|_p\le r$
and that one of the following must occur:
\begin{enumerate}
    \item
    $\lambda_k<2\lambda_m:=
    \frac{2\epsilon_{\rm opt}}
    {r^{\frac{p^2}{p+1}}R^{\frac{p}{p+1}}}\ell^p$, or
    \item
    $\|\widehat{x}_k-y_{k-1}\|_p
    \ge\frac34r$.
\end{enumerate}
Suppose $\lambda_k\ge\lambda_m$ for every $k$ and
\[A_{k+1}^{p-1}=2^{p-1}\lambda_{k+1}a_{k+1}^p,\qquad
A_{k+1}=A_k+a_{k+1}.\]
Let $K\le K_{\max}$ be the stopping time of the algorithm,
which stops when $A_k\ge A_{\max}$ or $k=K_{\max}$.
Then, for sufficiently large $C_p$, there exist constants
$c_{1,p},c_{2,p}>0$ depending only on $p$ such that
\[
\E\left[
A_K\left(F(x_K)-F(x^\star)
-c_{1,p}\epsilon_{\rm opt}\right)
+c_{2,p}\sum_{i\le K}\lambda_iA_i
\|\widehat{x}_i-y_{i-1}\|_p^p
\right]
\le A_0\left(F(x_0)-F(x^\star)\right)
+C_pR^p.
\]
\end{lemma}

\begin{proof}
Define the filtration
\[
\mathcal F_k=\sigma(x_1,v_1,A_1,\zeta_1,\ldots,x_k,v_k,A_k,\zeta_k),
\]
where $\zeta_i$ is the internal randomness of \textsc{NextLambda}, and note that $A_{k+1}$, $y_k$, and $\widehat{x}_{k+1}$ are deterministic when conditioned on $x_k,v_k,A_k,\zeta_k$. Following the proof of \Cref{lemma:biasedreducedK_maxbound}, we define
\[
M_k=A_k\left(F(x_k)-F(x^\star)-\frac{\epsilon_{\rm opt}}{20}\right)
+2D_\Phi(x^\star,v_k)
+\frac{1-2^{1-q}}{2q}\sum_{i\le k}\lambda_iA_i
\|\widehat{x}_i-y_{i-1}\|_p^p
\]
and show it is a supermartingale adapted to the filtration $\mathcal F_k$. From Lemma~\ref{lemma:biasedreduced_potentialdecreasebound}, we have
\begin{align}
\E[M_{k+1}\mid\mathcal F_k]
&\le A_k(F(x_k)-F(x^\star)) + 2D_\Phi(x^\star,v_k)
-\frac{1-2^{1-q}}q\lambda_{k+1}A_{k+1}
\|\widehat{x}_{k+1}-y_k\|_p^p\notag\\
&\quad+\mu_{k+1}-A_{k+1}\frac{\epsilon_{\rm opt}}{20}
+\frac{1-2^{1-q}}{2q}\sum_{i\le k+1}\lambda_iA_i
\|\widehat{x}_i-y_{i-1}\|_p^p,
\end{align}
where
\[
\mu_{k+1}:=A_{k+1}\varphi_{k+1}
+\frac{2^{2-q}a_{k+1}^q}{q}\sigma_{k+1}^q
+2Ra_{k+1}\beta_{k+1}.
\]
Substituting the values of $\varphi_{k+1}$, $\sigma_{k+1}^q$, and $\beta_{k+1}$ into the definition of $\mu_{k+1}$ gives
\[
\mu_{k+1}
=\frac{\lambda_{k+1}A_{k+1}r^p}{C_p\ell^{p+1}}
+\frac{2^{2-q}a_{k+1}^q\lambda_{k+1}^qr^p}{qC_p\ell^{p+1}}
+\frac{2a_{k+1}\epsilon_{\rm opt}}{C_p}.
\]
Recall that $A_{k+1}^{p-1}=2^{p-1}\lambda_{k+1}a_{k+1}^p$, therefore $a_{k+1}^q\lambda_{k+1}^q=\lambda_{k+1}A_{k+1}/2$ and
\[
\mu_{k+1}
=\frac{1+2^{1-q}/q}{C_p\ell^{p+1}}
\lambda_{k+1}A_{k+1}r^p
+\frac{2a_{k+1}\epsilon_{\rm opt}}{C_p}
=\frac{2^{(p-1)/p}(1+2^{1-q}/q)
 a_{k+1}\lambda_{k+1}^{1+1/p}A_{k+1}^{1/p}r^p}
{C_p\ell^{p+1}}
+\frac{2a_{k+1}\epsilon_{\rm opt}}{C_p}.
\]
If $\|\widehat{x}_{k+1}-y_k\|_p\ge3r/4$, then for sufficiently large $C_p$ and $\ell\ge1$,
\[
\mu_{k+1}\le
\frac{1-2^{1-q}}{2q}\lambda_{k+1}A_{k+1}
\|\widehat{x}_{k+1}-y_k\|_p^p
+\frac{a_{k+1}\epsilon_{\rm opt}}{20}.
\]
Else, if $\lambda_{k+1}<2\lambda_m$, we have
\[
\mu_{k+1}
<\frac{2^{(p-1)/p}(1+2^{1-q}/q)
 a_{k+1}(2\lambda_m)^{1+1/p}A_{k+1}^{1/p}r^p}
{C_p\ell^{p+1}}
+\frac{2a_{k+1}\epsilon_{\rm opt}}{C_p}.
\]
Now, note that $\lambda_mA_{\max}=C_p(R/r)^{p^2/(p+1)}\ell^p\ge1$ for sufficiently large $C_p$ when $0<r\le R$ and $\ell\ge1$, so the \textsc{NextLambda} lower bound $\lambda_k\ge\lambda_m$ implies $\lambda_k\ge1/A_{\max}$. From the definition of the stopping time, $A_{K-1}\le A_{\max}$. If $a_K>A_{\max}$, then $A_K=A_{K-1}+a_K\le2a_K$, and the weight equation gives
\[
a_K^p=\frac{A_K^{p-1}}{2^{p-1}\lambda_K}
\le\frac{(2a_K)^{p-1}A_{\max}}{2^{p-1}}
=a_K^{p-1}A_{\max},
\]
a contradiction. Therefore $a_K\le A_{\max}$, and $A_k\le A_K\le2A_{\max}$ for every $k\le K$. Thus, when $\lambda_{k+1}<2\lambda_m$, the bound on $\mu_{k+1}$ becomes
\begin{align*}
\mu_{k+1}
&<\frac{2^{(p-1)/p}(1+2^{1-q}/q)
 a_{k+1}(2\lambda_m)^{1+1/p}(2A_{\max})^{1/p}r^p}
{C_p\ell^{p+1}}
+\frac{2a_{k+1}\epsilon_{\rm opt}}{C_p}\\
&=a_{k+1}\epsilon_{\rm opt}
\left(\frac{4\cdot2^{1/p}(1+2^{1-q}/q)}{C_p^{(p-1)/p}}
+\frac2{C_p}\right)\\
&\le a_{k+1}\frac{\epsilon_{\rm opt}}{20}
\le\frac{1-2^{1-q}}{2q}\lambda_{k+1}A_{k+1}
\|\widehat{x}_{k+1}-y_k\|_p^p
+a_{k+1}\frac{\epsilon_{\rm opt}}{20},
\end{align*}
where the penultimate inequality holds for sufficiently large $C_p$.
Therefore, for every $k<K$, we have
\[
\mu_{k+1}\le\frac{1-2^{1-q}}{2q}\lambda_{k+1}A_{k+1}
\|\widehat{x}_{k+1}-y_k\|_p^p
+a_{k+1}\frac{\epsilon_{\rm opt}}{20}.
\]
Noting that $\E|M_k|<\infty$ and substituting the bound on $\mu_{k+1}$ into (21), we get
\begin{align*}
\E[M_{k+1}\mid\mathcal F_k]
&\le A_k\left(F(x_k)-F(x^\star)-\frac{\epsilon_{\rm opt}}{20}\right)
+2D_\Phi(x^\star,v_k)
+\frac{1-2^{1-q}}{2q}\sum_{i\le k}\lambda_iA_i
\|\widehat{x}_i-y_{i-1}\|_p^p
=M_k.
\end{align*}
Therefore $M_k$ is a supermartingale adapted to $\mathcal F_k$. Since $K\le K_{\max}$ is a bounded stopping time adapted to $\mathcal F_k$, by the optional stopping theorem for supermartingales we have
\begin{align*}
\E M_K&\le M_0
=A_0\left(F(x_0)-F(x^\star)-\frac{\epsilon_{\rm opt}}{20}\right)
+2D_\Phi(x^\star,v_0)\\
&\le A_0(F(x_0)-F(x^\star))+C_pR^p,
\end{align*}
where $v_0=x_0$, and $2D_\Phi(x^\star,x_0)\le C_pR^p$ holds when the mirror potential is centered at $x_0$ (or equivalently after translating $x_0$ to the origin). This completes the proof.
\end{proof}

\begin{proposition}
\label{proposition:prop2_equivalent}
Let $F:\mathbb R^d\to\mathbb R$ be convex and
$L$-Lipschitz with respect to $\|\cdot\|_p$, and let
$x^\star\in\mathbb B_p(R)$.
There is an implementation of \textsc{BallAccel}
(see Proposition~\ref{proposition:init_prop_guarantee})
taking parameters $r\in(0,R]$ and
$\epsilon_{\rm opt}\in(0,LR/e]$ with the following guarantee.
Define $\kappa$, $K$, and $\lambda^\star$ as in
Proposition~\ref{proposition:init_prop_guarantee}. For a sufficiently large constant $C_{\rm ba}>0$
depending only on $p$, \textsc{BallAccel} runs in at most $C_{\rm ba}K\log\kappa
\log\left(\frac{R\kappa}{r}\right)$ iterations and produces $x\in\mathbb R^d$ such that $\mathbb EF(x)\le F(x^\star)+\epsilon_{\rm opt}$.
For general $F$, the ball optimization oracles below
are applied to
\[
x\mapsto F(x)+L(\|x\|_p-R)_+.
\]
Their calls satisfy:
\begin{enumerate}
    \item Each iteration makes at most $C_{\rm ba}
    \log^2\left(\frac{R\kappa}{r}\right)$ calls to a
    $\left(\frac{\lambda r^p}{C_{\rm ba}},\lambda\right)$
    -ball optimization oracle, with $\lambda\in\left[
    \frac{\lambda^\star}{C_{\rm ba}},\frac{C_{\rm ba}L}{r^{p-1}}\right]$. 
    \item For each
    $j\in[\lceil\log_2K+C_{\rm ba}\rceil]$, at most
    $C_{\rm ba}^2\,2^{-j}K\log\kappa
    \log\left(\frac{R\kappa}{r}\right)$
    iterations query a
    $\left(\frac{\lambda r^p}
    {C_p\,2^j\log^{2p}(R\kappa/r)},\lambda
    \right)$ -ball optimization oracle, for some
    $\lambda\in\left[
    \frac{\lambda^\star}{C_{\rm ba}},\frac{C_{\rm ba}L}{r^{p-1}}\right]$.
\end{enumerate}

Moreover, if $F$ admits a global minimizer in $\mathbb B_p(R)$,
the same guarantees and oracle-call bounds hold for an
implementation whose ball optimization oracles are applied
directly to $F$, without adding a radial penalty and without
projecting the final output.
\end{proposition}

\begin{proof}
We follow the organization of Appendix~C of
\cite{carmon2023resqueing}.
For the general implementation, it suffices to take
\[
x^\star\in\argmin_{x\in\mathbb B_p(R)}F(x).
\]
Apply the method to
$x\mapsto F(x)+L(\|x\|_p-R)_+$.
This function is convex and $2L$-Lipschitz, and
\[
F(x)+L(\|x\|_p-R)_+
\ge
F\left(
\frac{x}{\max\{1,\|x\|_p/R\}}
\right)
\ge F(x^\star).
\]
Thus its objective error is nonnegative everywhere,
and projecting the final output onto $\mathbb B_p(R)$
gives the required guarantee for $F$.
Replacing $L$ by $2L$ changes the bounds only by
constant factors.
For a stochastic-subgradient implementation, add a
subgradient of $L(\|x\|_p-R)_+$ to the original oracle.
Its dual norm is at most $L$, so Minkowski's inequality
bounds the resulting $q$-th moment by $(2L)^q$.
Initialize $x_0=v_0=0$, using Lipschitz bound $2L$
in the initialization of $A_0$, and retain
\[A_i=A_{i-1}+a_i,\qquad
A_i^{p-1}=2^{p-1}\lambda_i a_i^p.\] 
Use the mirror update over $\mathbb B_p(R)$ and define
\[
\widehat x_i
=
\operatorname*{argmin}_{x\in\mathbb R^d}
\left\{
F(x)+L(\|x\|_p-R)_+
+\frac{\lambda_i}{p}\|x-y_{i-1}\|_p^p
\right\}.
\]
Lemma~\ref{lemma:biasedreduced_potentialdecreasebound}
applies with the same $2Ra_i\beta_i$ term:
its radius bound is needed only for $x^\star$ and
the mirror iterates, while $\widehat x_i$ minimizes
the objective over $\mathbb R^d$.
Use oracle implementations that return points in their
queried radius-$r$ balls for the primal, line-search,
and MorGradEst calls.
Each execution stops at the first completed iteration
with $A_i\ge A_{\max}$, or after
$\lceil C_pK\log(R\kappa/r)\rceil$ iterations, where
\[
A_{\max}=\frac{C_{\rm ba}R^p}{\epsilon_{\rm opt}}.
\]
Here $C_{\rm ba}$ is sufficiently large and depends only on $p$.
Write $\tau$ for this stopping index.
An execution is accepted only if $A_\tau\ge A_{\max}$.
We first bound the calls required in each iteration.
The primal updates and the mirror constraint give
\[
\|x_i\|_p\le R+ir,
\qquad
\|x_i-v_i\|_p\le2R+ir.
\]
Each execution has $O(K\log(R\kappa/r))$ iterations.
Since $\|x_i-v_i\|_p\le2R+ir$ and $Kr\le R$, we have
\[
\|x_i-v_i\|_p
\le2R+O(Kr\log(R\kappa/r))
=O(R\log(R\kappa/r)).
\]
Using this bound in
Theorem~\ref{theorem:nextlambda_guarantee}
gives $O(\log(R\kappa/r))$ distance tests per iteration.
The mirror iterates remain in $\mathbb B_p(R)$.

For each test, use independent calls with expected
objective error $\lambda r^p/C_p$.
By Lemma~\ref{lem:bo-amplification},
$O(\log(R\kappa/r))$ repetitions followed by
\textsc{Aggregate} make the conditional failure probability
at most
\[
\frac{1}{
C_p\kappa K\lceil\log_2(128\kappa)\rceil
\log^2(R\kappa/r)}.
\]
A union bound over all tests and executions gives
total failure probability at most $1/(C_p\kappa)$.
Abort an execution before exceeding the stated test budget.
When all tests are accurate, this budget is sufficient.
This proves item~1 without function-value queries.

For the analysis, consider an execution without the caps
on oracle calls, in which inaccurate distance tests
are replaced by their exact outcomes.
Fresh primal and MorGradEst randomness retains its
conditional guarantees; no conditioning on future
test accuracy or future caps is used.
Theorem~\ref{theorem:nextlambda_guarantee} then gives
\[
\|\widehat x_i-y_{i-1}\|_p\le r,
\qquad
\lambda_i<2\lambda_{\min}
\quad\text{or}\quad
\|\widehat x_i-y_{i-1}\|_p\ge\frac{3r}{4},
\]
where $\lambda_{\min}$ differs from $\lambda^\star$
by a constant factor.
For the additional calls, use the construction of
Lemma~\ref{lem:p-morgrad-ball-feasible} with
\[
j_{\max}=\lceil\log_2K\rceil,
\qquad
T_0=C_{\rm ba}^{4p}\log^{2p}(R\kappa/r),
\qquad
\delta_j=\frac{r}{(2^jT_0)^{1/p}}.
\]
That lemma and an independent primal call give
\begin{align*}
\beta_i
&\le
\frac{C_p\lambda_i r^{p-1}}
{\sqrt{KC_{\rm ba}^{4p}}\,\log^p(R\kappa/r)},\\
\sigma_i^q
&\le
\frac{C_p\lambda_i^qr^p}
{C_{\rm ba}^{2pq}\log^{pq}(R\kappa/r)}
\left(
1+\sum_{j=1}^{\lceil\log_2K\rceil}
2^{-j(1-q/2)}
\right),\\
\varphi_i
&\le
\frac{C_p\lambda_i r^p}
{C_{\rm ba}^{4p}\log^{2p}(R\kappa/r)}.
\end{align*}
The sum in the second bound is bounded by a constant
depending only on $p$ when $p>2$, and is
$O(\log(R\kappa/r))$ when $p=2$. Using $a_i^q\lambda_i^q=\lambda_iA_i/2$ and
Young's inequality, we obtain, for any sufficiently
small fixed $c_p>0$ and sufficiently large $C_{\rm ba}$,
\[
A_i\varphi_i
+\frac{2^{2-q}a_i^q}{q}\sigma_i^q
+2Ra_i\beta_i
\le
\frac{c_p\lambda_iA_i r^p}
{\log^{p+1}(R\kappa/r)}
+
\frac{C_pR^p}
{C_{\rm ba}^{2p^2}K^{p/2}\log(R\kappa/r)}.
\]
If $\|\widehat x_i-y_{i-1}\|_p\ge3r/4$,
the first term is absorbed by a fixed fraction
of the negative displacement term.
Otherwise $\lambda_i<2\lambda_{\min}$.
The final-increment argument gives
$A_\tau\le2A_{\max}$, and
$K^{p+1}r^p=R^p$ implies
\[
\frac{\lambda_iA_i r^p}
{a_i\log^{p+1}(R\kappa/r)}
\le
C_pC_{\rm ba}^{1/p}\epsilon_{\rm opt}
\left(
\frac{\log\kappa}{\log(R\kappa/r)}
\right)^{p+1}
\le
C_pC_{\rm ba}^{1/p}\epsilon_{\rm opt}.
\]
Choose $c_p$ sufficiently small, using
$\sum_{i\le\tau}a_i\le2A_{\max}$,
and then choose $C_{\rm ba}$ sufficiently large, using
$\tau\le\lceil C_pK\log(R\kappa/r)\rceil$
and $K^{1-p/2}\le1$.
The accumulated errors are then bounded by a fixed
multiple of $R^p$.
Applying bounded stopping to the conditional potential
inequality gives
\begin{equation}
\label{eq:p2-stopped-penalty-energy}
\begin{aligned}
\mathbb E\Bigg[
&A_\tau
\bigl(
F(x_\tau)+L(\|x_\tau\|_p-R)_+-F(x^\star)
\bigr)\\
&+\frac{1-2^{1-q}}{2q}
\sum_{i=1}^{\tau}
\lambda_iA_i\|\widehat x_i-y_{i-1}\|_p^p
\Bigg]
\le C_{0,p}R^p.
\end{aligned}
\end{equation}
Both terms on the left are nonnegative.
Since the error constants are chosen after $C_{\rm ba}$
and the iteration budget, $C_{0,p}$ can be fixed
independently of those choices.

We next show that an execution is accepted
with constant probability.
On $\{A_\tau<A_{\max}\}$,
$\tau=\lceil C_pK\log(R\kappa/r)\rceil$.
If at least half the iterations have
$\lambda_i<2\lambda_{\min}$, the weight identity gives
\[
A_\tau^{1/p}-A_0^{1/p}
\ge
\frac{1}{p2^{(p-1)/p}}
\sum_{i=1}^{\tau}\lambda_i^{-1/p}
\ge
\frac{\tau}{C_p\lambda_{\min}^{1/p}}.
\]
For a sufficiently large iteration budget,
this implies $A_\tau\ge A_{\max}$, a contradiction.
Thus at least half the iterations have displacement
at least $3r/4$.
Using
\[
\lambda_iA_i
=
\frac1{2^{p-1}}\left(\frac{A_i}{a_i}\right)^p,
\qquad
\sum_{i=1}^{\tau}\frac{a_i}{A_i}
\le
\log\frac{A_\tau}{A_0}
\le C_p\log(R\kappa/r),
\]
H\"older's inequality yields
\[
\sum_{i=1}^{\tau}
\lambda_iA_i\|\widehat x_i-y_{i-1}\|_p^p
\ge
\frac{(3r/4)^p}{2^{p-1}}
\frac{(\tau/2)^{p+1}}
{\left(\sum_{i=1}^{\tau}a_i/A_i\right)^p}
\ge
\frac{c_pr^p\tau^{p+1}}{\log^p(R\kappa/r)}.
\]
Increasing the iteration budget makes the displacement
term in \eqref{eq:p2-stopped-penalty-energy}
exceed $4C_{0,p}R^p$ on $\{A_\tau<A_{\max}\}$.
Markov's inequality therefore gives
\[
\Pr(A_\tau\ge A_{\max})\ge\frac34.
\]
The same potential bound gives
\[
\mathbb E\left[
\bigl(
F(x_\tau)+L(\|x_\tau\|_p-R)_+-F(x^\star)
\bigr)
\mathbf1_{\{A_\tau\ge A_{\max}\}}
\right]
\le
\frac{C_{0,p}R^p}{A_{\max}}
\le\frac{\epsilon_{\rm opt}}{64},
\]
where $C_{\rm ba}\ge64C_{0,p}$.

To enforce deterministic bounds on the additional calls,
pre-sample independent geometric levels for all potential
iterations, revealing each only when its iteration is reached.
Level $0$ is always used.
For $t\ge1$, a level-$t$ call occurs only when the
sampled level is $t$ or $t+1$, with probability $O(2^{-t})$.
A Chernoff bound and a union bound over all levels
and executions give the per-execution bounds
\[
C_p2^{-t}K\log(R\kappa/r),
\qquad
0\le t\le\lceil\log_2K\rceil,
\]
except with total probability at most $1/(C_p\kappa)$.
Abort before exceeding a bound and treat aborted
executions as unsuccessful.
Combining this event with the distance-test events,
the implemented and analyzed executions agree except
with probability at most $1/(C_p\kappa)$,
after increasing constants.
Project each output onto $\mathbb B_p(R)$.
Its excess loss lies in $[0,2LR]$, so this comparison
adds at most $2LR/(C_p\kappa)$ to the expected excess loss.
Consequently, each implemented execution satisfies
\begin{align*}
\Pr(A_\tau\ge A_{\max})&\ge\frac12,\\
\mathbb E\Bigg[
\left(
F\left(
\frac{x_\tau}{\max\{1,\|x_\tau\|_p/R\}}
\right)-F(x^\star)
\right)
\mathbf1_{\{A_\tau\ge A_{\max}\}}
\Bigg]
&\le\frac{\epsilon_{\rm opt}}{32}.
\end{align*}
Independently restart from $0$ until the first accepted
execution and return its projected output.
After $\lceil\log_2(128\kappa)\rceil$ unsuccessful
executions, return $0$.
This selection uses no objective values, and
\[
\mathbb E[F(x)-F(x^\star)]
\le
\sum_{b=1}^{\lceil\log_2(128\kappa)\rceil}
2^{-(b-1)}\frac{\epsilon_{\rm opt}}{32}
+
LR\,2^{-\lceil\log_2(128\kappa)\rceil}
<
\epsilon_{\rm opt}.
\]
This proves the guarantee for every comparator
in $\mathbb B_p(R)$.
Finally, the total iteration budget is
\[
\lceil\log_2(128\kappa)\rceil
\lceil C_pK\log(R\kappa/r)\rceil
=
O(K\log\kappa\log(R\kappa/r)).
\]
Index the primal and level-$0$ calls by $j=1$,
and level $t\ge1$ by $j=t+1$.
The total number of iterations querying level $j$ is
\[O\left(2^{-j}K\log\kappa\log(R\kappa/r)
\right),\]
with objective accuracy
\[\frac{\lambda r^p}
{C_p2^j\log^{2p}(R\kappa/r)}.\]
We now prove the additional claim.
Choose
\[
x^\star\in
\operatorname*{argmin}_{x\in\mathbb R^d}F(x)
\cap\mathbb B_p(R).
\]
Run the same construction directly on $F$, initialize $A_0$
using its Lipschitz bound $L$, and define
\[
\widehat x_i
=
\operatorname*{argmin}_{x\in\mathbb R^d}
\left\{
F(x)+\frac{\lambda_i}{p}\|x-y_{i-1}\|_p^p
\right\}.
\]
All local oracle calls are now for $F$ itself.
The oracle-call estimates and the conditional potential
analysis remain unchanged up to constants depending only on $p$.
Since $x^\star$ is a global minimizer,
$F(x)-F(x^\star)\ge0$ for every $x\in\mathbb R^d$.
Let
\[
N=\left\lceil C_pK\log(R\kappa/r)\right\rceil.
\]
For the analyzed execution with exact distance tests,
the same bounded-stopping argument gives
\[
\begin{aligned}
\mathbb E\Bigg[
&A_\tau\bigl(F(x_\tau)-F(x^\star)\bigr)+\frac{1-2^{1-q}}{2q}
\sum_{i=1}^{\tau}
\lambda_iA_i\|\widehat x_i-y_{i-1}\|_p^p
\Bigg]
\le C_{0,p}R^p.
\end{aligned}
\]
Both terms are nonnegative.
The preceding acceptance argument therefore gives
\[
\Pr(A_\tau\ge A_{\max})\ge\frac34
\]
and, for sufficiently large $C_{\rm ba}$,
\[
\mathbb E\left[
\bigl(F(x_\tau)-F(x^\star)\bigr)
\mathbf1_{\{A_\tau\ge A_{\max}\}}
\right]
\le\frac{\epsilon_{\rm opt}}{64}.
\]

For both the implemented and analyzed executions,
ball-feasibility of the primal oracle outputs implies
$\|x_i\|_p\le R+ir$ for every completed iteration.
Consequently, for every possible output with $\tau\le N$,
\[
0\le F(x_\tau)-F(x^\star)
\le L\|x_\tau-x^\star\|_p
\le L(2R+Nr).
\]
Strengthen the coupling failure probability over all planned
executions to at most
\[
\frac{1}{C_p\kappa(2+Nr/R)}.
\]
This preserves the stated asymptotic oracle-call budgets.
Indeed, $Kr\le R$ implies
\[
2+Nr/R=O(\log(R\kappa/r)),
\]
so the additional logarithm in the amplification requirement
is absorbed by the existing
$O(\log(R\kappa/r))$ repetitions.
The geometric-level caps are adjusted by constant factors
in the same way.

The contribution of the coupling failures is now at most
\[L(2R+Nr)\frac{1}{C_p\kappa(2+Nr/R)}=
\frac{LR}{C_p\kappa}
=\frac{\epsilon_{\rm opt}}{C_p}.\]
Thus, after increasing constants, each implemented execution
satisfies
\[\Pr(A_\tau\ge A_{\max})\ge\frac12\]
and
\[\mathbb E\left[
\bigl(F(x_\tau)-F(x^\star)\bigr)
\mathbf1_{\{A_\tau\ge A_{\max}\}}\right]
\le\frac{\epsilon_{\rm opt}}{32}.\]
Independently restart until the first accepted execution,
returning its unprojected output.
After $\lceil\log_2(128\kappa)\rceil$ unsuccessful executions,
return $0$.
Since $F(0)-F(x^\star)\le LR$, the same restart calculation
proves
\[
\mathbb E[F(x)-F(x^\star)]<\epsilon_{\rm opt}.
\]
The total iteration and oracle-call bounds are unchanged
up to constants depending only on $p$.
Since a global minimizer has objective value no larger
than any comparator in $\mathbb B_p(R)$, the guarantee
also holds for every such comparator.
\end{proof}

\begin{proposition}[based on results in \cite{ghadimi2012optimal}]
\label{proposition:p-acsa-broo}
Let $f,g$ satisfy the assumptions of Problem~1. For $\bar x\in\mathcal X$, interpret local minimization as being over
$\mathcal X\cap\mathbb B_p(\bar x,r)$, as in
Proposition~\ref{proposition:prop2_equivalent}.
When $\rho=r$, Algorithm~\ref{alg:ac-sa}
is a $(\phi,\lambda)$-ball optimization oracle for $\widehat f_\rho$
for $0<\phi\le\lambda r^p/2$. It makes
\[O\!\left(1+\log\frac{\lambda r^p}{\phi}+\left(\frac{\left(C_pL/\rho+(p-1)\lambda r^{p-2}\right)^p}{\lambda^2\phi^{p-2}}\right)^{\!\frac1{p+2}}+\frac{L^p}{\lambda\phi^{p-1}}\right)\]
total queries to $g$ with constant query depth, and has an additional
computational cost of
\[O\!\left(1+\log\frac{\lambda r^p}{\phi}+\left(\frac{\left(C_pL/\rho+(p-1)\lambda r^{p-2}\right)^p}{\lambda^2\phi^{p-2}}\right)^{\!\frac1{p+2}}\right)\]
sequential proximal-update depth and
\[O\!\left(1+\log\frac{\lambda r^p}{\phi}+\left(\frac{\left(C_pL/\rho+(p-1)\lambda r^{p-2}\right)^p}{\lambda^2\phi^{p-2}}
\right)^{\!\frac1{p+2}}+\frac{L^p}{\lambda\phi^{p-1}}\right)\]
vector-arithmetic work. 
\end{proposition}
\begin{proof}
We follow the proof of Proposition~4 in Appendix~D of
\cite{carmon2023resqueing}.
Let \[L_\Psi=\frac{C_pL}{\rho}+(p-1)\lambda r^{p-2}, \qquad
\Psi(x)=\widehat f_\rho(x)
+\frac{\lambda}{p}\|x-\bar x\|_p^p.\]
The function $\Psi$ is $L_\Psi$-smooth on $\mathcal X\cap\mathbb B_p(\bar x,r)$ and
$\lambda$-uniformly convex relative to \[V(x,y)=D_{\|\cdot-\bar x\|_p^p/p}(y,x),\] with
$V(x,y)\ge\nu\|y-x\|_p^p/p$
and $\nu=2^{2-p}$. Let
\[x^\star\in\arg\min_{x\in
\mathcal X\cap\mathbb B_p(\bar x,r)}\Psi(x),\qquad \phi_k=\frac{\lambda r^p}{2^k}.\]

To justify the minibatch moment term below, condition on the past
before a minibatch is used, and let $x$ be its current iterate.
Since $\rho=r$ and $x\in\mathcal X\cap\mathbb B_p(\bar x,r)$,
Lemma~\ref{lem:GREG-q moment} implies that the $N_k$ independent
$p$-GREG estimates
$\widetilde{\nabla}_{\bar x}^{g,j}\widehat f_\rho(x)$,
$j=1,\ldots,N_k$, are unbiased and have $q$-th moments at most $3L^q$.
Since $q=p/(p-1)\in(1,2]$, applying the von Bahr--Esseen inequality
coordinatewise and summing over coordinates gives
\[\begin{aligned}
&\mathbb E\left\|
\frac{1}{N_k}\sum_{j=1}^{N_k}
\left(
\widetilde{\nabla}_{\bar x}^{g,j}\widehat f_{\rho,p}(x)
-\nabla\widehat f_{\rho,p}(x)
\right)\right\|_q^q\le
\frac{C_p}{N_k^q}
\sum_{j=1}^{N_k}
\mathbb E\left\|
\widetilde{\nabla}_{\bar x}^{g,j}\widehat f_{\rho,p}(x)
-\nabla\widehat f_{\rho,p}(x)
\right\|_q^q
\le \frac{C_pL^q}{N_k^{q-1}}.
\end{aligned}\]
Applying Theorem 16 (Equation (72)) conditionally on the initial point of epoch k, with $\mu=\lambda$ and
the parameters of Algorithm~\ref{alg:ac-sa},
gives
\begin{align}
\mathbb E[\Psi(x_{T_k}^{\mathrm{ag}})
-\Psi(x^\star)]\le C_p\left[
\frac{L_\Psi^{p-1}
\mathbb E V(x_0,x^\star)}
{\nu M_{f,k}^{p-2}T_k^p}
+
\frac{M_{f,k}^q+L^q/N_k^{q-1}}
{(\nu\lambda)^{q-1}T_k^{q-1}}
\right].
\label{eq:acsa-prop4-short}
\end{align}
We prove by induction that
\[\mathbb E[\Psi(x_{T_k}^{\mathrm{ag}})
-\Psi(x^\star)]\le\phi_k.\]
For the first epoch,
$x_0=\bar x$, and hence
$V(x_0,x^\star)\le r^p/p$.
For subsequent epochs, first-order optimality at the constrained
minimizer $x^\star$ and relative uniform convexity give
\[\begin{aligned}
\Psi(x_0)-\Psi(x^\star) &\ge
\langle\nabla\Psi(x^\star),x_0-x^\star\rangle
+\lambda V(x^\star,x_0)\ge \lambda V(x^\star,x_0).
\end{aligned}\]
By proof of Lemma 22
and Fact~3, we have, for all $a,b\in\mathbb R$,
\[D_{|\cdot|^p/p}(a,b)\le C_pD_{|\cdot|^p/p}(b,a),\]
where $C_p$ depends only on $p$.
Translating by $\bar x$ and summing coordinatewise therefore gives
\[\begin{aligned}
V(x_0,x^\star)=D_{\|\cdot-\bar x\|_p^p/p}(x^\star,x_0)\le C_pD_{\|\cdot-\bar x\|_p^p/p}(x_0,x^\star)=C_pV(x^\star,x_0).
\end{aligned}\]
Consequently,
\[\mathbb E V(x_0,x^\star) \le\frac{C_p}{\lambda}
\mathbb E[\Psi(x_0)-\Psi(x^\star)]
\le\frac{C_p\phi_{k-1}}{\lambda}.\]
Substituting the choices of
$M_{f,k},T_k,N_k$ from
Algorithm~\ref{alg:ac-sa} into
\eqref{eq:acsa-prop4-short}, and choosing
the $p$-dependent constants sufficiently
large, yields
\[\mathbb E[\Psi(x_{T_k}^{\mathrm{ag}})
-\Psi(x^\star)]
\le\frac{\phi_{k-1}}2=\phi_k.\]
Since
$K=\lceil\log_2(\lambda r^p/\phi)\rceil$,
the final output satisfies
\[\mathbb E[\Psi(x_{T_K}^{\mathrm{ag}})
-\Psi(x^\star)]\le\frac{\lambda r^p}{2^K}\le\phi.\]
All iterates lie in
$\mathcal X\cap\mathbb B_p(\bar x,r)$
by convexity of the feasible set.
Finally, the geometric decrease of $\phi_k$
and $p>2$ give
\[\sum_{k=1}^KT_k=O\left(
1+\log\frac{\lambda r^p}{\phi}
+\left(\frac{L_\Psi^p}{\lambda^2\phi^{p-2}}
\right)^{\frac1{p+2}}\right)\]
and
\[\sum_{k=1}^KT_kN_k=O\left(
1+\log\frac{\lambda r^p}{\phi}
+\left(\frac{L_\Psi^p}{\lambda^2\phi^{p-2}}
\right)^{\frac1{p+2}}
+\frac{L^p}{\lambda\phi^{p-1}}
\right).\]
All original stochastic-gradient queries
are issued in parallel, since the anchor
$\bar x$ and the batch sizes are fixed
in advance. This proves the claimed
oracle-query and proximal-update bounds.
\end{proof}
\subsection{EpochDA Solver}
\begin{algorithm}[H]
\caption{ParallelEpochDA$(f,g,\bar{x},r,\rho,\lambda,\bar{\phi})$}
\label{alg:epochda}
\begin{algorithmic}[1]
\State \textbf{Input:}
$f:\mathbb{R}^d\to\mathbb{R}$ and
$g:\mathbb{R}^d\to\mathbb{R}^d$
satisfying the assumptions of
Problem~\ref{problem:parallel_problem},
$\bar{x}\in\mathbb{R}^d$,
$r=\rho>0$, and $\lambda,\phi>0$
\State
$\eta_1\gets\frac{2^{2-p}}{4^{p-1}\lambda}$, \quad
$T_1\gets\left\lceil\frac{4^p}{(2^{2-p})^2}\right\rceil$,
\quad
$T\gets\left\lceil
\frac{C_pL^p}{\lambda \phi^{p-1}}\right\rceil$ where $C_p \geq \frac{2\cdot4^p\,2^{p-1}\,3^p}{(2^{2-p})^2(2^{p-1}-1)}$
\State Sample
$\xi_i\sim\mathcal{GG}_p(0,\rho^2I_d)$,
$i\in[2T]$, independently
\State Query $g(\bar{x}+\xi_i)$ for all $i\in[2T]$
in parallel, using independent oracle randomness
\State $x_1^0\gets\bar{x}$, \quad $k\gets1$
\While{$\sum_{j=1}^kT_j\le T$}
    \State $x_k^1\gets x_k^0$, \quad $\Phi_k(x)\gets\frac{1}{p}\|x-x_k^0\|_p^p$
    \For{$t\in[T_k]$}
        \State $i\gets\sum_{j=1}^{k-1}T_j+t$
        \State
        $\widetilde{\nabla}_{\bar{x}}^g\widehat{f}_\rho(x_k^t)\gets\frac{\gamma_{\rho, p}(x_k^t-\bar{x}-\xi_i)}     {\gamma_{\rho, p}(\xi_i)}
        g(\bar{x}+\xi_i)$
        \State
        $x_k^{t+1}\gets
        \operatorname*{argmin}_{x\in\mathbb{B}_p(\bar{x},r)}\left\{
        \eta_k\left\langle    \widetilde{\nabla}_{\bar{x}}^g\widehat{f}_\rho(x_k^t),x\right\rangle
        +\frac{\eta_k\lambda}{p}\|x-\bar{x}\|_p^p
        +D_{\Phi_k}(x,x_k^t)
        \right\}$
    \EndFor
    \State
    $x_{k+1}^0\gets
    \frac{1}{T_k}\sum_{t=1}^{T_k}x_k^{t+1}$, \quad
    $T_{k+1}\gets\left\lceil2^{p-1}T_k\right\rceil$, \quad
    $\eta_{k+1}\gets\frac{\eta_k}{2^{p-1}}$, \quad
    $k\gets k+1$
\EndWhile
\State \Return $x_k^0$
\end{algorithmic}
\end{algorithm}

Throughout, let $f,g$ satisfy Problem~\ref{problem:parallel_problem},
and let
\[
x^\star=\operatorname*{argmin}_{x\in\mathbb B_p(\bar x,r)}
\left\{\widehat f_\rho(x)+\frac\lambda p\|x-\bar x\|_p^p\right\}.
\]
Let $C_{\mathrm{res},p}\ge1$ be the constant in the
conditional $q$-moment bound for $\rho=r$.

\begin{lemma}
In view of Problem~\ref{problem:parallel_problem}, \label{lem:epochda-bregman-single-epoch}
Each epoch $k$ of Algorithm~\ref{alg:epochda} satisfies
\begin{equation}
\label{eq:epochda-bregman-single-epoch}
\begin{split}
&\mathbb E\left[\widehat f_\rho(x_{k+1}^0)
+\frac\lambda p\|x_{k+1}^0-\bar x\|_p^p\right]
-\widehat f_\rho(x^\star)
-\frac\lambda p\|x^\star-\bar x\|_p^p\le
\frac{\mathbb E\|x^\star-x_k^0\|_p^p}{p\eta_kT_k}
+\frac{(4L)^q}{q}
 \left(\frac{\eta_k}{2^{2-p}}\right)^{q-1}.
\end{split}
\end{equation}
\end{lemma}

\begin{proof}
Problem~\ref{problem:parallel_problem} gives $\|\mathbb E g(x)\|_q\le\mathbb E\|g(x)\|_q\le L$; hence $f$ and $\widehat f_\rho$ are convex and $L$-Lipschitz.
Let $\mathcal F_{k,t-1}$ be the history of the samples already
used, including $x_k^t$. By \Cref{lem:GREG-q moment},
\begin{align}
\mathbb E\!\left[
\widetilde\nabla^g_{\bar x}
\widehat f_\rho(x_k^t)
\mid\mathcal F_{k,t-1}\right]
=\nabla\widehat f_\rho(x_k^t) \quad,\quad \mathbb E\!\left[
\left\|\widetilde\nabla^g_{\bar x}\widehat f_\rho(x_k^t)
\right\|_q^q\mid\mathcal F_{k,t-1}\right] \le 3L^q
\label{eq:epochda-bregman-oracle}
\end{align}
The optimality condition of the update, convexity of the
regularizer, and the Bregman three-point identity imply
\begin{align}
&\left\langle\widetilde\nabla^g_{\bar x}\widehat f_\rho(x_k^t),x_k^{t+1}-x^\star\right\rangle+\frac\lambda p\|x_k^{t+1}-\bar x\|_p^p
-\frac\lambda p\|x^\star-\bar x\|_p^p\notag\\
&\quad\le\frac1{\eta_k}\bigl[
D_{\Phi_k}(x^\star,x_k^t)
-D_{\Phi_k}(x^\star,x_k^{t+1})
-D_{\Phi_k}(x_k^{t+1},x_k^t)\bigr]
\label{eq:epochda-bregman-three-point}
\end{align}
By Lipschitz continuity, convexity, and conditional unbiasedness,
\begin{align*}
\mathbb E[\widehat f_\rho(x_k^{t+1})-\widehat f_\rho(x^\star)] &\le\mathbb E\left[\left\langle\widetilde\nabla^g_{\bar x} \widehat f_\rho(x_k^t),x_k^t-x^\star\right\rangle
+L\|x_k^{t+1}-x_k^t\|_p\right]\\
&\le\mathbb E\left\langle
\widetilde\nabla^g_{\bar x}\widehat f_\rho(x_k^t),
x_k^{t+1}-x^\star\right\rangle
+\mathbb E\left[\left(\left
\|\widetilde\nabla^g_{\bar x}
\widehat f_\rho(x_k^t)\right\|_q+L\right)
\|x_k^{t+1}-x_k^t\|_p\right]
\end{align*}
Since $D_{\Phi_k}(u,v)\ge2^{2-p}\|u-v\|_p^p/p$,
Young's inequality gives
\begin{align*}
\left(\left\|\widetilde\nabla^g_{\bar x}\widehat f_\rho(x_k^t)\right\|_q+L\right)
\|x_k^{t+1}-x_k^t\|_p
-\frac{D_{\Phi_k}(x_k^{t+1},x_k^t)}{\eta_k}
\le\frac1q \left(\frac{\eta_k}{2^{2-p}}\right)^{q-1}
\left(\left\|\widetilde\nabla^g_{\bar x}\widehat
f_\rho(x_k^t)\right\|_q+L\right)^q
\end{align*}
Moreover, conditional Minkowski's inequality and
\eqref{eq:epochda-bregman-oracle} yield
\[\mathbb E\!\left[
\left(\left\|\widetilde\nabla^g_{\bar x}\widehat f_\rho(x_k^t)\right\|_q+L\right)^q\mid\mathcal F_{k,t-1}\right]
\le(4L)^q.\]
Combining these estimates with
\eqref{eq:epochda-bregman-three-point} gives
\begin{align*}
&\mathbb E\left[\widehat f_\rho(x_k^{t+1})
+\frac\lambda p\|x_k^{t+1}-\bar x\|_p^p\right]
-\widehat f_\rho(x^\star)
-\frac\lambda p\|x^\star-\bar x\|_p^p\\&\quad\le
\frac{\mathbb E D_{\Phi_k}(x^\star,x_k^t)
-\mathbb E D_{\Phi_k}(x^\star,x_k^{t+1})}{\eta_k}
+\frac{(4L)^q}{q}
\left(\frac{\eta_k}{2^{2-p}}\right)^{q-1}
\end{align*}
Summing over $t=1,\ldots,T_k$ telescopes the Bregman terms
because $\Phi_k$ is fixed within the epoch. As $x_k^1=x_k^0$,
\[D_{\Phi_k}(x^\star,x_k^1)
=\frac1p\|x^\star-x_k^0\|_p^p.\]
Dropping the nonnegative final divergence, dividing by $T_k$,
and using convexity at
$x_{k+1}^0=T_k^{-1}\sum_{t=1}^{T_k}x_k^{t+1}$ proves
\eqref{eq:epochda-bregman-single-epoch}.
\end{proof}

\begin{lemma}
\label{lem:epochda-bregman-epoch-induction}
After each completed epoch $k$ of Algorithm~\ref{alg:epochda},
\begin{equation}
\label{eq:epochda-bregman-ending-gap}
\begin{split}
\mathbb E\left[\widehat f_\rho(x_{k+1}^0)+\frac\lambda p\|x_{k+1}^0-\bar x\|_p^p\right]
-\widehat f_\rho(x^\star)
-\frac\lambda p\|x^\star-\bar x\|_p^p\le\frac{(4L)^q}{\lambda^{q-1}2^k}.
\end{split}
\end{equation}
\end{lemma}

\begin{proof}
Convexity of $\widehat f_\rho$, $p$-uniform convexity of the
regularizer, and constrained optimality of $x^\star$ imply,
for every $w\in\mathbb B_p(\bar x,r)$,
\begin{equation}
\label{eq:epochda-bregman-gap-distance}
\begin{split}
&\widehat f_\rho(w)+\frac\lambda p\|w-\bar x\|_p^p -\widehat f_\rho(x^\star) -\frac\lambda p\|x^\star-\bar x\|_p^p\ge
\frac{2^{2-p}\lambda}{p}\|w-x^\star\|_p^p.
\end{split}
\end{equation}
Since $x_1^0=\bar x$ and $\widehat f_\rho$ is $L$-Lipschitz,
\begin{align*}
\widehat f_\rho(\bar x)-\widehat f_\rho(x^\star)
-\frac\lambda p\|x^\star-\bar x\|_p^p\le\sup_{a\ge0}\left\{La-\frac\lambda p a^p\right\}
=\frac{L^q}{q\lambda^{q-1}}
\le\frac{(4L)^q}{\lambda^{q-1}}.
\end{align*}
The parameter choices in Algorithm~\ref{alg:epochda} give
\[\eta_k=\frac{2^{2-p}}
{4^{p-1}\lambda\,2^{(p-1)(k-1)}},\qquad T_k\ge \frac{4^p}{(2^{2-p})^2}\,2^{(p-1)(k-1)},\]
and hence $\eta_kT_k\ge4/(2^{2-p}\lambda)$.
Suppose by induction that the expected objective error at the start of epoch $k$ is at most
$(4L)^q/(\lambda^{q-1}2^{k-1})$,
the base case $k=1$ is verified above. By \eqref{eq:epochda-bregman-gap-distance},
\[\frac{\mathbb E\|x^\star-x_k^0\|_p^p}{p\eta_kT_k}
\le\frac14\frac{(4L)^q}
{\lambda^{q-1}2^{k-1}}.\]
Using $(p-1)(q-1)=1$, the other term in
Lemma~\ref{lem:epochda-bregman-single-epoch} is
\[\frac{(4L)^q}{q}
\left(\frac{\eta_k}{2^{2-p}}
\right)^{q-1}=\frac{(4L)^q}
{4q\lambda^{q-1}2^{k-1}}\]
Applying Lemma~\ref{lem:epochda-bregman-single-epoch} therefore gives
\begin{align*}
\mathbb E\left[\widehat f_\rho(x_{k+1}^0)+\frac\lambda p\|x_{k+1}^0-\bar x\|_p^p\right]-\widehat f_\rho(x^\star)-\frac\lambda p\|x^\star-\bar x\|_p^p
\le\left(\frac14+\frac1{4q}\right)\frac{(4L)^q}
{\lambda^{q-1}2^{k-1}}\le
\frac{(4L)^q}{\lambda^{q-1}2^k}.
\end{align*}
This is the required bound for epoch $k$ and the induction
hypothesis for the start of epoch $k+1$.
\end{proof}

\begin{proposition}[Proposition 3,  \cite{carmon2023resqueing}]\label{prop:proposition3}
Let $f,g$ satisfy the assumptions of \Cref{problem:parallel_problem}. When $\rho = r$, \Cref{alg:epochda} is a $(\phi,\lambda)$-ball optimization oracle for $\widehat f_\rho$ which makes
$O\left(\frac{L^p}{\lambda\phi^{p-1}}\right)$
total queries to $g$ with constant query depth, and an additional computational cost of
$O\left(\frac{L^p}{\lambda\phi^{p-1}}\right)$ depth and work.
\end{proposition}
\begin{proof}
Let $k$ be the epoch index when \Cref{alg:epochda} terminates, so that
the returned point is $x_k^0$. The recurrence
$T_{j+1}=\lceil2^{p-1}T_j\rceil$
implies
\[T_j\le 2^{(p-1)(j-1)}
\left(T_1+\frac{1}{2^{p-1}-1}
\right)\le \frac{2\cdot4^p}{(2^{2-p})^2}\,2^{(p-1)(j-1)}.\]
By the stopping condition,
\begin{align*}
T<\sum_{j=1}^{k}T_j\le
\frac{2\cdot4^p}{(2^{2-p})^2}
\frac{2^{(p-1)k}-1}{2^{p-1}-1}
\le\frac{2\cdot4^p\,2^{p-1}}
{(2^{2-p})^2(2^{p-1}-1)}\,2^{(p-1)(k-1)}.
\end{align*}
The initial-error bound and
Lemma~\ref{lem:epochda-bregman-epoch-induction} therefore give
\begin{align*}
&\mathbb E\left[
\widehat f_\rho(x_k^0)
+\frac{\lambda}{p}\|x_k^0-\bar x\|_p^p\right]-\min_{x\in\mathbb B_p(\bar x,r)}
\left\{\widehat f_\rho(x)
+\frac{\lambda}{p}\|x-\bar x\|_p^p\right\}\\
&\qquad\le
\frac{(4L)^q}{\lambda^{q-1}2^{k-1}}\le\left[\frac{2\cdot4^p\,2^{p-1}\,3^p}
{(2^{2-p})^2(2^{p-1}-1)}
\frac{L^p}{\lambda T}
\right]^{1/(p-1)}\le\phi,
\end{align*}
where the last inequality follows from $T=\lceil C_pL^p/(\lambda\phi^{p-1})
\rceil$.
The output lies in $\mathbb B_p(\bar x,r)$ by construction
and convexity, so the preceding error bound establishes
the required guarantee.
All $2T=O(1+L^p/(\lambda\phi^{p-1}))$ queries can be made
in one parallel round.
The algorithm performs at most $T$ Bregman proximal updates
and $O(T)$ vector operations, giving $O(T)$ additional
computational depth and work if each update and vector
operation has constant depth and work.
\end{proof}

\begin{lemma}[Lemma 2.3, \cite{chakrabarty2023parallel}]
\label{lem:radial-penalty}
Let $f,g$ satisfy the assumptions of
\Cref{problem:parallel_problem}, with $L>0$.
For any $c\in\mathbb R^d$ and $R>0$, the function \[x\mapsto f(x)+2L(\|x-c\|_p-R)_+\]
is convex and $3L$-Lipschitz with respect to $\|\cdot\|_p$.
For any $y$ with $\|y-c\|_p\ge R$, define
\[y^{c,R}:=c+\frac{R}{\|y-c\|_p}(y-c).\]
Then $\|y^{c,R}-c\|_p=R$ and
\[f(y^{c,R})\le f(y)+2L(\|y-c\|_p-R)_+-L(\|y-c\|_p-R).\]
The penalized function attains a global minimum, and every
global minimizer $x^\star$ satisfies
\[\|x^\star-c\|_p\le R,
\qquad f(x^\star)+2L(\|x^\star-c\|_p-R)_+ =f(x^\star) =\min_{x\in\mathbb B_p(c,R)}f(x).\]
Moreover, adding a deterministic subgradient of
$2L(\|x-c\|_p-R)_+$ to $g(x)$ gives a stochastic
subgradient oracle for the penalized function with
$q$-th moment at most $(4L)^q$, using one query to $g$.
\end{lemma}
\begin{proof}
The proof follows the argument of Lemma~2.3 in
\cite{chakrabarty2023parallel}.
We first show that $f(x)+2L(\|x-c\|_p-R)_+$ is convex
and $3L$-Lipschitz.
The assumptions of \Cref{problem:parallel_problem} give
\[\|\mathbb E g(x)\|_q \le \bigl(\mathbb E\|g(x)\|_q^q\bigr)^{1/q}\le L.\]
Since $\mathbb E g(x)\in\partial f(x)$, the subgradient
inequality implies that $f$ is $L$-Lipschitz with respect
to $\|\cdot\|_p$.
The function $x\mapsto 2L(\|x-c\|_p-R)_+$ is convex
and $2L$-Lipschitz.
Adding $f$ therefore gives a convex and $3L$-Lipschitz
function.
Next, consider $y\in\mathbb R^d$ with $\|y-c\|_p\ge R$.
By the definition of $y^{c,R}$,
\[\|y^{c,R}-c\|_p=R, \qquad
\|y^{c,R}-y\|_p=\|y-c\|_p-R.\]
In particular, $2L(\|y^{c,R}-c\|_p-R)_+=0$.
Moreover,
\begin{align*}
f(y)+2L(\|y-c\|_p-R)_+ &=
f(y)+2L(\|y-c\|_p-R)\\ &\ge
f(y^{c,R})-L\|y^{c,R}-y\|_p
+2L(\|y-c\|_p-R)\\ &=
f(y^{c,R})+L(\|y-c\|_p-R),
\end{align*}
where the inequality follows from the $L$-Lipschitz
continuity of $f$.
This yields the desired result, as it implies that the value
of $f(x)+2L(\|x-c\|_p-R)_+$ at any point outside the ball
is larger than its value at some point satisfying
$\|x-c\|_p\le R$, where
$f(x)+2L(\|x-c\|_p-R)_+=f(x)$.
Further, since the ball is compact and $f$ is continuous,
the minimum of $f$ over the ball is attained.
For the stochastic oracle, define
\[v(x):=\begin{cases}
0, & \|x-c\|_p\le R,\\[2mm]
\displaystyle \frac{(|x_i-c_i|^{p-2}(x_i-c_i))_{i=1}^d}
{\|x-c\|_p^{p-1}},& \|x-c\|_p>R.
\end{cases}\]
Then $v(x)\in\partial[(\|\cdot-c\|_p-R)_+](x)$
and $\|v(x)\|_q\le1$, so
\[\mathbb E[g(x)+2Lv(x)]\in
\partial\bigl[f+2L(\|\cdot-c\|_p-R)_+\bigr](x).\]
By Minkowski's inequality,
\[\bigl(\mathbb E\|g(x)+2Lv(x)\|_q^q\bigr)^{1/q}
\le \bigl(\mathbb E\|g(x)\|_q^q\bigr)^{1/q}
+2L\|v(x)\|_q\le3L.\]
Computing $g(x)+2Lv(x)$ requires only one query to $g$.
\end{proof}

\begin{theorem}[Parallel EpochDA-based solver]
\label{theorem:epochda_parallel_guarantee}
\textsc{BallAccel} (Proposition~\ref{proposition:prop2_equivalent})
using parallel \textsc{EpochDA} (Algorithm \ref{alg:epochda}) as a ball optimization
oracle solves Problem \ref{problem:parallel_problem} with expected error
$\epsilon_{\rm opt}$, with
\[O\!\left(d^{\frac{p}{2(p+1)}}\kappa^{\frac{p}{p+1}}
\log^4(d\kappa)\right)
\text{ query depth}\]
and
\[O\!\left(d^{\frac{p}{2(p+1)}}\kappa^{\frac{p}{p+1}}
\log^4(d\kappa)+d^{\frac{p(p-2)}{2(p+1)}}\kappa^{\frac{p(2p-1)}{p+1}}\log^{2p(p-1)+2}(d\kappa)\right)
\text{ total queries},\]
and an additional computational cost of
\[O\!\left(d^{\frac{p}{2(p+1)}}\kappa^{\frac{p}{p+1}}
\log^4(d\kappa)+d^{\frac{p(p-2)}{2(p+1)}}\kappa^{\frac{p(2p-1)}{p+1}}\log^{2p(p-1)+2}(d\kappa)\right)\text{ depth}\]
and
\[O\!\left(\left(d^{\frac{p}{2(p+1)}}\kappa^{\frac{p}{p+1}}
\log^4(d\kappa)+d^{\frac{p(p-2)}{2(p+1)}}\kappa^{\frac{p(2p-1)}{p+1}}\log^{2p(p-1)+2}(d\kappa)
\right)\cdot d\right)\text{ work}.\]
\end{theorem}
\begin{proof}
Let
$x^\star\in\argmin_{x\in\mathbb B_p(R)}f(x)$, $\widehat x$ be the point returned by \textsc{BallAccel}
and define
\[f_{\rm pen}(x)=f(x)+2L(\|x\|_p-R)_+.\]
Let $F=\widehat{(f_{\mathrm{pen}})}_\rho$ and choose
$r=\rho=\frac{\epsilon_{\rm opt}}{C_pL\sqrt d}$,
where $C_p$ is a sufficiently large constant depending
only on $p$.
Since $f_{\rm pen}$ is $3L$-Lipschitz and
$\mathbb E\|\xi_\rho\|_p\le\rho\sqrt d$, this choice gives
\[
|F(x)-f_{\rm pen}(x)|
\le 3L\mathbb E\|\xi_\rho\|_p
\le \frac{\epsilon_{\rm opt}}4
\]
for every $x\in\mathbb R^d$.
We now show that $F$ has a minimizer in $\mathbb B_p(3R)$.
To see this, let $y\in\mathbb R^d$ satisfy $\|y\|_p>3R$, and set $\widetilde y=\frac{R}{\|y\|_p}y$. By Lemma~\ref{lem:radial-penalty},
\begin{align*}
F(y) &\ge f_{\rm pen}(y)-\frac{\epsilon_{\rm opt}}4\ge f_{\rm pen}(\widetilde y)
+L(\|y\|_p-R)-\frac{\epsilon_{\rm opt}}4 \\
&\ge F(\widetilde y) +L(\|y\|_p-R)-\frac{\epsilon_{\rm opt}}2 >
F(\widetilde y),
\end{align*}
where the last inequality uses
$\|y\|_p>3R$ and $\epsilon_{\rm opt}\le LR$.
Since $F$ is continuous, it attains its minimum over
$\mathbb B_p(3R)$, and the  inequality shows
that this is also a global minimum.
Since $F$ has a global minimizer in $B_p(3R)$, we use the
unpenalized implementation in Proposition~\ref{proposition:prop2_equivalent}
with radius $3R$, Lipschitz bound $3L$, and target accuracy
$\epsilon_{\mathrm{opt}}/2$. 
Let $g_{\mathrm{pen}}$ be the stochastic subgradient oracle for
$f_{\mathrm{pen}}$ from Lemma~\ref{lem:radial-penalty}. Each
evaluation of $g_{\mathrm{pen}}$ uses one query to $g$, and
$\mathbb E\|g_{\mathrm{pen}}(z)\|_q^q\le (4L)^q$. Thus, with
$L_{\mathrm{pen}}=4L$, Proposition~\ref{prop:proposition3}
implements each ball optimization call for $F$ using
Algorithm~\ref{alg:epochda} on
$(f_{\mathrm{pen}},g_{\mathrm{pen}})$. The local subproblem
has objective
\[
x\mapsto F(x)+\frac{\lambda}{p}\|x-\bar x\|_p^p,
\qquad x\in B_p(\bar x,r).
\]
Rescaling the Lipschitz parameter, radius,
and target accuracy changes the corresponding
BallAccel parameters only by constant factors
depending on $p$.
Hence, we retain the notation
\[
\kappa=\frac{LR}{\epsilon_{\rm opt}},
\qquad
K=\left(\frac Rr\right)^{\frac p{p+1}}
=\Theta\left(
d^{\frac{p}{2(p+1)}}\kappa^{\frac p{p+1}}
\right),
\qquad
\lambda^\star=
\frac{\epsilon_{\rm opt}K^p}{R^p}\log^p\kappa
\]
in the following complexity bounds, with
constant-factor changes absorbed into the
$O(\cdot)$ notation.
For the query depth, each ball optimization oracle
implemented by Algorithm~\ref{alg:epochda}
has constant query depth.
Proposition~\ref{proposition:prop2_equivalent}
requires at most
$O(K\log\kappa\log(R\kappa/r))$ iterations,
with at most $O(\log^2(R\kappa/r))$
ball optimization oracle calls per iteration.
Since $\log(R\kappa/r)=O(\log(d\kappa))$,
the query depth is
\[O\left( K\log\kappa
\log^3\left(\frac{R\kappa}{r}\right)\right)=O\left(
d^{\frac{p}{2(p+1)}}
\kappa^{\frac p{p+1}}
\log^4(d\kappa)\right).\]
By Proposition~\ref{prop:proposition3} with $L_{\mathrm{pen}}=4L$,
each $(\phi,\lambda)$-BROO uses
$O\left(1+\frac{L^p}{\lambda\phi^{p-1}}\right)$ queries to $g$.
We may upper bound every query cost using
$\lambda=\Omega_p(\lambda^\star)$.
Writing
$\ell:=\log(R\kappa/r)$ and
$j_{\max}:=\lceil\log_2K+C_{\rm ba}\rceil$,
and applying both requirements of
Proposition~\ref{proposition:prop2_equivalent},
the total number of stochastic-gradient queries
is at most
\begin{align*}
&O\Bigg(
K\log\kappa\,\ell^3
+\frac{K\log\kappa\,\ell^3L^p}
{(\lambda^\star)^pr^{p(p-1)}}
+K\log\kappa\,\ell
\sum_{j=1}^{j_{\max}}2^{-j}
\frac{L^p(2^j\ell^{2p})^{p-1}}
{(\lambda^\star)^pr^{p(p-1)}}
\Bigg)\\
&\qquad=
O\Bigg(
K\log\kappa\,\ell^3
+\frac{K\log\kappa\,\ell^3L^p}
{(\lambda^\star)^pr^{p(p-1)}}
+\frac{K\log\kappa\,\ell^{2p(p-1)+1}L^p}
{(\lambda^\star)^pr^{p(p-1)}}
\sum_{j=1}^{j_{\max}}2^{j(p-2)}
\Bigg).
\end{align*}
For $p\ge2$, the sum contains $O(\ell)$ terms,
each bounded by $O(K^{p-2})$, and hence
\[
\sum_{j=1}^{j_{\max}}2^{j(p-2)}
=O(K^{p-2}\ell).
\]
Moreover, by $K^{p+1}r^p=R^p$,
\[
\frac{KL^p}{(\lambda^\star)^pr^{p(p-1)}}
=\frac{\kappa^p}{\log^{p^2}\kappa}.
\]
Therefore, the total number of queries is bounded by
\begin{align*}
&O\Bigg(
K\log\kappa\,\ell^3
+\frac{\kappa^p\log\kappa\,\ell^3}
{\log^{p^2}\kappa}
+\frac{\kappa^pK^{p-2}\log\kappa\,
\ell^{2p(p-1)+2}}{\log^{p^2}\kappa}
\Bigg)\\
&\qquad=
O\left(
K\log^4(d\kappa)
+\kappa^pK^{p-2}
\log^{2p(p-1)+2}(d\kappa)
\right)\\
&\qquad=
O\Bigg(
d^{\frac{p}{2(p+1)}}\kappa^{\frac p{p+1}}
\log^4(d\kappa)
+
d^{\frac{p(p-2)}{2(p+1)}}
\kappa^{\frac{p(2p-1)}{p+1}}
\log^{2p(p-1)+2}(d\kappa)
\Bigg),
\end{align*}
where we used
$\log\kappa/\log^{p^2}\kappa\le1$ and $\kappa^pK^{p-2}=
\Theta\left(d^{\frac{p(p-2)}{2(p+1)}}\kappa^{\frac{p(2p-1)}{p+1}}\right)$.
Finally, Proposition~\ref{proposition:prop2_equivalent}
gives
\[
\mathbb EF(\widehat{x})
\le\min_{x\in\mathbb R^d}F(x)
+\frac{\epsilon_{\rm opt}}2
\le F(x^\star)+\frac{\epsilon_{\rm opt}}2.
\]
Since $f(x)\le f_{\rm pen}(x)$ and
$f_{\rm pen}$ is $3L$-Lipschitz,
\[
f(x)\le F(x)+3L\mathbb E\|\xi_\rho\|_p
\]
for all $x\in\mathbb R^d$.
Therefore,
\begin{align*}
\mathbb Ef(\widehat{x})
&\le
\mathbb EF(\widehat{x})
+3L\mathbb E\|\xi_\rho\|_p\\
&\le
F(x^\star)+\frac{\epsilon_{\rm opt}}2
+3L\mathbb E\|\xi_\rho\|_p\\
&\le
f(x^\star)+\frac{\epsilon_{\rm opt}}2
+6L\mathbb E\|\xi_\rho\|_p\\
&\le
f(x^\star)+\epsilon_{\rm opt}.
\end{align*}
\end{proof}

\subsection{$\ell_p$-AC-SA Solver}

\begin{algorithm}[H]
\caption{$\ell_p$-AC-SA$(f,\bar{x},r,\rho,\lambda,\bar{\phi})$}
\label{alg:ac-sa}
\begin{algorithmic}[1]
\Require $f:\mathbb{R}^d\to\mathbb{R}$ and $g:\mathbb{R}^d\to\mathbb{R}^d$ satisfying the assumptions of Problem~1, $\bar{x}\in\mathcal X$, $r,\rho,\lambda>0$, $0<\bar\phi\le\lambda r^p/2$,
$V(x,y)=D_{\|\cdot-\bar x\|_p^p/p}(y,x)$ and $\nu=2^{2-p}$.
\State $K\leftarrow\lceil\log_2(\lambda r^p/\bar\phi)\rceil$.
\State For $k\in[K]$, set
\[
\begin{aligned}
M_{f,k}&\leftarrow(\left(C_pL/\rho+(p-1)\lambda r^{p-2}\right)\nu\lambda\left(\lambda r^p/2^{k-1}\right)^p)^{1/(p+2)},\\
 T_k&\leftarrow\left\lceil C_p\left[1+
 \left(\frac{\left(C_pL/\rho+(p-1)\lambda r^{p-2}\right)^p}{(\nu\lambda)^2\left(\lambda r^p/2^{k-1}\right)^{p-2}}
 \right)^{1/(p+2)}\right]\right\rceil,\\
 N_k&\leftarrow\left\lceil C_p\max\left\{1,
 \frac{L^p}{\nu\lambda T_k\left(\lambda r^p/2^{k-1}\right)^{p-1}}\right\}\right\rceil.
\end{aligned}
\]
\State Sample $\xi_i\sim \mathcal{GG}_p(0,\rho^2 I_d)$, $i\in[N]$ independently, for
$N=\sum_{k\in[K]}T_kN_k$.
\State Query $g(\bar{x}+\xi_i)$ for all $i\in[N]$ in parallel

\State $x_0^{\mathrm{ag}}\leftarrow \bar{x}$, $x_0\leftarrow \bar{x}$

\For{$k\in[K]$}
    \For{$t\in[T_k]$}
        \State
        $\alpha_t\leftarrow p/(t+p-1)$,
        $\Gamma_t\leftarrow\Gamma(p+1)\Gamma(t)/\Gamma(t+p)$.
        \State $\displaystyle
        \gamma_{t,k}\leftarrow
        \frac{2\left(\sup_{s\ge1}\alpha_s^p/\Gamma_s\right)
        \left(C_pL/\rho+(p-1)\lambda r^{p-2}\right)^{p-1}}
        {\nu M_{f,k}^{p-2}}\Gamma_t$.
        \State
        $x_t^{\mathrm{md}}
        \leftarrow\frac{(1-\alpha_t)(\lambda+\gamma_{t,k})}
        {\gamma_{t,k}+(1-\alpha_t^2)\lambda}
        x_{t-1}^{\mathrm{ag}}+
        \frac{\alpha_t[(1-\alpha_t)\lambda+\gamma_{t,k}]}
        {\gamma_{t,k}+(1-\alpha_t^2)\lambda}
        x_{t-1}.$
        \State
        $N_{T,[k-1]}\leftarrow
        \sum_{k'\in[k-1]}T_{k'}N_{k'}.$
        \State
        $\displaystyle\begin{aligned}
        \widehat{\nabla} f(x_t^{\mathrm{md}})\leftarrow{}&
        \frac{1}{N_k}\sum_{n\in[N_k]}
        \frac{\gamma_{\rho, p}(x_t^{\mathrm{md}}-\bar{x}
        -\xi_{N_{T,[k-1]}+(t-1)N_k+n})}
        {\gamma_{\rho, p}(\xi_{N_{T,[k-1]}+(t-1)N_k+n})}\\
        &\qquad\cdot g(\bar{x}+\xi_{N_{T,[k-1]}+(t-1)N_k+n}).
        \end{aligned}$
        \State $\displaystyle\begin{aligned}
        x_t\gets\arg\min_{x\in\mathcal X\cap B_p(\bar{x},r)}
        \Big\{&\langle\alpha_t[\widehat{\nabla}f(x_t^{\mathrm{md}})
        +\lambda J_p(x_t^{\mathrm{md}}-\bar{x})],x\rangle\\
        &+\alpha_t\lambda V(x_t^{\mathrm{md}},x)\\
        &+[(1-\alpha_t)\lambda+\gamma_{t,k}]V(x_{t-1},x)\Big\}.
        \end{aligned}$
        \State
        $x_t^{\mathrm{ag}}\leftarrow
        \alpha_t x_t+(1-\alpha_t)x_{t-1}^{\mathrm{ag}}.$
    \EndFor
    \State $x_0^{\mathrm{ag}}\leftarrow x_{T_k}^{\mathrm{ag}}$, $x_0\leftarrow x_{T_k}^{\mathrm{ag}}$
\EndFor
\State \Return $x_{T_K}^{\mathrm{ag}}$
\end{algorithmic}
\end{algorithm}

\begin{definition}
\label{def:acsa-problem}
Consider the class of stochastic convex optimization (SCO) problems
\begin{equation}\label{acsa-opt-problem}
\Psi^*:=\min_{x\in X}\{\Psi(x):=f(x)+\mathcal X(x)\},
\end{equation}
where $X\subseteq\mathbb R^d$ is a nonempty closed convex set,
$\mathcal X:X\to\mathbb R$ is a simple convex function, and
$f:X\to\mathbb R$ is convex. Suppose $\omega:X\to\mathbb R$ is differentiable
and its Bregman divergence
\[
V(x,y):=\omega(y)-\omega(x)-\langle\nabla\omega(x),y-x\rangle
\]
satisfies $V(x,y)\ge\nu\|y-x\|_p^p/p$ for some $\nu>0$.
Assume that $\lambda\ge0$, $\mu\ge0$, $L\ge0$, and $M_f>0$,
and that, for every $x,y\in X$ and $f'(x)\in\partial f(x)$,
\begin{equation}
\label{eq:acsa-model}
\mu V(x,y)
\le f(y)-f(x)-\langle f'(x),y-x\rangle
\le \frac L2\|y-x\|_p^2
+M_f\|y-x\|_p.
\end{equation}
Define the lower model
\begin{equation}\label{eq:acsa-lower-model}
l_\Psi(z,x):=f(z)+\langle f'(z),x-z\rangle
+\mu V(z,x)+\mathcal X(x).
\end{equation}
\end{definition}

\begin{lemma}
\label{lemma:acsa-helperlemma1}
    Let the convex function $p : X \to \mathbb{R}$, the points $\tilde{x}, \tilde{y} \in X$ and the scalars $\mu_1, \mu_2 \geq 0$ be given. Let $\omega : X \to \mathbb{R}$ be a differentiable convex function and $V(x,z)$ be defined in Definition~\ref{def:acsa-problem} (Bregman divergence). If
\[
u^* \in \operatorname{Argmin}\left\{
p(u) + \mu_1 V(\tilde{x},u) + \mu_2 V(\tilde{y},u) : u \in X
\right\},
\]
then for any $u \in X$, we have
\[
p(u^*) + \mu_1 V(\tilde{x},u^*) + \mu_2 V(\tilde{y},u^*)
\leq
p(u) + \mu_1 V(\tilde{x},u) + \mu_2 V(\tilde{y},u)
- (\mu_1+\mu_2)V(u^*,u).
\]
\end{lemma}

\begin{proof}
    The proof is the same as in \cite{ghadimi2012optimal}, it holds for any differentiable convex function $\omega$.
\end{proof}

\begin{lemma}
\label{lemma:acsa-helperlemma2}
Let $x_t^{\mathrm{ag}}:=(1-\alpha_t)x_{t-1}^{\mathrm{ag}}+\alpha_t x_t$ for some $\alpha_t\in[0,1]$ and $(x_{t-1}^{\mathrm{ag}},x_t)\in X\times X$. Suppose that $f$ satisfies the upper-model bound in \eqref{eq:acsa-model}.
Then, for any $z\in X$ and any $f'(z)\in \partial f(z)$, we have
\[\Psi(x_t^{\mathrm{ag}})\leq
(1-\alpha_t)\Psi(x_{t-1}^{\mathrm{ag}})+\alpha_t\left[
f(z)+\langle f'(z),x_t-z\rangle+\mathcal X(x_t)\right]
+ \frac{L}{2}\|x_t^{\mathrm{ag}}-z\|_p^2+M_f\|x_t^{\mathrm{ag}}-z\|_p.\]
\end{lemma}
\begin{proof}
    The proof is the same as in \cite{ghadimi2012optimal}.
\end{proof}

Let $\mathcal F_{t-1}:=\sigma(\xi_1,\ldots,\xi_{t-1})$ denote
 the oracle history, with $\mathcal F_0$ trivial. Assume that
 $x_t^{md}$ is $\mathcal F_{t-1}$-measurable and that the stochastic
 oracle has a selected subgradient $g(x)\in\partial f(x)$ satisfying
\begin{equation}\label{eq:acsa-oracle}
\mathbb E[G(x_t^{md},\xi_t)\mid\mathcal F_{t-1}]
=g(x_t^{md}),\qquad
\mathbb E[\|G(x_t^{md},\xi_t)-g(x_t^{md})\|_q^q
\mid\mathcal F_{t-1}]\le\sigma^q.
\end{equation}
 Throughout the proof use $f'(x_t^{md})=g(x_t^{md})$ and write
\begin{equation}\label{eq:acsa-gradient-error}
G_t:=G(x_t^{md},\xi_t),\qquad
\delta_t:=G_t-g(x_t^{md}).
\end{equation}
The following proposition establishes a basic recursion for the generic AC-SA algorithm.

\begin{proposition}
\label{proposition:acsa-basic-recursion}
Let $(x_{t-1},x_{t-1}^{ag})\in X\times X$ be given. Also let $(x_t^{md},x_t,x_t^{ag})\in X\times X\times X$ be computed according to (3.1), (3.2) and (3.3) in \cite{ghadimi2012optimal}. If $K_t:=\nu(\mu+\gamma_t)-(L^{p-1}/M_f^{p-2})\alpha_t^p>0$, then for any $x\in X$, we have
\begin{align}
\Psi(x_t^{ag})+\mu V(x_t,x)
&\le (1-\alpha_t)\left[\Psi(x_{t-1}^{ag})+\mu V(x_{t-1},x)\right]
+\alpha_t l_\Psi(x_t^{md},x) \notag\\
&\quad +\gamma_t\left[V(x_{t-1},x)-V(x_t,x)\right]+\Delta_t(x),
\label{eq:acsa-basic-recursion}
\end{align}
where $l_\Psi(z,x)$ is defined in \eqref{eq:acsa-lower-model},
\begin{equation}\label{eq:acsa-delta}
\Delta_t(x):=\alpha_t\langle \delta_t,x-x_{t-1}^+\rangle+
\frac{\alpha_t^q(2M_f+\|\delta_t\|_q)^q}{qK_t^{q-1}},
\end{equation}
\begin{equation}
x_{t-1}^+:=\frac{\alpha_t\mu}{\mu+\gamma_t}x_t^{md}
+\frac{(1-\alpha_t)\mu+\gamma_t}{\mu+\gamma_t}x_{t-1}.
\end{equation}
\end{proposition}

\begin{proof}
    We first establish some basic relations among the search points $x_t^{ag}$, $x_t^{md}$ and $x_t$. Denote $d_t:=x_t^{ag}-x_t^{md}$. It follows from (3.1) and (3.3) in \cite{ghadimi2012optimal} that
\begin{align}
d_t
&=\alpha_t x_t+(1-\alpha_t)x_{t-1}^{ag}-x_t^{md} \notag\\
&=\alpha_t\left(x_t-\frac{\alpha_t\mu}{\mu+\gamma_t}x_t^{md}
-\frac{(1-\alpha_t)\mu+\gamma_t}{\mu+\gamma_t}x_{t-1}\right)
=\alpha_t(x_t-x_{t-1}^+),
\label{eq:acsa-search-point-relation}
\end{align}
which, in view of the convexity of $\|\cdot\|^p_p$ and the uniform convexity of $\omega$, implies that
\begin{align}
\frac{\nu(\mu+\gamma_t)}{p\alpha_t^p}\|d_t\|^p_p
&\le \frac{\nu(\mu+\gamma_t)}{p(\mu+\gamma_t)}
\left[\alpha_t\mu\|x_t-x_t^{md}\|^p_p+((1-\alpha_t)\mu+\gamma_t)\|x_t-x_{t-1}\|^p_p\right] \notag\\
&\le \alpha_t\mu V(x_t^{md},x_t)+((1-\alpha_t)\mu+\gamma_t)V(x_{t-1},x_t).
\label{eq:acsa-distance-control}
\end{align}
Using the above result and Lemma \ref{lemma:acsa-helperlemma2} (with $z=x_t^{md}$), we have
\begin{align}
\Psi(x_t^{ag})
&\le (1-\alpha_t)\Psi(x_{t-1}^{ag})
+\alpha_t\left[f(x_t^{md})+\langle f'(x_t^{md}),x_t-x_t^{md}\rangle+\mathcal X(x_t)\right] + \frac{L}{2}\|d_t\|_p^2 +M_f\|d_t\|_p \notag\\
&=(1-\alpha_t)\Psi(x_{t-1}^{ag})
+\alpha_t\left[f(x_t^{md})+\langle f'(x_t^{md}),x_t-x_t^{md}\rangle+\mathcal X(x_t)\right]+\frac{\nu(\mu+\gamma_t)}{p\alpha_t^p}\|d_t\|^p_p \notag\\
&\quad -\frac{\nu(\mu+\gamma_t)}{p\alpha_t^p}\|d_t\|^p_p+ \frac{L}{2}\|d_t\|_p^2 + M_f\|d_t\|_p \notag\\
&\leq (1-\alpha_t)\Psi(x_{t-1}^{ag})
+\alpha_t\left[f(x_t^{md})+\langle f'(x_t^{md}),x_t-x_t^{md}\rangle+\mathcal X(x_t)+\mu V(x_t^{md},x_t)\right] \notag\\
&\quad +((1-\alpha_t)\mu+\gamma_t)V(x_{t-1},x_t)
-\frac{\nu(\mu+\gamma_t)}{p\alpha_t^p}\|d_t\|^p_p+ \frac{L}{2}\|d_t\|_p^2 +M_f\|d_t\|_p.
\label{eq:acsa-objective-upper-bound}
\end{align}
By using Lemma \ref{lemma:acsa-helperlemma1} with $p(u)=\alpha_t[\langle G_t,u\rangle+\mathcal X(u)]$, $\mu_1=\alpha_t\mu$, $\mu_2=(1-\alpha_t)\mu+\gamma_t$, $\tilde x=x_t^{md}$ and $\tilde y=x_{t-1}$, we have
\begin{align}
&\alpha_t\left[f(x_t^{md})+\langle G_t,x_t-x_t^{md}\rangle+\mathcal X(x_t)+\mu V(x_t^{md},x_t)\right]
+((1-\alpha_t)\mu+\gamma_t)V(x_{t-1},x_t) \notag\\
&\le \alpha_t\left[f(x_t^{md})+\langle G_t,x-x_t^{md}\rangle+\mathcal X(x)+\mu V(x_t^{md},x)\right]
+((1-\alpha_t)\mu+\gamma_t)V(x_{t-1},x)-(\mu+\gamma_t)V(x_t,x) \notag\\
&=\alpha_t l_\Psi(x_t^{md},x)+\alpha_t\langle\delta_t,x-x_t^{md}\rangle
+((1-\alpha_t)\mu+\gamma_t)V(x_{t-1},x)-(\mu+\gamma_t)V(x_t,x),
\label{eq:acsa-proximal-comparison}
\end{align}
for any $x\in X$,where the last equality follows from
\eqref{eq:acsa-lower-model}
and \eqref{eq:acsa-gradient-error}.
Combining \eqref{eq:acsa-objective-upper-bound} and \eqref{eq:acsa-proximal-comparison}, and using the fact that $f'(x_t^{md})=G_t-\delta_t$ due to \eqref{eq:acsa-gradient-error}, we obtain
\begin{align}
\Psi(x_t^{ag})
&\le (1-\alpha_t)\left[\Psi(x_{t-1}^{ag})+\mu V(x_{t-1},x)\right]
+\alpha_t l_\Psi(x_t^{md},x)+\gamma_t\left[V(x_{t-1},x)-V(x_t,x)\right] \notag\\
&\quad -\mu V(x_t,x)
\underbrace{-\frac{\nu(\mu+\gamma_t)}{p\alpha_t^p}\|d_t\|^p_p
+ \frac{L}{2}\|d_t\|_p^2 +M_f\|d_t\|_p+\alpha_t\langle\delta_t,x-x_t\rangle}_{U_t},
\quad \forall x\in X.
\end{align}
It remains to show that the term $U_t$ defined above is bounded by
$\Delta_t(x)$ as defined in the proposition. Indeed, we have, by the last identity in \eqref{eq:acsa-search-point-relation},
\begin{align}
U_t
&=-\frac{\nu(\mu+\gamma_t)}{p\alpha_t^p}\|d_t\|^p_p + \frac{L}{2}\|d_t\|_p^2
+M_f\|d_t\|_p-\langle\delta_t,d_t\rangle
+\langle\delta_t,d_t+\alpha_t(x-x_t)\rangle \notag\\
&\leq -\frac{\nu(\mu+\gamma_t)}{p\alpha_t^p}\|d_t\|^p_p + \frac{L^{p-1}}{pM_f^{p-2}}\|d_t\|_p^p
+2M_f\|d_t\|_p-\langle\delta_t,d_t\rangle
+\alpha_t\langle\delta_t,x-x_{t-1}^+\rangle \notag\\
&\leq -\frac{\nu(\mu+\gamma_t)-(L^{p-1}/M_f^{p-2})\alpha_t^p}{p\alpha_t^p}\|d_t\|^p_p +(2M_f+\|\delta_t\|_q)\|d_t\|_p+\alpha_t\langle\delta_t,x-x_{t-1}^+\rangle \notag\\
&\leq \frac{\alpha_t^q(2M_f+\|\delta_t\|_q)^q}{qK_t^{q-1}}
+\alpha_t\langle\delta_t,x-x_{t-1}^+\rangle
=\Delta_t(x),
\label{eq:acsa-remainder-bound}
\end{align}
where the last inequality follows by the Young's inequality w.r.t. $\|d_t\|$.
\end{proof}

Lemma \ref{lemma:next_acsa_lemma} below follows from Proposition~\ref{proposition:acsa-basic-recursion} by taking summation of the relations in \eqref{eq:acsa-basic-recursion}.

\begin{lemma}
\label{lemma:next_acsa_lemma}
Let $\{x_t^{ag}\}_{t\geq 1}$ be computed by the generic AC-SA algorithm for solving
$F_{X,X}(L,M_f,\sigma,\mu)$. Also assume that $\{\alpha_t\}_{t\geq 1}$ and
$\{\gamma_t\}_{t\geq 1}$ are chosen so that $K_t>0$ for all $t\ge1$. We have
\begin{flalign}
\Psi(x_t^{ag})+\mu V(x_t,x)
-\Gamma_t\sum_{\tau=1}^t
\left[
\frac{\alpha_\tau}{\Gamma_\tau}
l_\Psi(x_\tau^{md},x)
\right]
&\leq
\Gamma_t\sum_{\tau=1}^t
\frac{\gamma_\tau}{\Gamma_\tau}
\left[
V(x_{\tau-1},x)-V(x_\tau,x)
\right]
+
\Gamma_t\sum_{\tau=1}^t
\frac{\Delta_\tau(x)}{\Gamma_\tau},
\label{eq:acsa-cumulative-recursion}
\end{flalign}
for any $x\in X$ and any $t\geq 1$, where $l_\Psi(z,x)$ and $\Delta_\tau(x)$ are defined in
\eqref{eq:acsa-lower-model} and \eqref{eq:acsa-delta}, respectively.
\end{lemma}

\begin{proof}
Dividing both sides of relation \eqref{eq:acsa-basic-recursion} by $\Gamma_t$, and using the definition of $\Gamma_t$
in \eqref{eq:acsa-gamma-sequence} and the fact that $\alpha_1=1$, we have
\begin{flalign*}
\frac{1}{\Gamma_t}
\left[
\Psi(x_t^{ag})+\mu V(x_t,x)
\right]
&\leq
\frac{1}{\Gamma_{t-1}}
\left[
\Psi(x_{t-1}^{ag})+\mu V(x_{t-1},x)
\right]
+
\frac{\alpha_t}{\Gamma_t}
l_\Psi(x_t^{md},x)
\\
&\quad+
\frac{\gamma_t}{\Gamma_t}
\left[
V(x_{t-1},x)-V(x_t,x)
\right]
+
\frac{\Delta_t(x)}{\Gamma_t},
\quad \forall t\geq 2,
\end{flalign*}
and
\begin{flalign*}
\frac{1}{\Gamma_1}
\left[
\Psi(x_1^{ag})+\mu V(x_1,x)
\right]
&\leq
\frac{\alpha_1}{\Gamma_1}
l_\Psi(x_1^{md},x)
+
\frac{\gamma_1}{\Gamma_1}
\left[
V(x_0,x)-V(x_1,x)
\right]
+
\frac{\Delta_1(x)}{\Gamma_1}.
\end{flalign*}
Summing up the above inequalities, we have
\begin{flalign*}
\frac{1}{\Gamma_t}
\left[
\Psi(x_t^{ag})+\mu V(x_t,x)
\right]
&\leq
\sum_{\tau=1}^t
\frac{\alpha_\tau}{\Gamma_\tau}
l_\Psi(x_\tau^{md},x)
+
\sum_{\tau=1}^t
\frac{\gamma_\tau}{\Gamma_\tau}
\left[
V(x_{\tau-1},x)-V(x_\tau,x)
\right]
+
\sum_{\tau=1}^t
\frac{\Delta_\tau(x)}{\Gamma_\tau},
\end{flalign*}
which clearly implies \eqref{eq:acsa-cumulative-recursion}.
\end{proof}

\begin{theorem}
Consider the generic AC-SA algorithm for $F_{X,X}(L,M_f,\sigma,\mu)$.
Assume $\alpha_1=1$, $0<\alpha_t<1$ for $t\ge2$, and that
$\{\alpha_t\}_{t\ge1}$ and $\{\gamma_t\}_{t\ge1}$ are chosen such that
\begin{flalign}
&K_t:=\nu(\mu+\gamma_t)
-\frac{L^{p-1}}{M_f^{p-2}}\alpha_t^p>0,
\label{eq:acsa-parameter-positive}
\end{flalign}
\begin{flalign}
&\gamma_1/\Gamma_1=\gamma_2/\Gamma_2=\ldots,
\label{eq:acsa-parameter-ratio}
\end{flalign}
where
\begin{flalign}
\Gamma_t:=
\begin{cases}
1, & t=1,\cr
(1-\alpha_t)\Gamma_{t-1}, & t\geq 2.
\end{cases}
\label{eq:acsa-gamma-sequence}
\end{flalign}
Then,

\noindent
a) under the stochastic-oracle conditions \eqref{eq:acsa-oracle}, we have
\begin{flalign}
\E[\Psi(x_t^{ag})-\Psi^*]
&\leq B_e(t):=
\Gamma_t\gamma_1V(x_0,x^*)+\Gamma_t\sum_{\tau=1}^t\frac{2^{q-1}
\alpha_\tau^q[(2M_f)^q+\sigma^q]}{q\Gamma_\tau K_\tau^{q-1}}
\label{eq:acsa-expected-error}
\end{flalign}
for any $t\geq 1$, where $x^*$ is an arbitrary optimal solution of \eqref{acsa-opt-problem}.
\end{theorem}

\begin{proof}
We show this following the original proof of Theorem~1 in \cite{ghadimi2012optimal}.
Observe that by the definition of $\Gamma_t$ and the fact that
$\alpha_1=1$, we have
\begin{align}
\sum_{\tau=1}^{t}\frac{\alpha_\tau}{\Gamma_\tau}
&=\frac{\alpha_1}{\Gamma_1}
+\sum_{\tau=2}^{t}\frac{1}{\Gamma_\tau}\left(1-\frac{\Gamma_\tau}{\Gamma_{\tau-1}}\right)=\frac{1}{\Gamma_1}
+\sum_{\tau=2}^{t}\left(\frac{1}{\Gamma_\tau}-\frac{1}{\Gamma_{\tau-1}}\right)
=\frac{1}{\Gamma_t}.
\label{eq:acsa-weight-identity}
\end{align}
The lower-model inequality follows directly from
\eqref{eq:acsa-model}: for every $x,z\in X$,
\[l_\Psi(z,x)=f(z)+\langle f'(z),x-z\rangle+\mu V(z,x)+\mathcal X(x)
\le f(x)+\mathcal X(x)=\Psi(x).\]
Using this inequality and the previous observation, we obtain
\begin{align}
\Gamma_t\sum_{\tau=1}^{t}
 \left[\frac{\alpha_\tau}{\Gamma_\tau}
 l_\Psi(x_\tau^{md},x)\right]
&\le\Gamma_t\sum_{\tau=1}^{t}
 \left[\frac{\alpha_\tau}{\Gamma_\tau}\Psi(x)\right]
=\Psi(x),\qquad\forall x\in X.
\label{eq:acsa-lower-model-average}
\end{align}
Moreover, by the condition
$\gamma_1/\Gamma_1=\gamma_2/\Gamma_2=\cdots$, we have
\begin{align}
&\Gamma_t\sum_{\tau=1}^{t}\frac{\gamma_\tau}{\Gamma_\tau}
 \left[V(x_{\tau-1},x)-V(x_\tau,x)\right]
=\Gamma_t\frac{\gamma_1}{\Gamma_1}
 \left[V(x_0,x)-V(x_t,x)\right]
 \le\Gamma_t\gamma_1V(x_0,x),
\label{eq:acsa-bregman-telescoping}
\end{align}
where the last inequality follows from $\Gamma_1=1$ and
$V(x_t,x)\ge0$.
Substituting \eqref{eq:acsa-lower-model-average} and \eqref{eq:acsa-bregman-telescoping} into the
cumulative recursion \eqref{eq:acsa-cumulative-recursion}, and using $\mu V(x_t,x)\ge0$, gives
\begin{align}
\Psi(x_t^{ag})-\Psi(x)
&\le \Psi(x_t^{ag})+\mu V(x_t,x)-\Psi(x)\le\Gamma_t\gamma_1V(x_0,x)+\Gamma_t\sum_{\tau=1}^{t}
\frac{\Delta_\tau(x)}{\Gamma_\tau},\qquad\forall x\in X.
\label{eq:acsa-objective-recursion}
\end{align}
Note that $(x_\tau^{md},x_{\tau-1},x_{\tau-1}^{ag})$
is determined by the history
$\mathcal F_{\tau-1}=\sigma(\xi_1,\ldots,\xi_{\tau-1})$.
Consequently, $x_{\tau-1}^+$ is
$\mathcal F_{\tau-1}$-measurable.
Taking $x=x^\star$ and using
conditional unbiasedness of the stochastic oracle, we obtain
\begin{align}
\mathbb E\!\left[
\left\langle\delta_\tau,x^\star-x_{\tau-1}^+\right\rangle
 \middle|\mathcal F_{\tau-1}\right]
&=\left\langle
 \mathbb E[\delta_\tau\mid\mathcal F_{\tau-1}],
 x^\star-x_{\tau-1}^+\right\rangle=0.
\label{eq:acsa-noise-mean-zero}
\end{align}
By the conditional $q$-moment assumption in \eqref{eq:acsa-oracle},
\[\mathbb E[\|\delta_\tau\|_q^q\mid\mathcal F_{\tau-1}]\le\sigma^q.\]
Furthermore, the elementary inequality
$(a+b)^q\le2^{q-1}(a^q+b^q)$, valid for $a,b\ge0$,
implies
\begin{align*}
\mathbb E\!\left[(2M_f+\|\delta_\tau\|_q)^q\right] &\le2^{q-1}\!\left((2M_f)^q+\mathbb E\|\delta_\tau\|_q^q\right)\\
&\le2^{q-1}
\!\left((2M_f)^q+\sigma^q\right).
\end{align*}
Recalling the definition of $\Delta_\tau(x)$ in the preceding
proposition and observing that $q-1=1/(p-1)$, we conclude that
\begin{align*}
\mathbb E\!\left[\Delta_\tau(x^\star)\right]
&=\alpha_\tau\mathbb E\!\left[
\left\langle\delta_\tau,x^\star-x_{\tau-1}^+\right\rangle
\right]+\frac{\alpha_\tau^q}
{qK_\tau^{q-1}}\mathbb E\!\left[(2M_f+\|\delta_\tau\|_q)^q\right]\\
&\le\frac{2^{q-1}\alpha_\tau^q
\bigl((2M_f)^q+\sigma^q\bigr)}
{qK_\tau^{q-1}}.
\end{align*}
Finally, taking expectations on both sides of \eqref{eq:acsa-objective-recursion} with
$x=x^\star$, and applying the previous estimate, we have
\begin{align*}
\mathbb E\!\left[\Psi(x_t^{ag})-\Psi^\star\right]
&\le\Gamma_t\gamma_1V(x_0,x^\star)+\Gamma_t\sum_{\tau=1}^{t}
\frac{\mathbb E[\Delta_\tau(x^\star)]}{\Gamma_\tau}\\
&\le\Gamma_t\gamma_1V(x_0,x^\star)+\frac{2^{q-1}\Gamma_t}{q}
\sum_{\tau=1}^{t}
\frac{\alpha_\tau^q\bigl((2M_f)^q+\sigma^q\bigr)}
{\Gamma_\tau K_\tau^{q-1}},
\end{align*}
as claimed.
\end{proof}

\begin{theorem}
\label{theorem:p-parallel-acsa}
Assume that the conditions of
Propositions~\ref{proposition:prop2_equivalent} and
\ref{proposition:p-acsa-broo} hold,
then \textsc{BallAccel}
(Proposition~\ref{proposition:prop2_equivalent}),
using parallel $\ell_p$-AC-SA
(Algorithm~\ref{alg:ac-sa}) as a ball optimization oracle,
solves Problem~\ref{problem:parallel_problem}
with expected error $\epsilon_{\rm opt}$, with
\[
 O\!\left(
 d^{\frac{p}{2(p+1)}}\kappa^{\frac{p}{p+1}}
 \log^4(d\kappa)
 \right)
 \quad\text{query depth}
\]
and
\[O\!\left(\left(d^{\frac{p}{2(p+1)}}\kappa^{\frac{p}{p+1}}
+d^{\frac{p}{2(p+2)}}
\kappa^{\frac{2p}{p+2}}
\right)\log^{2p}(d\kappa)
+d^{\frac{p(p-2)}{2(p+1)}}
\kappa^{\frac{p(2p-1)}{p+1}}
\log^{2p(p-1)+2}(d\kappa)
\right)\quad\text{total queries},\]
and an additional computational cost of
\[O\!\left(\left(d^{\frac{p}{2(p+1)}}\kappa^{\frac{p}{p+1}}
+d^{\frac{p}{2(p+2)}}
\kappa^{\frac{2p}{p+2}}
\right)\log^{2p}(d\kappa)
\right)\quad\text{depth}\]
and 
\[O\!\left(\left[\left(
d^{\frac{p}{2(p+1)}}
\kappa^{\frac{p}{p+1}}
+d^{\frac{p}{2(p+2)}}
\kappa^{\frac{2p}{p+2}}
\right)\log^{2p}(d\kappa)
+d^{\frac{p(p-2)}{2(p+1)}}
\kappa^{\frac{p(2p-1)}{p+1}}
\log^{2p(p-1)+2}(d\kappa)
\right]d\right)\quad\text{work}
\]
The computational bounds use the vector-operation cost convention of Problem~\ref{problem:parallel_problem} and
Proposition~\ref{proposition:p-acsa-broo}.
\end{theorem}

\begin{proof}
Let
$x^\star\in\argmin_{x\in\mathbb B_p(R)}f(x)$,
and define
\[
f_{\rm pen}(x)=f(x)+2L(\|x\|_p-R)_+,
\qquad
F=\widehat{f_{\rm pen}}_\rho,
\qquad
r=\rho=\frac{\epsilon_{\rm opt}}{C_p'L\sqrt d}.
\]
By Lemma~\ref{lem:radial-penalty},
$f_{\rm pen}$ is convex and $3L$-Lipschitz,
and its stochastic subgradient oracle has dual
$q$-th moment at most $(4L)^q$.
Convexity and smoothing give
$F(x)\ge f_{\rm pen}(x)$.
Choosing $C_p'$ sufficiently large, the smoothing estimate yields
\[0\le F(x)-f_{\rm pen}(x)
\le 3L\mathbb E\|\xi_\rho\|_p
\le\frac{\epsilon_{\rm opt}}4,
\qquad x\in\mathbb R^d.\]
By the Lipschitz continuity of $f$ and projection onto
$\mathbb B_p(R)$,
\[
f_{\rm pen}(x)
\ge f(x^\star)+L(\|x\|_p-R)_+,
\qquad x\in\mathbb R^d.
\]
Consequently, $F\ge f_{\rm pen}$ is coercive and attains
its global minimum.
For any global minimizer $z$ of $F$,
\[
L(\|z\|_p-R)_+
\le F(z)-f(x^\star)
\le F(x^\star)-f(x^\star)
\le\frac{\epsilon_{\rm opt}}4.
\]
Since $\epsilon_{\rm opt}\le LR$, this implies
\[
\|z\|_p
\le R+\frac{\epsilon_{\rm opt}}{4L}
\le2R.
\]
Thus $F$ has a global minimizer in $\mathbb B_p(2R)$. Apply the implementation in
Proposition~\ref{proposition:prop2_equivalent}
to $F$, with comparison point $x^\star$,
comparison radius $2R$, Lipschitz parameter $3L$,
and target accuracy $\epsilon_{\rm opt}/2$.
All ball optimization oracles are applied directly to $F$.
Each required regularized subproblem therefore has objective
\[
x\mapsto F(x)+\frac{\lambda}{p}\|x-\bar x\|_p^p.
\]
Proposition~\ref{proposition:p-acsa-broo}, applied to
$f_{\rm pen}$ and its stochastic subgradient oracle,
implements these calls using
Algorithm~\ref{alg:ac-sa}, with $L$ replaced by $3L$. The BallAccel parameters corresponding to this invocation
differ from the original $\kappa$, $K$, and $\lambda^\star$
only by constant factors depending on $p$.
In the complexity estimates below, we retain their definitions
from Proposition~\ref{proposition:init_prop_guarantee}.
In particular,
\[K=\Theta\left(
d^{\frac{p}{2(p+1)}}
\kappa^{\frac p{p+1}}\right),
\qquad\log(R\kappa/r)=O(\log(d\kappa)).\]
Each $\ell_p$-AC-SA ball optimization call has constant query depth.
Thus Proposition~\ref{proposition:prop2_equivalent}
gives total query depth
\[O\left(
K\log\kappa\log^3(R\kappa/r)
\right)=O\left(
d^{\frac{p}{2(p+1)}}
\kappa^{\frac p{p+1}}
\log^4(d\kappa)\right).\]
For the computational depth, all calls use
$\lambda=\Omega_p(\lambda^\star)$, and hence
\[\frac{L}{\lambda r^{p-1}}
\le\frac{C_p\kappa}{K^{1/p}\log^p\kappa}\le
\frac{C_p\kappa}{K^{1/p}}.\]
By Proposition~\ref{proposition:p-acsa-broo},
each call in the first requirement of
Proposition~\ref{proposition:prop2_equivalent}
requires
\[
O\left(
1+\left(\frac{\kappa}{K^{1/p}}\right)^{\frac p{p+2}}
\right)
\]
sequential updates.
A level-$j$ call in the second requirement requires
\[
O\left(
1+j+\log\log(R\kappa/r)
+\left(1+\frac{\kappa}{K^{1/p}}\right)^{\frac p{p+2}}
2^{\frac{j(p-2)}{p+2}}
\log^{\frac{2p(p-2)}{p+2}}(R\kappa/r)
\right)
\]
updates.
Using the call counts in the two requirements and
\[
\sum_{j\ge1}2^{-j}(1+j)=O(1),
\qquad
\sum_{j\ge1}2^{-\frac{4j}{p+2}}=O(1),
\]
the combined computational depth is
\[
O\left(
\left[
K+\kappa^{\frac p{p+2}}K^{\frac{p+1}{p+2}}
\right]\log^{2p}(d\kappa)
\right).
\]
Here we used
\[
4\le2p,
\qquad
2+\frac{2p(p-2)}{p+2}\le2p.
\]
Substituting the expression for $K$ gives the stated
computational-depth bound. The contribution of the
$L^p/(\lambda\phi^{p-1})$ term in
Proposition~\ref{proposition:p-acsa-broo}
is summed exactly as in the proof of
Theorem~\ref{theorem:epochda_parallel_guarantee},
and is at most
\[
O\left(
\kappa^pK^{p-2}
\log^{2p(p-1)+2}(d\kappa)
\right).
\]
Adding the update bound above gives total queries
\[
O\left(
\left[
K+\kappa^{\frac p{p+2}}K^{\frac{p+1}{p+2}}
\right]\log^{2p}(d\kappa)
+
\kappa^pK^{p-2}\log^{2p(p-1)+2}(d\kappa)
\right).
\]
The identities
\[
\kappa^{\frac p{p+2}}K^{\frac{p+1}{p+2}}
=
\Theta\left(
d^{\frac{p}{2(p+2)}}
\kappa^{\frac{2p}{p+2}}
\right)
\]
and
\[
\kappa^pK^{p-2}
=
\Theta\left(
d^{\frac{p(p-2)}{2(p+1)}}
\kappa^{\frac{p(2p-1)}{p+1}}
\right)
\]
give the stated query bound.
Under the stated computational cost convention,
multiplication by $d$ gives the work bound. Finally, the unpenalized implementation of
Proposition~\ref{proposition:prop2_equivalent}
returns $\widehat x$ satisfying
\[
\mathbb E F(\widehat x)
\le F(x^\star)+\frac{\epsilon_{\rm opt}}2.
\]
Since $f\le f_{\rm pen}\le F$ and
$f_{\rm pen}(x^\star)=f(x^\star)$,
\[
\begin{aligned}
\mathbb E f(\widehat x)
\le\mathbb E F(\widehat x)\le F(x^\star)+\frac{\epsilon_{\rm opt}}2\le f(x^\star)+\frac{3\epsilon_{\rm opt}}4\le f(x^\star)+\epsilon_{\rm opt}.
\end{aligned}
\]
This completes the proof.
\end{proof}

\section{$\ell_p$ Parallel Stochastic Convex Optimization: Lower Bounds}  
\label{section:parallel_lower_bounds}
We consider minimizing a convex function $f:\R^d\to\R$ that is
$1$-Lipschitz with respect to $\|\cdot\|_p$. To fix notation, let
\[
B_p^d=\{x\in\R^d:\|x\|_p\le1\},
\qquad
B_q^d=\{s\in\R^d:\|s\|_q\le1\}.
\]
Then $f$ is $1$-Lipschitz with respect to $\|\cdot\|_p$ if and only if all
its subgradients lie in $B_q^d$. An algorithm accesses $f$ only through a
\emph{local oracle}, meaning that the answer at a point $x$ depends only on the restriction
of $f$ to an arbitrarily small neighborhood of $x$. This class contains the
usual (sub)gradient oracle as well as its noisy variants, and our lower bound
holds for any fixed local oracle.

We consider \emph{$K$-parallel} algorithms, whereby in each round the algorithm makes up to $K$ queries in parallel,
chosen based on the answers from previous rounds but not on the answers to
other queries in the same round. Following
\cite{balkanski2018parallelization}, we measure complexity through
\emph{adaptivity}, i.e., the number of sequential rounds $M$. Throughout
the section, we focus on $K\le\mathrm{poly}(d)$. Otherwise, as noted by
\cite{diakonikolas2019lower}, with exponentially many queries per round, one
can query an $\eps$-net of the feasible set and return the best point found.

Our hard instance is constructed on a restricted region of the unit ball
$B_p^d$ in which every coordinate is of order $d^{-1/p}$. Specifically, we
define the feasible region
\[
X:=\left\{x\in\R^d: \|x\|_p\le1,\quad \|x\|_\infty\le4d^{-1/p}\right\},
\]
and require that all queries and the final output of the algorithm
belong to $X$.

\begin{theorem}\label{thm:main}
There are universal constants $c,C>0$ with the following property.
Let $M\ge2$, $K\ge1$, $\rho\in(0,1/4)$, and suppose $d\ge16M$.
Further, assume that
\begin{align}
M^{-1/p} &\ge C\eps,
\label{eq:main-accuracy}\\
M^{\frac32-\frac1p}
\sqrt{\frac{\log(16KM^2/\rho)}{d}} &\le c.
\label{eq:main-variance}
\end{align}
Then for every randomized $K$-parallel algorithm whose queries and
output belong to $X$, and every fixed local oracle, there exists a
convex function $f:\R^d\to\R$, $1$-Lipschitz with respect to $\ell_p$,
such that the algorithm returns an $\eps$-minimizer of $f$ over $X$ in
fewer than $M$ rounds with probability at most $\rho$. Consequently,
up to logarithmic factors,
\[M=\widetilde\Omega\left(
\min\left\{\eps^{-p},d^{\frac{p}{3p-2}}\right\}\right).\]
\end{theorem}
\subsection{Random block directions}
\label{sec:block-directions}

The directions below follow the disjoint-support construction of
\cite{diakonikolas2019lower} for $p\ge2$, with two key differences: the supports
are drawn at random and each $z_i$ comes with a
direction $h_i$ that is used to construct the wall function.
For the rest of the section let $L=\left\lfloor\frac{d}{8M}\right\rfloor$ and
choose pairwise disjoint sets $J_1,\ldots,J_M\subset[d]$, each of size $L$,
uniformly among all ordered disjoint collections. Independently choose
Rademacher signs $r_{ij}\in\{-1,+1\}$ and set
\[
z_i=L^{-1/q}\sum_{j\in J_i}r_{ij}e_j,\qquad
h_i=L^{-1/p}\sum_{j\in J_i}r_{ij}e_j,\qquad
x^*=-M^{-1/p}\sum_{i=1}^M h_i.
\]
Since the blocks are disjoint, $\|z_i\|_q=\|h_i\|_p=\|x^*\|_p=1$ and
$\langle z_i,x^*\rangle=-M^{-1/p}$ for every $i\in[M]$. Moreover, $x^*\in X$:
$d\ge16M$ implies $L\ge\frac{d}{16M}$, and therefore
\[
\|x^*\|_\infty=(ML)^{-1/p}\le\left(\frac{16}{d}\right)^{1/p}\le4d^{-1/p}.
\]

\subsection{The wall function}
\label{sec:wall}
As in \cite{bubeck2019complexity}, who introduced a wall in the $\ell_2$ setting, our objective also contains a wall $W$. For $t\in[M]$, let
$R_t:=[d]\setminus\bigcup_{i<t}J_i$ denote the coordinates outside the first
$t-1$ blocks. The wall is a convex function that is small at $x^*$ and, at round
$t$, dominates all future affine pieces at every query $x\in X$ with
$\|x_{R_t}\|_p>4M^{-1/p}$ (Lemmas~\ref{lem:shielding} and~\ref{lem:hiding}).
At such queries the objective carries no information about the future
directions, so concentration is needed only for the remaining queries. Fix a parameter $\theta>0$, to be chosen in Lemma~\ref{lem:all-hidden}. For
$t\in[M]$ let
\begin{equation}\label{eq:Ct}
\mathcal C_t := \left\{ s\in B_q^d:
\|s\|_\infty\le L^{-1/q},
\quad |\langle h_i,s\rangle|\le\theta
\ \text{for all }i<t\right\},
\end{equation}
and let
\begin{equation}\label{eq:Cfull}
\mathcal C:=\left\{s\in B_q^d:
\|s\|_\infty\le L^{-1/q},
\quad|\langle h_i,s\rangle|\le\theta
\ \text{for all }i\le M\right\}.
\end{equation}
The caps $|\langle h_i,s\rangle|\le\theta$ keep the slope $s$ nearly orthogonal
to the block directions, which makes the wall small at $x^*$.
We define the \emph{truncated wall}
\begin{equation}\label{eq:Wt}
W_t(x)
:=
\max_{s\in\mathcal C_t}
\left(
\ip{x}{s}
-\frac{M^{-1/p}}{32}\norm[2]{s}^2
\right)
-\frac{33}{32}M^{-1/p},
\end{equation}
which depends only on $h_1,\ldots,h_{t-1}$, and the \emph{full wall}
\begin{equation}\label{eq:W}
W(x)
:=
\max_{s\in\mathcal C}
\left(
\ip{x}{s}
-\frac{M^{-1/p}}{32}\norm[2]{s}^2
\right)
-\frac{33}{32}M^{-1/p}.
\end{equation}
The quadratic term is used only to make the maximizing slope unique and
continuous, so that $W_t=W$ on a neighborhood of each query
(Lemma~\ref{lem:wall-consistency}).

Now, we establish the two properties of the wall that we need: it is
large on queries with large residual norm on $R_t$, and small at $x^*$.
\begin{lemma}\label{lem:shielding}
Let $R_t:=[d]\setminus\bigcup_{i<t}J_i$.
For every $x\in X$ satisfying
$\|x_{R_t}\|_p>4M^{-1/p}$,
we have
\begin{equation}\label{eq:shield}
W_t(x) \ge \|x_{R_t}\|_p -\frac{17}{16}M^{-1/p}.
\end{equation}
Moreover, if $M\theta\le1$, then
\[W(x^*)\le -\frac1{32}M^{-1/p}.\]
\end{lemma}

\begin{proof}
Fix $x\in X$ with $\|x_{R_t}\|_p>4M^{-1/p}$ and define the
$\ell_p$-norming functional of $x_{R_t}$ by
\[s_j:=\frac{\operatorname{sgn}(x_j)|x_j|^{p-1}}
{\|x_{R_t}\|_p^{p-1}}
\mathbf 1_{\{j\in R_t\}}.\]
Then $\|s\|_q=1,\ \langle x,s\rangle
=\|x_{R_t}\|_p$.
We claim that $s\in\mathcal C_t$. Since $x\in X$, $\|x\|_\infty\le4d^{-1/p}$.
Hence
\begin{align*}
\|s\|_\infty= \left(\frac{\|x_{R_t}\|_\infty}{\|x_{R_t}\|_p}\right)^{p-1}
<\left(\frac{4d^{-1/p}}{4M^{-1/p}}
\right)^{p-1}=\left(\frac{M}{d}\right)^{1/q}.
\end{align*}
Since $L\le d/(8M)$,
\[
L^{-1/q}
\ge
\left(\frac{8M}{d}\right)^{1/q}
>
\left(\frac{M}{d}\right)^{1/q},
\]
so $\|s\|_\infty\le L^{-1/q}$. Moreover, $s$ is supported on
$R_t$, whereas every $h_i$ with $i<t$ is supported on $J_i$.
Therefore, $\langle h_i,s\rangle=0$ for $i<t$. Thus $s\in\mathcal C_t$. Since $q\le2$ and $\|s\|_q=1$, $\|s\|_2\le 1$.
Using $s$ as a feasible candidate in \eqref{eq:Wt},
\begin{align*}
W_t(x)\ge
\langle x,s\rangle
-\frac{M^{-1/p}}{32}\|s\|_2^2
-\frac{33}{32}M^{-1/p}\ge
\|x_{R_t}\|_p
-\frac{34}{32}M^{-1/p}=
\|x_{R_t}\|_p
-\frac{17}{16}M^{-1/p},
\end{align*}
which proves \eqref{eq:shield}.
We next bound the full wall at $x^*$. Let $s\in\mathcal C$. Since
\[x^*=-M^{-1/p}\sum_{i=1}^M h_i,\]
the cap constraints give
\begin{align*}
\langle x^*,s\rangle
=-M^{-1/p}\sum_{i=1}^M\langle h_i,s\rangle\le
M^{-1/p}\sum_{i=1}^M|\langle h_i,s\rangle|\le
M^{1-1/p}\theta.
\end{align*}
If $M\theta\le1$, then
$M^{1-1/p}\theta=M^{-1/p}(M\theta)
\le M^{-1/p}.$
The quadratic term is nonpositive, so for all $s \in \mathcal C$, it satisfies
\[\langle x^*,s\rangle
-\frac{M^{-1/p}}{32}\|s\|_2^2
\le M^{-1/p}.\]
Taking the maximum over $s\in\mathcal C$ and subtracting
$33M^{-1/p}/32$ yields
\[W(x^*)\le-\frac1{32}M^{-1/p}.\]
\end{proof}

\begin{lemma}\label{lem:wall-consistency}
For every $t\in[M]$, the truncated wall $W_t$ and the full wall $W$
are convex and $1$-Lipschitz with respect to $\ell_p$. Let
\[s_t(x):=\arg\max_{s\in\mathcal C_t}\left\{\langle x,s\rangle
-\frac{M^{-1/p}}{32}\|s\|_2^2
\right\}.\]
The maximizer is unique. Moreover, if
\begin{equation}\label{eq:cap-slack}
|\langle h_i,s_t(x)\rangle|
\le
\frac{\theta}{2}
\qquad
\text{for every }i=t,\ldots,M,
\end{equation}
then there exists an open neighborhood $U$ of $x$ such that
\[W_t(y)=W(y) \qquad \text{for every }y\in U.\]
\end{lemma}

\begin{proof}
For every fixed $s\in\mathcal C_t$, the map
\[x\mapsto \langle x,s\rangle
-\frac{M^{-1/p}}{32}\|s\|_2^2\]
is affine with slope $s\in B_q^d$. Hence the maximum over
$s\in\mathcal C_t$ is convex and $1$-Lipschitz with respect to
$\ell_p$. The same argument applies to $W$. Subtracting a constant would not affect both properties. Completing the square gives $s_t(x)=\Pi_{\mathcal C_t}(16M^{1/p}x)$, so the
maximizer is unique and, as the Euclidean projection is $1$-Lipschitz,
\begin{equation}\label{eq:st-continuity}
\|s_t(y)-s_t(x)\|_2 \le
\bigl\|16M^{1/p}y-16M^{1/p}
x\bigr\|_2= 16M^{1/p}\|y-x\|_2.
\end{equation}
Suppose \eqref{eq:cap-slack} holds. For every $i=t,\ldots,M$,
\begin{align*}
|\langle h_i,s_t(y)\rangle|
&\le
|\langle h_i,s_t(x)\rangle|
+
|\langle h_i,s_t(y)-s_t(x)\rangle|\\
&\le
\frac{\theta}{2}
+
\|h_i\|_2\|s_t(y)-s_t(x)\|_2\\
&\le
\frac{\theta}{2}
+
16M^{1/p}\|h_i\|_2\|y-x\|_2.
\end{align*}
Since there are finitely many indices, there exists $\delta>0$ such
that, whenever $\|y-x\|_2<\delta$,
\[
|\langle h_i,s_t(y)\rangle|<\theta
\qquad
\text{for every }i=t,\ldots,M.
\]
The vector $s_t(y)$ already satisfies the caps corresponding to
$i<t$ because $s_t(y)\in\mathcal C_t$. Thus $s_t(y)\in\mathcal C$.
Since $\mathcal C\subseteq\mathcal C_t$, the full and truncated
optimization problems have the same optimal value at such $y$.
Therefore $W_t(y)=W(y)$
throughout a neighborhood of $x$.
\end{proof}

\subsection{Hard objective and hiding future pieces}

For $\tau>0$, define
\begin{align}
f(x)&:=\max\left\{
\max_{1\le i\le M}\left(
\frac1{32}\langle z_i,x\rangle
-\frac{i\tau}{8}\right),W(x)
\right\},\label{eq:f}\\
f_t(x)&:=\max\left\{
\max_{1\le i\le t}\left(
\frac1{32}\langle z_i,x\rangle
-\frac{i\tau}{8}\right),W_t(x)
\right\}.
\label{eq:ft}
\end{align}

\begin{lemma}\label{lem:hiding}
Fix a round $t<M$ and a query $x\in X$. Assume
\eqref{eq:cap-slack}. If $\|x_{R_t}\|_p\le4M^{-1/p}$,
also assume
\begin{equation}\label{eq:hard-slack} |\langle z_i,x\rangle|\le\tau
\qquad \text{for every }i\ge t.
\end{equation}
Then $f=f_t$ on a neighborhood of $x$.
\end{lemma}

\begin{proof}
Lemma~\ref{lem:wall-consistency} gives $W=W_t$ on a neighborhood of
$x$. First suppose
$\|x_{R_t}\|_p>4M^{-1/p}.$
Every future $z_i$, $i\ge t$, is supported on $R_t$, and
$\|z_i\|_q=1$, so
\[\langle z_i,x\rangle=\langle z_i,x_{R_t}\rangle\le\|x_{R_t}\|_p.
\]
By Lemma~\ref{lem:shielding},
\[W_t(x)\ge\|x_{R_t}\|_p
-\frac{17}{16}M^{-1/p}.\]
Therefore, for every $i\ge t$,
\begin{align*}
W_t(x) -\frac1{32}\langle z_i,x\rangle
+\frac{i\tau}{8} &\ge
\frac{31}{32}\|x_{R_t}\|_p
-\frac{17}{16}M^{-1/p}>
\frac{45}{16}M^{-1/p}.
\end{align*}
Thus the wall strictly dominates every future affine piece, and the
strict inequality holds locally.
Now suppose $\|x_{R_t}\|_p\le4M^{-1/p}$.
For every $i>t$, \eqref{eq:hard-slack} gives
\begin{align*}
\left(\frac1{32}\langle z_t,x\rangle
-\frac{t\tau}{8}\right)-\left(
\frac1{32}\langle z_i,x\rangle
-\frac{i\tau}{8}\right)\ge
-\frac{\tau}{32}-\frac{t\tau}{8}
-\frac{\tau}{32}+\frac{(t+1)\tau}{8}
=\frac{\tau}{16}>0.
\end{align*}
Hence every future affine piece lies strictly below the $t$-th
piece. Combining the affine and wall parts proves the claim.
\end{proof}

\subsection{Random-block concentration}

\begin{lemma}\label{lem:bernstein}
Let $R\subseteq[d]$ have cardinality $n\ge7d/8$, and let $J$ be a
uniformly random $L$-subset of $R$, where $L=\left\lfloor\frac{d}{8M}\right\rfloor$.
Let $(r_j)_{j\in J}$ be independent Rademacher signs and define
\[z=L^{-1/q}\sum_{j\in J}r_je_j,
\qquad h=L^{-1/p}\sum_{j\in J}r_je_j.\]
Then there exists a universal constant $C>0$ such that, for every
$u\ge1$, the following hold.

\begin{enumerate}[label=(\roman*)]
\item
For every deterministic $x\in X$ satisfying $\|x_R\|_p\le4M^{-1/p}$,
with probability at least $1-2e^{-u}$,
\begin{equation}\label{eq:hard-concentration}
|\langle z,x\rangle|
\le C M^{-1/p}\left[M^{\frac12-\frac1p}\sqrt{\frac{u}{d}}+\frac{Mu}{d}\right].
\end{equation}

\item
For every deterministic $s\in B_q^d$ satisfying
$\|s\|_\infty\le L^{-1/q}$, with probability at least $1-2e^{-u}$,
\begin{equation}\label{eq:wall-concentration}
|\langle h,s\rangle| \le
C\left[\sqrt{\frac{u}{d}}
+\frac{Mu}{d}\right].
\end{equation}
\end{enumerate}
\end{lemma}

\begin{proof}
Fix a deterministic vector $a\in\R^R$ and write
\[Q_J(a):=\sum_{j\in J}a_j^2.\]
Since $J$ is a uniform sample without replacement from $R$, Hoeffding's
comparison theorem allows us to compare the moment generating
function of $Q_J(a)$ with that of a with-replacement sample. Let
$Y_1,\ldots,Y_L$ be independent and uniformly distributed over the
finite population $\{a_j^2:j\in R\}$. Then
\[\E Y_\ell=\frac{\|a\|_2^2}{n},
\qquad 0\le Y_\ell\le\|a\|_\infty^2,
\] and
\[\operatorname{Var}(Y_\ell) \le
\frac{\|a\|_\infty^2\|a\|_2^2}{n}.\]
The standard Bernstein bound for the independent sample gives: for every $t>0$, with probability at least $1-e^{-t}$,
\begin{equation}\label{eq:sampled-energy}
Q_J(a) \le 2\left[ \frac{L}{n}\|a\|_2^2 +\|a\|_\infty^2t
\right].
\end{equation}
Conditioning on $J$, since the signs are independent
Rademachers, we have
\[\Pr\left(
\left|\sum_{j\in J}r_ja_j\right|\ge r\,\middle|\,J\right)
\le2\exp\left(-\frac{r^2}{2Q_J(a)}\right)\]
Applying \eqref{eq:sampled-energy} with $t=2u$ and then the Rademacher
tail bound, we obtain, with probability at least $1-2e^{-u}$,
\begin{equation}\label{eq:signed-block}
\left|\sum_{j\in J}r_ja_j\right|
\le C\left[\sqrt{\frac{L}{n}}\|a\|_2\sqrt u + \|a\|_\infty u\right].
\end{equation}
For part~(i), apply \eqref{eq:signed-block} with $a=x_R$ and multiply
by $L^{-1/q}$. Since $p\ge2$,
\[\|x_R\|_2\le n^{\frac12-\frac1p}\|x_R\|_p \le 4n^{\frac12-\frac1p}M^{-1/p}.\]
Hence the first term satisfies
\begin{align*}
L^{-1/q}\sqrt{\frac{L}{n}}\|x_R\|_2\sqrt u
\le C M^{-1/p}
L^{\frac1p-\frac12}n^{-1/p}\sqrt u\le C M^{-1/p}M^{\frac12-\frac1p}
\sqrt{\frac{u}{d}},
\end{align*}
where we used
\[\frac{d}{16M}\le L\le\frac{d}{8M},
\qquad n\ge\frac{7d}{8}.\]
For the second term, because $x\in X$, $\|x_R\|_\infty\le 4d^{-1/p}$,
and therefore
\begin{align*}
L^{-1/q}\|x_R\|_\infty u
\le C\left(\frac{M}{d}\right)^{1/q}d^{-1/p}u=
C\frac{M^{1/q}u}{d}=
C M^{-1/p}\frac{Mu}{d}.
\end{align*}
This proves \eqref{eq:hard-concentration}. For part~(ii), apply \eqref{eq:signed-block} with $a=s_R$ and multiply
by $L^{-1/p}$. Since $q\le2$,
\begin{align*}
\|s_R\|_2^2=\sum_{j\in R}|s_j|^{2-q}|s_j|^q\le\|s\|_\infty^{2-q}\|s\|_q^q\le L^{-(2-q)/q},
\end{align*}
so $\|s_R\|_2 \le L^{\frac1p-\frac12}$.
Thus $L^{-1/p}\sqrt{\frac{L}{n}}\|s_R\|_2\sqrt u \le
C\sqrt{\frac{u}{d}}$,
while
\[L^{-1/p}\|s_R\|_\infty u
\le L^{-1/p}L^{-1/q}u =
\frac{u}{L} \le C\frac{Mu}{d}.\]
which proves \eqref{eq:wall-concentration}.
\end{proof}

\subsection{The high-probability event}

\begin{lemma}\label{lem:all-hidden}
With the choices
\begin{equation}\label{eq:tau}
\tau=C M^{-1/p}\left[M^{\frac12-\frac1p}\sqrt{\frac{\log(16KM^2/\rho)}{d}}
+\frac{M\log(16KM^2/\rho)}
{d}\right]
\end{equation}
and
\begin{equation}\label{eq:theta}
\theta=C\left[\sqrt{\frac{\log(16KM^2/\rho)}{d}}+\frac{M\log(16KM^2/\rho)}{d}
\right],
\end{equation}
and with $C$ sufficiently large, the following statements hold
simultaneously with probability at least $1-\rho$:
\begin{enumerate}[label=(\roman*)]
\item
for every round $t<M$, every query $x$ in that round satisfying
$\|x_{R_t}\|_p\le4M^{-1/p}$, and every $i\ge t$,
\[|\langle z_i,x\rangle|\le\tau;\]
\item
for every round $t<M$, every query $x$ in that round, and every
$i\ge t$, \[|\langle h_i,s_t(x)\rangle|\le\theta/2;\]
\item
after round $M-1$, the final output satisfies the analogous two
bounds with the fresh last block $J_M$.
\end{enumerate}
\end{lemma}

\begin{proof}
For each $t\in[M]$, let $\mathcal F_{t-1}$ denote the $\sigma$-field
generated by all randomness revealed before round $t$. In particular,
$J_1,\ldots,J_{t-1}$ and their signs are
$\mathcal F_{t-1}$-measurable, and therefore
\[R_t=[d]\setminus\bigcup_{j<t}J_j
\]
is fixed under this conditioning. Since
\[|R_t|=d-(t-1)L\ge d-ML \ge\frac{7d}{8},\]
Lemma~\ref{lem:bernstein} applies on $R_t$. Let
\[
x_{t,1},\ldots,x_{t,K_t},
\qquad
K_t\le K,
\]
be the queries in round $t$. Conditioned on $\mathcal F_{t-1}$, these
queries are fixed. Moreover, every query belongs to $X$. For every
$i\ge t$, the conditional marginal distribution of $J_i$ is uniform
over the $L$-subsets of $R_t$, and its signs are fresh. Thus
Lemma~\ref{lem:bernstein} applies to each fixed triple $(t,k,i)$.
Set
\[u:=\log\left(\frac{16KM^2}{\rho}\right).\]
If $\|(x_{t,k})_{R_t}\|_p\le4M^{-1/p}$,
part~(i) of Lemma~\ref{lem:bernstein}, together with the choice of
$\tau$, gives
\[\Pr\left(|\langle z_i,x_{t,k}\rangle|>\tau
\,\middle|\,\mathcal F_{t-1}
\right)\le 4e^{-u}=\frac{\rho}{4KM^2}.\] If the residual norm is larger, the corresponding hard-direction
event is irrelevant.
For the wall bounds, conditioned on $\mathcal F_{t-1}$, the set
$\mathcal C_t$ and the query $x_{t,k}$ are fixed. By
Lemma~\ref{lem:wall-consistency}, $s_t(x_{t,k})$ is uniquely
determined, and
\[
\|s_t(x_{t,k})\|_q\le1,
\qquad
\|s_t(x_{t,k})\|_\infty\le L^{-1/q}.
\]
Part~(ii) of Lemma~\ref{lem:bernstein}, together with the choice of
$\theta$, gives
\[\Pr\left(|\langle h_i,s_t(x_{t,k})\rangle|>\theta/2
\,\middle|\,\mathcal F_{t-1}
\right)\le\frac{\rho}{4KM^2}.\]
The same argument applies to the final output $\widehat x$ after
round $M-1$. By assumption $\widehat x\in X$, and conditioned on
$\mathcal F_{M-1}$ the block $J_M$ and its signs are still fresh.
Therefore the analogous hard-direction and wall bounds each fail
with probability at most $\rho/(4KM^2)$.
There are fewer than $KM^2$ hard-direction events and fewer than
$KM^2$ wall events, together with two final-output events. Since
$M\ge2$ and $K\ge1$,
\[2KM^2+2\le 4KM^2.\]
A union bound therefore gives total failure probability at most $\rho$. 
\end{proof}

\subsection{Proof of Theorem~\ref{thm:main}}

\begin{proof}[Proof of Theorem~\ref{thm:main}]
First consider a deterministic algorithm. The function $f$ is convex and $1$-Lipschitz with respect to $\ell_p$: every affine hard piece has dual norm $1/32$, while Lemma~\ref{lem:wall-consistency} gives that the wall is $1$-Lipschitz.
On the event in Lemma~\ref{lem:all-hidden},
Lemma~\ref{lem:hiding} shows that $f$ and $f_t$ agree on a neighborhood of every round-$t$ query.
Induction over the rounds therefore shows that the complete transcripts coincide through round $M-1$.
Let \[u=\log\left(\frac{16KM^2}{\rho}\right).\]
The assumption \eqref{eq:main-variance}, with $c$ sufficiently small relative to the constant in \eqref{eq:tau}, implies
\begin{equation}\label{eq:smallness}
M\tau\le\frac1{64}M^{-1/p}.
\end{equation} Thus \eqref{eq:smallness} follows by choosing $c$ sufficiently small.
The same assumption also implies $M\theta\le1$. 
We next bound the optimum. By Lemma~\ref{lem:shielding},
\[W(x^*)\le-\frac1{32}M^{-1/p}.\]
Also, for every $i$,
\[\frac1{32}\langle z_i,x^*\rangle
-\frac{i\tau}{8}=-\frac1{32}M^{-1/p}
-\frac{i\tau}{8}\le -\frac1{32}M^{-1/p}.\]
Since $x^*\in X$,
\begin{equation}\label{eq:fstar}
f^*:=\min_{x\in X}f(x)\le
f(x^*)\le-\frac1{32}M^{-1/p}.
\end{equation}
Consider any query $x$ made in round $t<M$. If  
$\|x_{R_t}\|_p>4M^{-1/p}$,
then Lemmas~\ref{lem:hiding} and \ref{lem:shielding} give
\[f(x)=f_t(x)\ge W_t(x)\ge \|x_{R_t}\|_p-\frac{17}{16}M^{-1/p}
>\frac{47}{16}M^{-1/p}.\]
If instead $\|x_{R_t}\|_p\le4M^{-1/p}$,
Lemma~\ref{lem:all-hidden} gives
$|\langle z_t,x\rangle|\le\tau$, and therefore
\[f(x)=f_t(x)\ge\frac1{32}\langle z_t,x\rangle-\frac{t\tau}{8}\ge
-\frac{\tau}{32}-\frac{M\tau}{8}.\]
Combining this with \eqref{eq:smallness} and \eqref{eq:fstar}, in
either case,
\begin{equation}\label{eq:query-gap}
f(x)-f^*\ge\frac1{64}M^{-1/p}.
\end{equation}
It remains to treat the final output $\widehat x$ after round $M-1$.
Since $\widehat x$ is determined by the first $M-1$ blocks, $J_M$ is still uniform over the $L$-subsets of $R_M$, with independent
signs. We now apply part~(iii) of Lemma~\ref{lem:all-hidden}.
If $\|\widehat x_{R_M}\|_p>4M^{-1/p}$,
then part~(iii) of Lemma~\ref{lem:all-hidden} gives
\[|\langle h_M,s_M(\widehat x)\rangle|\le\theta/2.\]
Hence Lemma~\ref{lem:wall-consistency} gives
$W(\widehat x)=W_M(\widehat x)$,
and Lemma~\ref{lem:shielding} yields
\[f(\widehat x)\ge W(\widehat x)
=W_M(\widehat x)>\frac{47}{16}M^{-1/p}.\]
If instead $\|\widehat x_{R_M}\|_p\le4M^{-1/p}$,
part~(iii) gives $|\langle z_M,\widehat x\rangle|
\le\tau$,
so the $M$-th affine piece gives
\[f(\widehat x)\ge-\frac{\tau}{32}
-\frac{M\tau}{8}.\]
Thus \eqref{eq:query-gap} also holds for the final output.
Choose the constant $C$ in \eqref{eq:main-accuracy} sufficiently large that
\[\frac1{64}M^{-1/p}>\eps.\]
Then, on the event of Lemma~\ref{lem:all-hidden}, no output produced in fewer than $M$ rounds is an $\eps$-minimizer. This event has probability at least $1-\rho$. Finally, let the algorithm be randomized, and let $\omega$ denote its random
choices, which are independent of the blocks and signs. For every fixed $\omega$
the algorithm is deterministic, so the argument above shows that it returns an
$\eps$-minimizer in fewer than $M$ rounds with probability at most $\rho$ over
the blocks and signs. Averaging over $\omega$, the same bound holds over
$\omega$, the blocks and the signs jointly. Hence, for some fixed choice of
blocks and signs, the probability over $\omega$ alone is at most $\rho$, and the
corresponding $f$ is the function in the theorem.The condition \eqref{eq:main-accuracy} gives
\[M\lesssim\eps^{-p},\]
while \eqref{eq:main-variance} gives, up to logarithmic factors,
\[M \lesssim d^{\frac{p}{3p-2}}.\]
This proves the theorem.
\end{proof}

\section{Conclusion}
In this work, we studied parallel SCO in $\ell_p$ geometries for $p \in [2, \infty)$. We first introduced a class of \odc algorithms, and presented various applications of these tools. We then showed how to use these algorithms to develop two efficient $\ell_p$ parallel algorithms, which in favorable geometries, improves upon the corresponding sequential stochastic query depth complexity. We finally present a dimension-dependent lower bound on the number of adaptive rounds required for parallel $\ell_p$ nonsmooth convex optimization, improving upon that of \cite{diakonikolas2019lower}. One interesting extension to this work would be to consider how these tools could be used in a private SCO setting, where the goal is to produce an iterate that has at most $\eps_{\mathrm{opt}}$ expected population risk more than the optimal risk. One work that has studied private SCO in $\ell_p$ geometries is \cite{bassily2021non}, and it would be interesting to improve upon the upper and lower bounds presented here. 
\subsection{AI Methodology Statement.} We used ChatGPT 5.6 and Claude 
to assist with checking and aiding some steps of our proofs. 
In particular, the addition of the quadratic term in the wall function was suggested by ChatGPT. Generative AI tools were not used for devising the main ideas and conclusions of the paper. The authors have reviewed all AI-assisted work and take full responsibility for the contents of the paper. 

\newpage 
\bibliographystyle{alpha}
\bibliography{references}
\newpage

\appendix
\section{Helper Facts}
We present a few useful properties and facts about Bregman divergences, which we use in our analysis. 
\begin{property}
\label{property:threepointproperty}
(Three Point Property). For all $a, b, c \in \mathbb{R}^d$ and convex $\Phi$, it holds that \begin{flalign*}
    &\langle \nabla \Phi(a) - \nabla \Phi(b), c - a \rangle = D_\Phi(c, b) - D_\Phi(c, a) - D_\Phi(a, b).
\end{flalign*}
\end{property}

\begin{fact}
    \label{fact:bregmanlowerbound}
    Let $\Phi(x) = \frac{1}{p} ||x||_p^p$ for $p \ge 2$. For all $x, y \in \mathbb{R}^d$, it holds that \begin{flalign*}
        &D_\Phi(x, y) \ge \frac{2^{2 - p}}{p} ||x - y||_p^p.
    \end{flalign*}
\end{fact}

\begin{proof}
    Fix $x, y \in \mathbb{R}^d$. It holds that \begin{flalign*}
        &D_\Phi(x, y) = \int_{0}^1 \langle \nabla \Phi(y + t(x - y)) - \nabla \Phi(y), x - y \rangle dt 
    \end{flalign*} Since $\Phi(x) = \frac{1}{p} \sum_{i = 1}^d |x_i|^p$, we analyze (slightly overloading notation here) for $\psi(t) = |t|^{p - 2} t$: \begin{flalign*}
        &I_\psi(x, y) = \int_0^1 (\psi(y + t(x - y)) - \psi(y))(x - y) dt
    \end{flalign*} We now show that for all $a, b \in \mathbb{R}$, it holds that $(\psi(a) - \psi(b))(a - b) \ge 2^{2 - p}|a - b|^p$. If $a = b$, the result trivially holds, so we only consider the setting in which $a \ne b$. Since $\psi'(t) = (p - 1)|t|^{p - 2}$, it holds that \begin{flalign*}
        &(\psi(a) - \psi(b))(a - b) = (a - b)^2 (p - 1) \int_{0}^1 |b + s(a - b)|^{p - 2} ds
    \end{flalign*} Letting $u = b/(a - b)$ so $b + (a - b)s = (a- b)(u + s)$, we write \begin{flalign*}
        &\int_{0}^1 |b + s(a - b)|^{p - 2} ds = |a - b|^{p - 2} \int_{0}^1 |u + s|^{p - 2} ds
    \end{flalign*} Letting $F(u) = \int_0^1 |u + s|^{p - 2} ds$, observe that the minimum of $F$ occurs at $u = -\frac{1}{2}$, since $|x|^{p - 2}$ is convex for $p \ge 2$. Therefore, \begin{flalign*}
        &|a - b|^{p - 2} \int_{0}^1 |u + s|^{p - 2} \ge |a - b|^{p - 2} \int_{0}^1 |s-\frac{1}{2}|^{p - 2} = |a - b|^{p - 2} \int_{-1/2}^{1/2} |u|^{p - 2} du = |a - b|^{p - 2}\cdot \frac{2^{2 - p}}{p - 1}
    \end{flalign*} which implies that $(\psi(a) - \psi(b))(a - b) \ge 2^{2 - p}|a - b|^p$. Further, it holds that \begin{flalign*}
        &(\psi(y + t(x - y)) - \psi(y))(x - y) = \frac{1}{t} \cdot (\psi(y + t(x - y)) - \psi(y))(tx - ty) \ge 2^{2 - p} |x - y|^p t^{p - 1}
    \end{flalign*} which means that $I_\psi(x, y) \ge 2^{2 - p}|x - y|^p\int_{0}^1 t^{p - 1} dt = \frac{2^{2 - p}}{p} |x - y|^p$. As $D_\Phi = \sum_{i = 1}^d I_\psi(x_i, y_i)$, the claim follows. 
\end{proof}

\begin{fact}
    \label{fact:bregmanupperbound}
    Let $\Phi(x) = \frac{1}{p} ||x||_p^p$ for some $p \ge 2$. For all $x, y \in \mathbb{R}^d$, it holds that \begin{flalign*}
        &D_\Phi(x, y) \le \frac{p - 1}{2} (||x||_p + 2||y||_p)^{p - 2} ||x - y||_p^2.
    \end{flalign*}
\end{fact}

\begin{proof}
    We first prove a smoothness-like guarantee for $\Phi$: $||\nabla \Phi(x) - \nabla \Phi(y)||_q \le (p - 1)(||x||_p + ||y||_p)^{p - 2} ||x - y||_p$. Since $\nabla \Phi(x)_i = \psi(x_i), \psi(t) = |t|^{p - 2}t$, we have that for all $a, b \in \mathbb{R}$:  \begin{flalign*}
        |\psi(a) - \psi(b)| &= (p - 1)|a - b| \int_0^1 |b + s(a - b)|^{p - 2} ds \\ 
        &\le (p - 1)|a - b| \int_0^1 ((1 - s)|b| + s|a|)^{p - 2} ds \\ 
        &\le (p - 1)|a - b| \int_0^1 (|a| + |b|)^{p - 2} ds = (p - 1)|a - b|(|a| + |b|)^{p - 2}
    \end{flalign*} Thus, for each coordinate: $|\nabla \Phi(x_i) - \nabla \Phi(y_i)| \le (p - 1)(|x_i| + |y_i|)^{p - 2}|x_i - y_i|$ so we have that \begin{flalign*}
        &||\nabla \Phi(x) - \nabla \Phi(y)||_q = (\sum_i |\nabla \Phi(x_i) - \nabla \Phi(y_i)|^q)^{1/q} \le (p - 1)(\sum_i (|x_i| + |y_i|)^{(p - 2)q} |x_i - y_i|^q)^{1/q}
    \end{flalign*} Applying Hölder's inequality, we have that \begin{flalign*}
        \sum_i (|x_i| + |y_i|)^{(p - 2)q} |x_i - y_i|^q &\le (\sum_i (|x_i| + |y_i|)^{\frac{(p - 2)q \cdot (p - 1)}{p - 2}})^{\frac{p - 2}{p - 1}} \cdot  (\sum_i |x_i - y_i|^{q(p - 1)})^{\frac{1}{p - 1}} \\ 
        &\le (\sum_i (|x_i| + |y_i|)^p)^{\frac{p - 2}{p - 1}} \cdot (\sum_i|x_i - y_i|^p)^{\frac{1}{p - 1}} 
    \end{flalign*} which implies that \begin{flalign*}
        (\sum_i (|x_i| + |y_i|)^{(p - 2)q} |x_i - y_i|^q)^{1/q} &\le (\sum_i (|x_i| + |y_i|)^p)^{\frac{p - 2}{q(p - 1)}} \cdot (\sum_i|x_i - y_i|^p)^{\frac{1}{q(p - 1)}} \\ 
        &= (\sum_i (|x_i| + |y_i|)^p)^{\frac{p - 2}{p}} \cdot (\sum_i |x_i - y_i|^p)^{\frac{1}{p}} \\ 
        &= ((\sum_i(|x_i| + |y_i|)^p)^{1/p})^{p - 2} \cdot ||x - y||_p \\ 
        &\le ((\sum_i |x_i|^p)^{1/p} + (\sum_i |y_i|^p)^{1/p})^{p - 2} \cdot ||x - y||_p \\ 
        &= (||x||_p + ||y||_p)^{p - 2} \cdot ||x - y||_p
    \end{flalign*} where we used Minkowski's inequality in the last inequality. Therefore, it holds that \begin{flalign*}
        D_\Phi(x, y) &= \int_0^1 \langle \nabla \Phi(y + t(x - y)) - \nabla \Phi(y), x - y \rangle dt \\ 
        &\le \int_0^1 ||\nabla \Phi(y + t(x - y)) - \nabla \Phi(y)||_q \cdot ||x - y||_p dt \\
        &\le \int_0^1(p - 1)\cdot (||y + t(x - y)||_p + ||y||_p)^{p - 2} \cdot ||t(x - y)||_p \cdot ||x - y||_p dt \\
        &= (p - 1)||x - y||_p^2 \int_0^1 t (||y + t(x - y)||_p + ||y||_p)^{p - 2}dt \\ 
        &\le (p - 1)||x - y||_p^2 \cdot (||x||_p + 2||y||_p)^{p - 2} \int_0^1 t dt \\ 
        &= \frac{p - 1}{2} (||x||_p + 2||y||_p)^{p - 2} \cdot ||x - y||_p^2
    \end{flalign*} which finishes the proof. 
\end{proof}
Now, we present a few useful facts regarding our Generalized Gaussian ($\mathcal{GG}_p$) distribution. 

\begin{fact}
\label{fact:convolved-objective-smoothness}
\label{fact:GG_bounded_distance_convolution_f}
Suppose that $f:\mathbb{R}^d\to\mathbb{R}$ is convex and
$L$-Lipschitz with respect to $\|\cdot\|_p$.
Then its generalized Gaussian convolution $\hat f_{\rho,p}$
(Definition~\ref{def:convolution}) satisfies the following:
\begin{enumerate}
    \item $\hat f_{\rho,p}$ is $(C'_pL/\rho)$-smooth with respect to
    $\|\cdot\|_p$, for some constant $C'_p>0$ depending only on $p$.
    \item For every $x\in\mathbb{R}^d$,
    \[
        f(x)\le\hat f_{\rho,p}(x)\le f(x)+L\rho\sqrt d.
    \]
    The dependence on $d$ in this bound is asymptotically tight.
\end{enumerate}
\end{fact}

\begin{proof}
First consider $x,x'\in\mathbb{R}^d$ with
$\|x-x'\|_p\le\rho/2$.
Use the same center $\overline{x}=x$ in both estimators from
Definition~\ref{defn:estimator definition}.
By the unbiasedness assertion in
Lemma~\ref{lem:GREG-q moment}, for
$\xi\sim\mathcal{GG}_p(0,\rho^2)$,
\[
\begin{aligned}
&\nabla\hat f_{\rho,p}(x)-\nabla\hat f_{\rho,p}(x')=
\mathbb{E}_\xi\left[
    \frac{
        \gamma_{\rho,p}(-\xi)
        -\gamma_{\rho,p}(x'-x-\xi)
    }{\gamma_{\rho,p}(\xi)}
    g(x+\xi)
\right].
\end{aligned}
\]
The triangle inequality, the bound $\|g(x+\xi)\|_q\le L$,
and Cauchy--Schwarz give
\[
\begin{aligned}
\|\nabla\hat f_{\rho,p}(x)-\nabla\hat f_{\rho,p}(x')\|_q
&\le
L\,\mathbb{E}_\xi\left|
    \frac{
        \gamma_{\rho,p}(-\xi)
        -\gamma_{\rho,p}(x'-x-\xi)
    }{\gamma_{\rho,p}(\xi)}
\right|\\
&\le
L\left(
\mathbb{E}_\xi\left|
    \frac{
        \gamma_{\rho,p}(-\xi)
        -\gamma_{\rho,p}(x'-x-\xi)
    }{\gamma_{\rho,p}(\xi)}
\right|^2
\right)^{1/2}\\
&\le
\frac{\sqrt{2}\,e^{(p-1)/2}C_pL}{\rho}
\|x-x'\|_p.
\end{aligned}\]
The last inequality follows from the second bound in
Lemma~\ref{lemma:GGguaranteespart2} with $m=2$ and
$\overline{x}=x$. For arbitrary $x,x'\in\mathbb{R}^d$, choose an integer $N\ge1$
such that $\|x-x'\|_p/N\le\rho/2$, and define
\[
    x_i=x+\frac{i}{N}(x'-x),
    \qquad i=0,\ldots,N.
\]
Then $x_0=x$, $x_N=x'$, and
$\|x_i-x_{i-1}\|_p=\|x-x'\|_p/N\le\rho/2$.
Applying the local bound to each consecutive pair yields
\[
\begin{aligned}
\|\nabla\hat f_{\rho,p}(x)-\nabla\hat f_{\rho,p}(x')\|_q
&\le
\sum_{i=1}^{N}
\|\nabla\hat f_{\rho,p}(x_i)-\nabla\hat f_{\rho,p}(x_{i-1})\|_q\\
&\le
\frac{\sqrt{2}\,e^{(p-1)/2}C_pL}{\rho}
\sum_{i=1}^{N}\|x_i-x_{i-1}\|_p\\
&=
\frac{\sqrt{2}\,e^{(p-1)/2}C_pL}{\rho}
\|x-x'\|_p.
\end{aligned}
\]
Thus the first assertion holds with
$C'_p=\sqrt{2}\,e^{(p-1)/2}C_p$.

For the second property, we use the randomized smoothing argument
in~\cite{duchi2012randomized}.
Let $\xi\sim\mathcal{GG}_p(0,\rho^2)$.
By symmetry, $\mathbb{E}[\xi]=0$, so convexity and Jensen's
inequality imply
\[
    f(x)=f\bigl(x+\mathbb{E}[\xi]\bigr)
    \le\mathbb{E}[f(x+\xi)]
    =\hat f_{\rho,p}(x).
\]
On the other hand, Lipschitz continuity and Cauchy--Schwarz give
\[\begin{aligned}
\hat f_{\rho,p}(x)\le f(x)+L\,\mathbb{E}\|\xi\|_p\le f(x)+L\bigl(\mathbb{E}\|\xi\|_p^2\bigr)^{1/2}
\end{aligned}\]
Since $\gamma_{\rho,p}(z)$ is proportional to
$\exp(-\|z\|_p^2/(2\rho^2))$, radial integration gives
\[\begin{aligned}
\mathbb{E}\|\xi\|_p^2=
\frac{\int_0^\infty r^{d+1}
\exp\left(-\frac{r^2}{2\rho^2}\right)\,dr}{
\int_0^\infty r^{d-1}
\exp\left(-\frac{r^2}{2\rho^2}\right)\,dr}=
2\rho^2\frac{\Gamma(d/2+1)}{\Gamma(d/2)}=d\rho^2.
\end{aligned}\]
Consequently,
$\hat f_{\rho,p}(x)\le f(x)+L\rho\sqrt d$.
To see this bound is tight, we take $f(x)=L\|x\|_p$.
Notice that this yields
\[\begin{aligned}
\hat f_{\rho,p}(0)-f(0)
=L\,\mathbb{E}\|\xi\|_p=\sqrt{2}L\rho\,
\frac{\Gamma((d+1)/2)}{\Gamma(d/2)}
\sim L\rho\sqrt d
\qquad\text{as }d\to\infty.
\end{aligned}
\]
\end{proof}

\begin{fact}
\label{fact:expectedexpbound}
    Let $p \in [2, \infty)$ and let $\Phi(x) = \frac{1}{2}||x||_p^2$. Furthermore, let $\gamma_{\rho, p}(x) = Z_{\rho, p}^{-1} \exp(-\frac{\Phi(x)}{\rho^2})$ where $Z_{\rho, p}^{-1}$ is the appropriate normalization constant. Then, for every $u\in\mathbb R^d$,
\[\mathbb E_{\xi\sim\gamma_{\rho,p}}
\left[\exp\left(\frac{\langle\nabla\Phi(\xi),u\rangle}{\rho^2}\right)
\right]\leq\exp\left(\frac{p-1}{2\rho^2}\|u\|_p^2\right)\]
\end{fact}

\begin{proof}
    Since $\Phi(x)=\frac12\|x\|_p^2$ is $(p-1)$-smooth with respect to
$\|\cdot\|_p$, for every $x,h\in\mathbb R^d$,
\[\Phi(x+h)\leq\Phi(x)+
\langle\nabla\Phi(x),h\rangle
+\frac{p-1}{2}\|h\|_p^2\] Taking $x=\xi$ and $h=-u$ gives
\[\Phi(\xi-u)\leq\Phi(\xi)-
\langle\nabla\Phi(\xi),u\rangle
+\frac{p-1}{2}\|u\|_p^2,\]
and therefore
\[\langle\nabla\Phi(\xi),u\rangle
\leq\Phi(\xi)-\Phi(\xi-u)+
\frac{p-1}{2}\|u\|_p^2\] Hence
\begin{align*}
\mathbb E_{\xi\sim\gamma_{\rho,p}}
\left[\exp\left(
\frac{\langle\nabla\Phi(\xi),u\rangle}{\rho^2}\right)\right] &=
\frac{1}{Z_{\rho,p}}\int_{\mathbb R^d}\exp\left(\frac{\langle\nabla\Phi(\xi),u\rangle-\Phi(\xi)}{\rho^2}
\right)d\xi \\
&\leq\exp\left(\frac{p-1}{2\rho^2}\|u\|_p^2\right)
\frac{1}{Z_{\rho,p}}\int_{\mathbb R^d}\exp\left(-\frac{\Phi(\xi-u)}{\rho^2}\right)d\xi.
\end{align*} By translation invariance of Lebesgue measure,
\[\int_{\mathbb R^d}\exp\left(
-\frac{\Phi(\xi-u)}{\rho^2}
\right)d\xi=\int_{\mathbb R^d}
\exp\left(-\frac{\Phi(\xi)}{\rho^2}
\right)d\xi=Z_{\rho,p}\]
Therefore,
\[\mathbb E_{\xi\sim\gamma_{\rho,p}}
\left[\exp\left(
\frac{\langle\nabla\Phi(\xi),u\rangle}{\rho^2}\right)\right]\leq
\exp\left(\frac{p-1}{2\rho^2}\|u\|_p^2\right)\]
 \end{proof}

\begin{algorithm}[H]
\caption{Aggregate}\label{alg:aggregate}
\begin{algorithmic}[1]
\State Input: $\mathcal{S} = \{y_i\}_{i \in [k]} \subset \mathbb{R}^d$, $R \ge 0$ such that for some $z \in \mathbb{R}^d$ at least $0.51k$ points $y_i \in \mathcal{S}$ have $||y_i - z||_p \le R/3$
\For{$i \in [k]$}
    \If{at least 0.51k points $y_j \in \mathcal{S}$ satisfy $||y_i - y_j||_p \le 2R/3$} \Return $\tilde{z} \gets y_i$
    \EndIf
\EndFor
\end{algorithmic}
\end{algorithm}

\begin{fact}
\label{fact:aggregateguarantee}
Let $\mathcal S=(y_1,\ldots,y_k)$ with $y_i\in\mathbb R^d$,
and let $R>0$.
Suppose that there exists $z\in\mathbb R^d$ such that
\[\left|\left\{i\in[k]:\|y_i-z\|_p\le R/3\right\}\right|
\ge0.51k.\]
Then $\textsc{Aggregate}(\mathcal S,R)$ runs in $O(k^2d)$ time
and returns $\tilde z\in\mathcal S$ satisfying
$\|\tilde z-z\|_p\le R$.
\end{fact}

\begin{proof}
Following \cite{kelner2023semi}, let
\[
T=\{i\in[k]:\|y_i-z\|_p\le R/3\}.
\]
For every $i,j\in T$, the triangle inequality gives
$\|y_i-y_j\|_p\le2R/3$.
Thus every $i\in T$ has at least $0.51k$ points within
distance $2R/3$, so the algorithm has a valid point to return.

Suppose that the algorithm returns $\tilde z=y_i$, and let
\[
\mathcal Q=\{j\in[k]:\|y_i-y_j\|_p\le2R/3\}.
\]
By the selection rule, $|\mathcal Q|\ge0.51k$.
Since $|T|\ge0.51k$,
\[
|\mathcal Q\cap T|
\ge|\mathcal Q|+|T|-k
\ge0.02k>0.
\]
Choose $j^\star\in\mathcal Q\cap T$.
Then
\[\|\tilde z-z\|_p
\le\|y_i-y_{j^\star}\|_p
+\|y_{j^\star}-z\|_p
\le\frac{2R}{3}+\frac{R}{3}
=R.\]
The running time is $O(k^2d)$ because the algorithm performs
at most $O(k^2)$ pairwise distance computations,
each costing $O(d)$ time.
\end{proof}
\end{document}